\documentclass[10pt,leqno,a4paper]{article}

\usepackage{texdraw,multicol,rotating}
\usepackage{amsmath}
\usepackage{epsfig}
\usepackage{texdraw}
\usepackage{amssymb,amsthm,amsfonts,amscd}
\input{txdtools.tex}
\theoremstyle{plain}
\newtheorem{thm}{Theorem}[section]

\theoremstyle{definition}
\newtheorem{defi}[thm]{Definition}

\newtheorem{example}[thm]{Example}
\theoremstyle{remark}
\newtheorem{remark}[thm]{Remark}
\numberwithin{equation}{section}

\newbox{\tmpa}
\newbox{\tmpb}

\newcommand{\nc}{\newcommand}
\nc{\Uq}{U_q} \nc{\Z}{\mathbf{Z}} \nc{\C}{\mathbf{C}}
\nc{\Q}{\mathbf{Q}}
\nc{\op}{\oplus} \nc{\ot}{\otimes} \nc{\pv}{P^{\vee}}
\nc{\ali}{\alpha_i} \nc{\B}{\mathbf{B}} \nc{\F}{\mathbf{F}}
\nc{\bP}{\mathbf{P}} \nc{\V}{\mathbf{V}} \nc{\La}{\Lambda}
\nc{\la}{\lambda}
\nc{\nbinom}[2]{\genfrac{}{}{0pt}{1}{{#1}}{{#2}}}
\nc{\qbinom}[2]{\left[\genfrac{}{}{0pt}{1}{{#1}}{{#2}}\right]}
\nc{\path}{\mathcal{P}} \nc{\fit}{\tilde{f}_i}
\nc{\eit}{\tilde{e}_i} \nc{\Y}{\mathbf{Y}} \nc{\A}{\mathbf{A}}
\nc{\ra}{\rightarrow} \nc{\vep}{\varepsilon} \nc{\vphi}{\varphi}
\nc{\g}{\mathfrak{g}} \nc{\h}{\mathfrak{h}} \nc{\oP}{\overline{P}}
\nc{\pathp}{\mathbf{p}}

\nc{\tris}{ \bsegment \move(0 0)\lvec(10 0)\lvec(10 10)\lvec(0
0)\ifill f:0.7 \esegment } \nc{\recs}{ \bsegment \move(0
0)\lvec(10 0)\lvec(10 5)\lvec(0 5)\lvec(0 0)\ifill f:0.7 \esegment
}
\nc{\hcvec}[5]{%
\getpos(#1 #3)\spx\spy \getpos(#2 #3)\epx\epy \getpos(#4
#5)\xoff\yoff \realadd \spx \xoff \twox \realadd \epx {-\xoff}
\thrx \realadd \spy \yoff \posy \move({\spx} {\spy}) \clvec
({\twox} {\posy})({\thrx} {\posy})({\epx} {\epy}) \rmove(0 0) }

\nc{\ahead}[2]{%
\cossin (0 0)({#1} {#2})\cosa\sina \bsegment
  \drawdim in \setunitscale 0.065
  \realmult {-0.5} \cosa \hcosa
  \realmult {-0.5} \sina \hsina
  \move({\hcosa} {\hsina}) \ravec({\cosa} {\sina})
\esegment }
\nc{\boxi}{%
{%
\savebox{\tmppic}{\begin{texdraw} \small \drawdim em \textref h:C
v:C \setunitscale 0.55 \htext(0 0){$i$} \move(-1 -1)\lvec(-1
1)\lvec(1 1)\lvec(1 -1)\lvec(-1 -1)
\end{texdraw}}%
\raisebox{-0.19\height}{\usebox{\tmppic}}%
}%
}
\nc{\boxj}{%
{%
\savebox{\tmppic}{\begin{texdraw} \small \drawdim em \textref h:C
v:C \setunitscale 0.55 \htext(0 0.1){$j$} \move(-1 -1)\lvec(-1
1)\lvec(1 1)\lvec(1 -1)\lvec(-1 -1)
\end{texdraw}}%
\raisebox{-0.19\height}{\usebox{\tmppic}}%
}%
}
\nc{\boxipo}{%
{%
\savebox{\tmppic}{\begin{texdraw} \small \drawdim em \textref h:C
v:C \setunitscale 0.55 \htext(0.15 0){$i\!\!+\!\!1$} \move(-1.4
-1)\lvec(-1.4 1)\lvec(1.4 1)\lvec(1.4 -1)\lvec(-1.4 -1)
\end{texdraw}}%
\raisebox{-0.19\height}{\usebox{\tmppic}}%
}%
} \everytexdraw{ \drawdim in \arrowheadsize l:0.065 w:0.03
\arrowheadtype t:F \textref h:C v:C }
\newsavebox{\tmppic}
\newsavebox{\tmpfig}
\newsavebox{\tmpdraw}
\newsavebox{\tmpfiga}
\newsavebox{\tmpfigb}
\newsavebox{\tmpfigc}
\newsavebox{\tmpfigd}
\newsavebox{\tmpfige}
\newsavebox{\tmpfigf}
\newsavebox{\tmpfigg}
\newsavebox{\tmpfigh}
\newsavebox{\tmpfigi}
\newsavebox{\tmpfigj}
\newsavebox{\tmpfigk}
\newsavebox{\tmpfigl}
\newsavebox{\tmpfigm}
\newsavebox{\tmpfign}
\newsavebox{\tmpfigo}
\newsavebox{\tmpfigp}
\newsavebox{\tmpfigq}
\newsavebox{\tmpfigr}
\newsavebox{\tmpfigs}
\newsavebox{\tmpfigt}
\newsavebox{\tmpfigu}
\newsavebox{\tmpfigv}
\newsavebox{\tmpfigw}
\newsavebox{\tmpfigx}
\newsavebox{\tmpfigy}
\newsavebox{\tmpfigz}

\newsavebox{\tmpfigaa}
\newsavebox{\tmpfigab}
\newsavebox{\tmpfigac}
\newsavebox{\tmpfigad}
\newsavebox{\tmpfigae}
\newsavebox{\tmpfigaf}
\newsavebox{\tmpfigag}
\newsavebox{\tmpfigah}
\newsavebox{\tmpfigai}
\newsavebox{\tmpfigaj}
\newsavebox{\tmpfigak}
\newsavebox{\tmpfigal}
\newsavebox{\tmpfigam}
\newsavebox{\tmpfigan}
\newsavebox{\tmpfigao}
\newsavebox{\tmpfigap}
\newsavebox{\tmpfigaq}
\newsavebox{\tmpfigar}
\newsavebox{\tmpfigas}
\newsavebox{\tmpfigat}
\newsavebox{\tmpfigau}
\newsavebox{\tmpfigav}
\newsavebox{\tmpfigaw}
\newsavebox{\tmpfigax}
\newsavebox{\tmpfigay}
\newsavebox{\tmpfigaz}

\newsavebox{\tmpfigba}
\newsavebox{\tmpfigbb}
\newsavebox{\tmpfigbc}
\newsavebox{\tmpfigbd}
\newsavebox{\tmpfigbe}
\newsavebox{\tmpfigbf}
\newsavebox{\tmpfigbg}
\newsavebox{\tmpfigbh}

\nc{\node}{\lcir r:1 }
\nc{\sline}{\bsegment\savepos(10 0)(*ex *ey)
            \move(1 0)\rlvec(8 0)
            \esegment\move(*ex *ey)}
\nc{\dline}{\bsegment\savepos(10 0)(*ex *ey)
            \move(0.93 0.4)\rlvec(8.14 0)\rmove(0 -0.8)\rlvec(-8.14 0)
            \esegment\move(*ex *ey)}
\nc{\uline}{\bsegment\savepos(0 10)(*ex *ey)
            \move(0 1)\rlvec(0 8)
            \esegment\move(*ex *ey)}
\nc{\lpoint}{\savecurrpos(*ex *ey)
             \rmove(2.5 1.5)\rlvec(-1.5 -1.5)\rlvec(1.5 -1.5)
             \move(*ex *ey)}
\nc{\rpoint}{\savecurrpos(*ex *ey)
             \rmove(-2.5 -1.5)\rlvec(1.5 1.5)\rlvec(-1.5 1.5)
             \move(*ex *ey)}
\nc{\bline}{\bsegment\savepos(10 0)(*ex *ey)
            \linewd 0.6 \move(1.1 0)\rlvec(7.8 0)
            \esegment\move(*ex *ey)}

\nc{\araise}[1]{\raisebox{4.5pt}{#1}}
\nc{\braise}[1]{\raisebox{12.1pt}{#1}}
\nc{\craise}[1]{\raisebox{8pt}{#1}}
\nc{\draise}[1]{\raisebox{12pt}{#1}}

\begin{document}

\title{\bf Young Wall Realization of Crystal Bases for Classical Lie Algebras}

\author{Seok-Jin Kang$^\star$
\begin{thanks}
{This research was supported by KOSEF Grant \# 98-0701-01-5-L and
the Young Scientist Award, Korean Academy of Science and
Technology}
\end{thanks},
Jeong-Ah Kim$^\diamond$
\begin{thanks}
{This research was supported by KOSEF Grant \# 98-0701-01-5-L and
BK21 Mathematical Sciences Division, Seoul National University}
\end{thanks},
Hyeonmi Lee$^\star$ $^\dag$ and Dong-Uy Shin$^\star$ $^\dag$ \cr {
} \cr {\small $^{\star}$ School of Mathematics, Korea Institute
for Advanced Study,} \cr {\small 207-43 Cheongryangri-dong,
Dongdaemun-gu,} \cr {\small Seoul 130-012, Korea} \cr {\small
$^{\diamond}$ Department of Mathematics, Seoul National
University,} \cr {\small Seoul 151-747, Korea} \cr {\small
\texttt{sjkang@kias.re.kr,}}
{\small\texttt{jakim@math.snu.ac.kr,}}\cr {\small
\texttt{hmlee@kias.re.kr,}} {\small\texttt{shindong@kias.re.kr}} }

\date{}
\maketitle

\begin{abstract}
In this paper, we give a new realization of crystal bases for
finite dimensional irreducible modules over classical Lie
algebras. The basis vectors are parameterized by certain Young
walls lying between highest weight and lowest weight vectors.
\end{abstract}
\maketitle

\vskip 1cm

\section*{Introduction}

The classical Lie algebras and their representations have been
the fundamental algebraic structure behind many branches of
mathematics and mathematical physics. Through the past 100 years,
it has been discovered that the representation theory of
classical Lie algebras has a close connection with the
combinatorics of Young tableaux and symmetric functions. (see,
for example, \cite{Fulton}, \cite{Mac}.) As can be found in
\cite{KasNak, Nak}, this connection can be explained in a
beautiful manner using the {\it crystal basis theory} for quantum
groups, and one can derive a lot of new and interesting results
in combinatorial representation theory.

\vskip 3mm The {\it quantum groups} are deformations of the
universal enveloping algebras of Kac-Moody algebras, and the
crystal bases can be viewed as bases at $q=0$ for the integrable
modules over quantum groups in the category ${\mathcal O}_{int}$.
The crystal bases are given a structure of colored oriented
graphs, called the {\it crystal graphs}, which reflect the
combinatorial structure of integrable modules in the category
${\mathcal O}_{int}$. Moreover, they have many nice combinatorial
features; for instance, they have a remarkably simple behavior
with respect to taking the tensor product.

\vskip 2mm For classical Lie algebras, Kashiwara and Nakashima
gave an explicit realization of crystal bases for finite
dimensional irreducible modules \cite{KasNak}. In their work,
crystal bases were characterized as the sets of semistandard Young
tableaux with given shapes satisfying certain additional
conditions. Motivated by their work, Kang and Misra discovered a
Young tableaux realization of crystal bases for finite dimensional
irreducible modules over the exceptional Lie algebra $G_{2}$
\cite{KM}. In \cite{Li1},  Littelmann gave another description of
crystal bases for finite dimensional simple Lie algebras using the
Lakshmibai-Seshadri monomial theory. His approach was generalized
(by himself) to the {\it path model theory} for all symmetrizable
Kac-Moody algebras \cite{Li2,Li3}. Littelmann's theory also gives
rise to colored oriented graphs, which turned out to be isomorphic
to the crystal graphs \cite{Kas96}.

\vskip 3mm In this paper, we give a new realization of crystal
bases for finite dimensional irreducible modules over classical
Lie algebras. The basis vectors are parameterized by certain {\it
Young walls} lying between the highest weight and lowest weight
vectors. The Young walls were introduced in \cite{HK} and
\cite{Kang} as a combinatorial scheme for realizing the crystal
bases for quantum affine algebras. They consist of colored blocks
with various shapes built on the given ground-states and can be
viewed as generalization of Young diagrams. The crystal bases for
basic representations for quantum affine algebras are
characterized as the sets of {\it reduced proper Young walls}
\cite{Kang}.

\vskip 3mm Let us briefly explain the main idea of our approach.
Let ${\mathfrak g}$ be a classical Lie algebra lying inside an
affine Lie algebra ${\widehat {\mathfrak g}}$ so that the Dynkin
diagram of ${\mathfrak g}$ can be obtained by removing the 0-node
from the Dynkin diagram of ${\widehat{\mathfrak g}}$. Consider the
crystal graph $B(\Lambda)$ of a basic representation $V(\Lambda)$
of ${\widehat{\mathfrak g}}$ consisting of reduced proper Young
walls. If we remove all the 0-arrows in $B(\Lambda)$, it is
decomposed into a disjoint union of infinitely many connected
components, each of which is isomorphic to the crystal graph
$B(\lambda)$ for a finite dimensional irreducible $\mathfrak
g$-module $V(\lambda)$ with highest weight $\lambda$. Conversely,
any crystal graph $B(\lambda)$ for a finite dimensional
irreducible $\mathfrak g$-module $V(\lambda)$ arises in this way.
That is, given a dominant integral weight $\lambda$ for $\frak g$,
there is a dominant integral weight $\Lambda$ of level $1$ for
$\widehat{\frak g}$ such that $B(\lambda)$ appears as a connected
component in $B(\Lambda)$ without $0$-arrows.

\vskip 2mm Thus the remaining task is to characterize these
connected components in $B(\Lambda)$. However, given a dominant
integral weight $\lambda$ for the classical Lie algebra $\mathfrak
g$, there are infinitely many connected components in $B(\Lambda)$
that are isomorphic to $B(\lambda)$. Among these, we choose the
characterization of $B(\lambda)$ corresponding to the connected
components having the least number of blocks.

\vskip 3mm  In \cite{KS}, from this Young wall realization of
crystal bases over classical Lie algebras using affine
combinatorial objects, Kim and Shin derived another tableaux
realization, which is different from the one given by Kashiwara
and Nakashima.  Moreover, using the result of our work, Lee gave a
realization of $A_n$ type Demazure crystals for certain highest
weights \cite{Lee}.

\vskip 5mm \noindent {\bf Acknowledgments.} \ We would like to
express our sincere gratitude to Professor Georgia Benkart for her
interest in this work and many valuable discussions.

\vskip 1cm

\section{Quantum groups and Young walls}
The basic notions on quantum groups and crystal bases may be found
in \cite{HK,Kas90,Kas91}. In this section, we mostly explain the
basic combinatorics of Young walls which were introduced in
\cite{HK,Kang}.

Let us fix basic notations :

$$\aligned
&\mathfrak{g} : \text{Kac-Moody algebra of finite classical type.}\\
&U_q(\mathfrak{g}) : \text{quantum classical algebra.}\\
&\widehat{\mathfrak{g}} : \text{Kac-Moody algebra of affine type.}\\
&U_q(\widehat{\mathfrak{g}}) : \text{quantum affine algebra.}\\
&I : \text{index set for simple roots of finite or affine Kac-Moody algebra.}\\
&\pv=\left\{\begin{array}{ll} \bigoplus_{i \in I} \Z
h_i &\text{for finite type}\\
\bigoplus_{i \in I} \Z h_i \oplus \Z d&\text{for affine type}
\end{array}\right. : \text{dual weight lattice.}
\endaligned
$$

$$\aligned
&\ali, \delta, \La_i : \text{simple root, null root, fundamental
weight.}\\
&P=\{\lambda\in \mathfrak h^* | \lambda(P^{\vee})\subset \Z \} :
      \text{weight lattice.}\\
&\eit, \fit : \text{Kashiwara operators.}
\endaligned
$$

\vskip 3mm The {\it Young walls} are built of colored blocks with
three different shapes:

\vskip 5mm
\hskip 1cm
\begin{center}

\begin{tabular}{rcl}
\raisebox{-0.35em}[1.25em][0.25em]{%
\begin{texdraw}
\drawdim em \setunitscale 0.12 \textref h:C v:C \move(0 0)\lvec(10
0)\lvec(10 10)\lvec(0 10)\lvec(0 0) \move(10 0)\lvec(15 5)\lvec(15
15)\lvec(5 15)\lvec(0 10) \move(10 10)\lvec(15 15)
\end{texdraw}%
} & : &
unit width, unit height, unit thickness,\\

\\

\raisebox{-0.35em}[1.25em][0.25em]{%
\begin{texdraw}
\drawdim em \setunitscale 0.12 \textref h:C v:C \move(0 0)\lvec(10
0)\lvec(10 10)\lvec(0 10)\lvec(0 0) \move(10 0)\lvec(12.5
2.5)\lvec(12.5 12.5)\lvec(2.5 12.5)\lvec(0 10) \move(10
10)\lvec(12.5 12.5)
\end{texdraw}%
} & : &
unit width, unit height, half-unit thickness,\\

\\

\raisebox{-0.35em}[1.25em][0.25em]{%
\begin{texdraw}
\drawdim em \setunitscale 0.12 \textref h:C v:C \move(0 0)\lvec(10
0)\lvec(10 5)\lvec(0 5)\lvec(0 0) \move(10 0)\lvec(15 5)\lvec(15
10)\lvec(5 10)\lvec(0 5) \move(10 5)\lvec(15 10)
\end{texdraw}%
} & : & unit width, half-unit height, unit thickness.
\end{tabular}

\end{center}

\vskip 3mm With these colored blocks, we build the walls of
thickness less than or equal to 1 unit which extend infinitely to
the left. Given a dominant integral weight $\Lambda$ of level 1
for the affine Lie algebra $\widehat{\mathfrak g}$, we fix a frame
called the {\it ground-state wall} of weight $\Lambda$, and build
the walls on this frame. For each type of quantum affine algebras,
we use different sets of colored blocks and ground-state walls,
whose description can be found in \cite{HK,Kang}.

\vskip 3mm The rules for building the walls are given as follows:
\begin{itemize}

\item [(1)] The walls must be built on top of the ground-state wall.

\item [(2)] The colored blocks should be stacked in the patterns given in
\cite{HK,Kang}.

\item [(3)] No block can be placed on top of a column of half-unit thickness.

\item [(4)] Except for the right-most column, there should be no free space
      to the right of any block.

\end{itemize}

By (4), the heights of the columns are weakly decreasing as we go
from right to left. For this reason, the walls built by the above
rules will be called the {\it Young walls}.

\vskip 3mm

In the following example, for the affine Lie algebra
$B_{n}^{(1)}$, we will illustrate the colored blocks, the
ground-state wall, and the pattern for building the walls. For
convenience, we will use the following notations:

\vskip 3mm

\begin{center}
\begin{tabular}{rcl}
\raisebox{-0.4\height}{
\begin{texdraw}
\drawdim em \setunitscale 0.12 \move(-10 0)\lvec(0 0)\lvec(0
10)\lvec(-10 10)\lvec(-10 0) \move(0 0)\lvec(5 5)\lvec(5
15)\lvec(-5 15)\lvec(-10 10) \move(0 10)\lvec(5 15) \htext(-5
5){$*$}
\end{texdraw}
} & $\longleftrightarrow$ & \raisebox{-0.4\height}{
\begin{texdraw}
\drawdim em \setunitscale 0.12 \move(-10 0)\lvec(0 0)\lvec(0
10)\lvec(-10 10)\lvec(-10 0) \htext(-5 5){$*$}
\end{texdraw}
}\\[4mm]
\raisebox{-0.4\height}{
\begin{texdraw}
\drawdim em \setunitscale 0.12 \move(0 0)\lvec(10 0)\lvec(10
5)\lvec(0 5)\lvec(0 0) \move(10 0)\lvec(15 5)\lvec(15 10)\lvec(5
10)\lvec(0 5) \move(10 5)\lvec(15 10) \htext(5 1.5){*}
\end{texdraw}
} & $\longleftrightarrow$ & \raisebox{-0.4\height}{
\begin{texdraw}
\drawdim em \setunitscale 0.12 \textref h:C v:C \move(0 0)\lvec(10
0)\lvec(10 5)\lvec(0 5)\lvec(0 0) \htext(5 1.5){*}
\end{texdraw}
}
\end{tabular}
\qquad
\begin{tabular}{rcl}
\raisebox{-0.4\height}{
\begin{texdraw}
\drawdim em \setunitscale 0.12 \move(-10 0)\lvec(0 0)\lvec(0
10)\lvec(-10 10)\lvec(-10 0) \move(0 0)\lvec(2.5 2.5)\lvec(2.5
12.5)\lvec(-7.5 12.5)\lvec(-10 10) \move(0 10)\lvec(2.5 12.5)
\lpatt(0.3 1) \move(0 0)\lvec(-2.5 -2.5)\lvec(-12.5 -2.5)\lvec(-10
0) \htext(-5 5){$*$}
\end{texdraw}
} & $\longleftrightarrow$ & \raisebox{-0.4\height}{
\begin{texdraw}
\drawdim em \setunitscale 0.12 \move(-10 0)\lvec(0 0)\lvec(0
10)\lvec(-10 10)\lvec(-10 0) \move(0 10)\lvec(-10 0) \htext(-2.5
2.5){$*$}
\end{texdraw}
}\\[4mm]
\raisebox{-0.4\height}{
\begin{texdraw}
\drawdim em \setunitscale 0.12 \move(-10 0)\lvec(0 0)\lvec(0
10)\lvec(-10 10)\lvec(-10 0) \move(0 0)\lvec(2.5 2.5)\lvec(2.5
12.5)\lvec(-7.5 12.5)\lvec(-10 10) \move(0 10)\lvec(2.5 12.5)
\lpatt(0.3 1) \move(2.5 2.5)\lvec(5 5)\lvec(2.5 5) \htext(-5 4){*}
\end{texdraw}
} & $\longleftrightarrow$ & \raisebox{-0.4\height}{
\begin{texdraw}
\drawdim em \setunitscale 0.12 \move(-10 0)\lvec(0 0)\lvec(0
10)\lvec(-10 10)\lvec(-10 0) \move(0 10)\lvec(-10 0) \htext(-6.5
6){*}
\end{texdraw}
}
\end{tabular}
\end{center}

\begin{example}
The walls for the affine Lie algebra $B_{n}^{(1)}$ are built of
the following data.

\vskip 2mm

(a) Colored blocks: \vskip 3mm \hskip 15mm
\raisebox{-0.33\height}[0.69\height][0.35\height]{%
\begin{texdraw}
\fontsize{9}{9}\selectfont \drawdim em \setunitscale 0.14 \move(0
0)\lvec(10 0)\lvec(10 10)\lvec(0 10)\lvec(0 0) \move(10
0)\lvec(12.5 2.5)\lvec(12.5 12.5)\lvec(2.5 12.5)\lvec(0 10)
\move(10 10)\lvec(12.5 12.5) \htext(5 5){$0$}
\end{texdraw}%
}\,,\quad
\raisebox{-0.33\height}[0.69\height][0.35\height]{%
\begin{texdraw}
\fontsize{9}{9}\selectfont \drawdim em \setunitscale 0.14 \move(0
0)\lvec(10 0)\lvec(10 10)\lvec(0 10)\lvec(0 0) \move(10
0)\lvec(12.5 2.5)\lvec(12.5 12.5)\lvec(2.5 12.5)\lvec(0 10)
\move(10 10)\lvec(12.5 12.5) \htext(5 5){$1$}
\end{texdraw}%
}\,,\quad \raisebox{-0.33\height}[0.69\height][0.35\height]{%
\begin{texdraw}
\fontsize{9}{9}\selectfont \drawdim em \setunitscale 0.14 \move(0
0)\lvec(10 0)\lvec(10 5)\lvec(0 5)\lvec(0 0) \move(10 0)\lvec(15
5)\lvec(15 10)\lvec(5 10)\lvec(0 5) \move(10 5)\lvec(15 10)
\htext(5 2.5){$_{n}$}
\end{texdraw}%
} \,,\quad \raisebox{-0.33\height}[0.69\height][0.35\height]{%
\begin{texdraw}
\fontsize{9}{9}\selectfont \drawdim em \setunitscale 0.14 \move(0
0)\lvec(10 0)\lvec(10 10)\lvec(0 10)\lvec(0 0) \move(10 0)\lvec(15
5)\lvec(15 15)\lvec(5 15)\lvec(0 10) \move(10 10)\lvec(15 15)
\htext(5 5){$j$}
\end{texdraw}%
} \quad ($j=2,\cdots,n-1$)

\vskip 5mm

(b) The ground-state wall of weight $\Lambda_{0}$:

\vskip 3mm\hskip 15mm $Y_{\La_0} =
\raisebox{-0.3\height}{\begin{texdraw} \fontsize{9}{9}\selectfont
\drawdim em \setunitscale 0.14 \move(-42 0)\lvec(0 0)\lvec(2.5
2.5)\lvec(2.5 12.5)\lvec(-39.5 12.5) \move(-42 10)\lvec(0 10)
\move(0 0)\lvec(0 10)\lvec(2.5 12.5) \move(-10 0)\lvec(-10
10)\lvec(-7.5 12.5) \move(-20 0)\lvec(-20 10)\lvec(-17.5 12.5)
\move(-30 0)\lvec(-30 10)\lvec(-27.5 12.5) \move(-40 0)\lvec(-40
10)\lvec(-37.5 12.5) \move(0 0)\lvec(-2.5 -2.5)\lvec(-44.5 -2.5)
\move(-10 0)\rlvec(-2.5 -2.5) \move(-20 0)\rlvec(-2.5 -2.5)
\move(-30 0)\rlvec(-2.5 -2.5) \move(-40 0)\rlvec(-2.5 -2.5)
\htext(-5 5){$1$} \htext(-15 5){$0$} \htext(-25 5){$1$} \htext(-35
5){$0$}
\end{texdraw}}=\raisebox{-0.3\height}{\begin{texdraw} \drawdim
em \setunitscale 0.13 \linewd 0.5 \move(-5 0)\lvec(40 0)\lvec(40
10)\lvec(-5 10)\move(0 0)\rlvec(0 10)\move(10 0)\rlvec(0
10)\lvec(0 0)\move(20 0)\rlvec(0 10)\lvec(10 0)\move(30 0)\rlvec(0
10)\lvec(20 0)\move(40 10)\lvec(30 0) \htext(7
3){{\tiny$0$}}\htext(17 3){{\tiny$1$}}\htext(27
3){{\tiny$0$}}\htext(37 3){{\tiny$1$}}
\end{texdraw}}$\,.

\vskip 5mm (c) The pattern for building the walls on
$Y_{\Lambda_{0}}$:

\savebox{\tmpfigk}{\begin{texdraw} \fontsize{7}{7}\selectfont
\drawdim em \setunitscale 1.8 \nc{\dtri}{ \bsegment \move(-1
0)\lvec(0 1)\lvec(0 0)\lvec(-1 0)\ifill f:0.7 \esegment } \move(0
0)\dtri \move(-1 0)\dtri \move(-2 0)\dtri \move(-3 0)\dtri \move(0
0)\rlvec(-4.3 0) \move(0 1)\rlvec(-4.3 0) \move(0 2)\rlvec(-4.3 0)
\move(0 3.5)\rlvec(-4.3 0) \move(0 4.5)\rlvec(-4.3 0) \move(0
5.5)\rlvec(-4.3 0) \move(0 6.5)\rlvec(-4.3 0) \move(0
8)\rlvec(-4.3 0) \move(0 9)\rlvec(-4.3 0) \move(0 10)\rlvec(-4.3
0) \move(0 11)\rlvec(-4.3 0) \move(0 0)\rlvec(0 11.3) \move(-1
0)\rlvec(0 11.3) \move(-2 0)\rlvec(0 11.3) \move(-3 0)\rlvec(0
11.3) \move(-4 0)\rlvec(0 11.3) \move(-1 0)\rlvec(1 1) \move(-2
0)\rlvec(1 1) \move(-3 0)\rlvec(1 1) \move(-4 0)\rlvec(1 1)
\move(-1 9)\rlvec(1 1) \move(-2 9)\rlvec(1 1) \move(-3 9)\rlvec(1
1) \move(-4 9)\rlvec(1 1) \move(0 5)\rlvec(-4.3 0) \htext(-0.3
0.25){$1$} \htext(-0.75 0.75){$0$} \htext(-0.5 1.5){$2$}
\vtext(-0.5 2.75){$\cdots$} \htext(-0.5 4){$n\!\!-\!\!1$}
\htext(-0.5 6){$n\!\!-\!\!1$} \htext(-0.5 8.5){$2$} \htext(-0.3
9.25){$1$} \htext(-0.75 9.75){$0$} \htext(-0.5 10.5){$2$}
\htext(-2.3 0.25){$1$} \htext(-2.75 0.75){$0$} \htext(-2.5
1.5){$2$} \vtext(-2.5 2.75){$\cdots$} \htext(-2.5
4){$n\!\!-\!\!1$} \htext(-2.5 6){$n\!\!-\!\!1$} \htext(-2.5
8.5){$2$} \htext(-2.3 9.25){$1$} \htext(-2.75 9.75){$0$}
\htext(-2.5 10.5){$2$} \htext(-1.3 0.25){$0$} \htext(-1.75
0.75){$1$} \htext(-1.5 1.5){$2$} \vtext(-1.5 2.75){$\cdots$}
\htext(-1.5 4){$n\!\!-\!\!1$} \htext(-1.5 6){$n\!\!-\!\!1$}
\htext(-1.5 8.5){$2$} \htext(-1.3 9.25){$0$} \htext(-1.75
9.75){$1$} \htext(-1.5 10.5){$2$} \htext(-3.3 0.25){$0$}
\htext(-3.75 0.75){$1$} \htext(-3.5 1.5){$2$} \vtext(-3.5
2.75){$\cdots$} \htext(-3.5 4){$n\!\!-\!\!1$} \htext(-3.5
6){$n\!\!-\!\!1$} \htext(-3.5 8.5){$2$} \htext(-3.3 9.25){$0$}
\htext(-3.75 9.75){$1$} \htext(-3.5 10.5){$2$} \htext(-0.5
4.75){$n$} \htext(-2.5 4.75){$n$} \htext(-1.5 5.25){$n$}
\htext(-3.5 5.25){$n$} \htext(-1.5 4.75){$n$} \htext(-3.5
4.75){$n$} \htext(-0.5 5.25){$n$} \htext(-2.5 5.25){$n$}
\vtext(-0.5 7.25){$\cdots$} \vtext(-1.5 7.25){$\cdots$}
\vtext(-2.5 7.25){$\cdots$} \vtext(-3.5 7.25){$\cdots$}
\end{texdraw}}%

\begin{center}
\raisebox{-\height}{\usebox{\tmpfigk}}
\end{center}

\end{example}

\vskip 3mm
\begin{defi}
Let $\Lambda$ be a dominant integral weight of level 1 for the
affine Lie algebra $\widehat{\mathfrak g}$.

(a) A column in a Young wall is called a {\it
full column} if its height is a multiple of the unit length and
its top is of unit thickness.

(b) For the classical quantum affine algebras
of type $A_{2n-1}^{(2)}$ $(n\ge 3)$, $D_n^{(1)}$ $(n\ge 4)$,
$A_{2n}^{(2)}$ $(n\ge 2)$, $D_{n+1}^{(2)}$ $(n\ge 2)$ and
$B_n^{(1)}$ $(n\ge 3)$, a Young wall is said to be {\it proper} if
none of the full columns have the same height.

(c) For the quantum affine algebras of type
$A_n^{(1)}$ $(n\ge 1)$, every Young wall is defined to be {\it
proper}.

\end{defi}

\vskip 3mm Let $\delta$ be the null root for the quantum affine
algebra $U_q(\widehat{\g})$ and write
\begin{equation*}
\begin{cases}
\delta = a_0\alpha_0 + a_1\alpha_1 + \cdots + a_n\alpha_n
&\text{for $\widehat{\g} = A_n^{(1)}, \cdots, B_n^{(1)}$},\\
2\delta = a_0\alpha_0 + a_1\alpha_1 + \cdots + a_n\alpha_n
&\text{for $\widehat{\g} = D_{n+1}^{(2)}$}.
\end{cases}
\end{equation*}
The part of a column consisting of $a_0$-many 0-blocks, $a_1$-many
1-blocks, $\cdots$, $a_n$-many $n$-blocks in some cyclic order is
called a {\it $\delta$-column}.

\vskip 3mm

\vskip 3mm

\begin{defi} (a) A column in a proper Young wall is said to {\it contain a
removable $\delta$} if we may remove a $\delta$-column from $Y$
and still obtain a proper Young wall.

\vskip 1mm (b) A proper Young wall is said to be {\it reduced} if
none of its columns contain a removable $\delta$.

\end{defi}

\vskip 3mm

Let $\F(\Lambda)$ be the set of all proper Young walls and let
$\Y(\Lambda)$ denote the set of all reduced proper Young walls.
Then we can define a crystal structure on $\F(\Lambda)$ so that it
may become a crystal graph for some integrable
$U_q(\widehat{\mathfrak g})$-module in the category
$\mathcal{O}_{int}$ \cite{Kang, KangKwon}. In this case, the set
$\Y(\Lambda)$ becomes a connected component in the crystal graph
$\F(\Lambda)$ and it is isomorphic to the crystal graph
$B(\Lambda)$ for the basic representation $V(\Lambda)$ of the
quantum affine algebra $U_q(\widehat{\mathfrak g})$. We briefly
explain the crystal structure of $\F(\Lambda)$. The main point is
how to define the action of Kashiwara operators $\eit$ and $\fit$
$(i=0, 1, \cdots, n)$ on proper Young walls.

\vskip 3mm
\begin{defi} (a) A block of color $i$  in a proper Young wall is called
a {\it removable $i$-block} if the wall remains a proper Young
wall after removing the block. A column in a proper Young wall is
called {\it $i$-removable} if the top of that column is a
removable $i$-block.

\vskip 1mm (b) A place in a proper Young wall where one may add an
$i$-block to obtain another proper Young wall is called an {\it
admissible $i$-slot}. A column in a proper Young wall is called
{\it $i$-admissible} if the top of that column is an admissible
$i$-slot.

\end{defi}

\vskip 3mm Fix $i\in I$ and let $Y=(y_k)_{k=0}^{\infty}\in
\F(\Lambda)$ be a proper Young wall.

\begin{itemize}
\item [(1)] To each column $y_k$ of $Y$, we assign its {\it $i$-signature} as
follows:
\begin{itemize}
\item [(a)] we assign $-\,-$ if the column $y_k$ is twice $i$-removable;
      (the $i$-block will be of half-unit height in this case).
\item [(b)] we assign $-$ if the column is once $i$-removable,
      but not $i$-admissible (the $i$-block may be of unit height or
      of half-unit height);
\item [(c)] we assign $-\,+$ if the column is once $i$-removable and
      once $i$-admissible (the $i$-block will be of half-unit height
      in this case);
\item [(d)] we assign $+$ if the column is once $i$-admissible, but not
      $i$-removable (the $i$-block may be of unit height or
      of half-unit height);
\item [(e)] we assign $+\,+$ if the column is twice $i$-admissible
      (the $i$-block will be of half-unit height in this case).
\end{itemize}
\item [(2)] From the (infinite) sequence of $+$'s and $-$'s, cancel out every
$(+,-)$-pair to obtain a finite sequence of $-$'s followed by
$+$'s, reading from left to right. This sequence is called the
{\it $i$-signature} of the proper Young wall $Y$.
\item [(3)] We define $\eit Y$ to be the proper Young wall obtained from $Y$
by removing the $i$-block corresponding to the right-most $-$ in
the $i$-signature of $Y$. We define $\eit Y = 0$ if there exists
no $-$ in the $i$-signature of $Y$.
\item [(4)] We define $\fit Y$ to be the proper Young wall obtained from $Y$
by adding an $i$-block to the column corresponding to the
left-most $+$ in the $i$-signature of $Y$. We define $\fit Y = 0$
if there exists no $+$ in the $i$-signature of $Y$.
\end{itemize}

\vskip 3mm Then we have:

\vskip 3mm
\begin{thm} {\rm \cite{Kang, KangKwon}}
{\rm(a)}\,  The set $\F(\Lambda)$ together with the Kashiwara
operators defined as above becomes a crystal graph for an
integrable $U_q(\widehat{\mathfrak g})$-module in the category
$\mathcal {O}_{int}$.

{\rm(b)}\, For all $i\in I$ and $Y\in \Y(\Lambda)$, we have
$$\eit Y\in \Y(\Lambda)\cup \{0\} \quad \text{and} \quad
\fit Y\in \Y(\Lambda)\cup \{0\}.$$ Moreover, there exists a
crystal isomorphism
\begin{equation*}
\Y(\Lambda) \stackrel{\sim} \longrightarrow B(\Lambda) \quad
\text{given by} \ \ Y_{\Lambda} \longmapsto u_{\Lambda},
\end{equation*}
where $u_{\Lambda}$ is the highest weight vector in $B(\Lambda)$.
\end{thm}

\vskip 3mm

\vskip 1cm
\section{Realization of Crystal Bases}
In this section, we will state the main result of this paper -- a
new realization of crystal bases for finite dimensional
irreducible modules over classical Lie algebras.

\vskip 3mm Let us explain the main idea of our approach. Let
$\frak g$ be a classical Lie algebra lying inside an affine Lie
algebra $\widehat{\frak g}$ so that the Dynkin diagram of $\frak
g$ can be obtained by removing the $0$-node from the Dynkin
diagram of $\widehat{\frak g}$. In this paper, we will focus on
the following pairs of a classical Lie algebra and an affine Lie
algebra:

$$A_n \subset A_n^{(1)},\,\, C_n \subset
A_{2n-1}^{(2)},\,\, B_n \subset B_n^{(1)},\,\, D_n \subset
D_n^{(1)}.$$

Fix such a pair $\frak g \subset \widehat{\frak g}$ and let
$\Lambda$ be a dominant integral weight of level 1 for the affine
Lie algebra $\widehat{\frak g}$. Then by Theorem 2.6, the crystal
graph $B(\Lambda)$ is realized as the set $\Y(\Lambda)$ of all
reduced proper Young walls built on the ground-state wall
$Y_{\Lambda}$. If we remove all the $0$-arrows in $\Y(\Lambda)$,
then it is decomposed into a disjoint union of infinitely many
connected components, each of which is isomorphic to the crystal
graph $B(\lambda)$ for some dominant integral weight $\lambda$ for
$\frak g$.

\vskip 2mm Conversely, any crystal graph $B(\lambda)$ for $\frak
g$ arises in this way. That is, given a dominant integral weight
$\lambda$ for $\frak g$, there is a dominant integral weight
$\Lambda$ of level $1$ for $\widehat{\frak g}$ such that
$B(\lambda)$ appears as a connected component in $B(\Lambda)$
without $0$-arrows. More precisely, we denote by $\la_i$
($i=1,\cdots,n$) and $\La_i$ ($i=0,1,\cdots,n$) the fundamental
weights for the quantum classical Lie algebras and the quantum
affine algebras, respectively and define the linear functionals
$\omega_i$ by

\vskip 2mm \hskip 3mm 1) $\frak g = A_n, \ C_n$:
\begin{equation*}
\omega_i = \la_i \quad \text{for} \ \ i=1, \cdots, n,
\end{equation*}

\hskip 3mm 2) $\frak g = B_n$:
\begin{equation*}
\omega_i = \begin{cases}
\la_i \quad & \text{for} \ \ i=1, \cdots, n-1, \\
2 \la_n \quad & \text{for} \ \ i=n,
\end{cases}
\end{equation*}

\hskip 3mm 3) $\frak g = D_n$:
\begin{equation*}
\omega_i = \begin{cases}
\la_i \quad & \text{for} \ \ i=1, \cdots, n-2, \\
\la_{n-1} + \la_n \quad & \text{for} \ \ i=n-1, \\
2 \la_n \quad & \text{for} \ \ i=n, \\
2 \la_{n-1} \quad & \text{for} \ \ i=n+1.
\end{cases}
\end{equation*}
Then for each dominant integral weight $\la$ for $\frak g$, we may
take the level $1$ dominant integral weight $\Lambda$ for
$\widehat{\frak g}$ as follows:

\vskip 3mm 1) $A_n\subset A_n^{(1)}$
$$\aligned
&\la=a_1\omega_1+\cdots+a_n\omega_n, \\
&\Lambda=\Lambda_i \,\,\,\,\,\text{if\,
$a_1+2a_2+\cdots+na_n\equiv i$\, mod $n+1$,}
\endaligned
$$

2) $C_n\subset A_{2n-1}^{(2)}$
$$\aligned
&\la=a_1\omega_1+\cdots+a_n\omega_n, \\
&\Lambda=\begin{cases} \Lambda_0 \quad &\text{if\,
$a_1+2a_2+\cdots+na_n$ is odd,}\\
\Lambda_1 \quad &\text{if\, $a_1+2a_2+\cdots+na_n$ is even,}
\end{cases}
\endaligned
$$

3) $B_n\subset B_n^{(1)}$
$$\aligned
&\la=a_1\omega_1+\cdots+a_n\omega_n+b\la_n, \\
&\Lambda=\begin{cases} \Lambda_0 \quad&\text{if\, $b=0$ and
$a_1+2a_2+\cdots+na_n$ is odd,}\\
\Lambda_1 \quad&\text{if\, $b=0$ and $a_1+2a_2+\cdots+na_n$ is even,}\\
\Lambda_n \quad&\text{if\, $b=1$,}
\end{cases}
\endaligned
$$

4) $D_n\subset D_n^{(1)}$
$$\aligned
&\la=a_1\omega_1+\cdots+a_{n+1}\omega_{n+1}+b_1\la_{n-1}+b_2\la_n, \\
&\Lambda=\begin{cases} \Lambda_0 \quad&\text{if\, $b_1=b_2=0$ and
$a_1+2a_2+\cdots+na_n+na_{n+1}$ is odd,}\\
\Lambda_1 \quad&\text{if\, $b_1=b_2=0$ and
$a_1+2a_2+\cdots+na_n+na_{n+1}$ is even,}\\
\Lambda_{n-1} \quad&\text{if\, $b_1=1$ and $b_2=0$,}\\
\Lambda_n \quad&\text{if\, $b_1=0$ and $b_2=1$.}
\end{cases}
\endaligned
$$

\vskip 3mm Now, we need to identify the highest weight vector
$u_{\lambda}$ for $B(\lambda)$ with some reduced proper Young wall
in $\Y(\Lambda)$ which is annihilated by all $\tilde{e_i}$ for
$i=1, \cdots, n$. However, given a dominant integral weight $\la$
for $\frak g$, there are infinitely many such Young walls in
$\Y(\Lambda)$. Equivalently, given $\lambda$, there are infinitely
many connected components of $\Y(\Lambda)$ without $0$-arrows that
are isomorphic to $B(\lambda)$. Thus the main task is to
characterize these connected components. Among these, we choose
the characterization of $B(\lambda)$ corresponding to the
connected components having the least number of blocks.

\savebox{\tmpfiga}{\begin{texdraw}\fontsize{6}{6}\selectfont
\drawdim em \setunitscale 0.9 \nc{\dtri}{ \bsegment \lvec(-2
0)\lvec(0 2)\lvec(0 0)\ifill f:0.7 \esegment }

\move(-2 0)\rlvec(15 0)\rlvec(0 13)\move(13 13)\rlvec(-2
0)\rlvec(0 -6)\move(13 11)\rlvec(-4 0)\rlvec(0 -4)\move(13
9)\rlvec(-6 0)\rlvec(0 -2)\move(4 4)\rlvec(0 -4)\move(13
2)\rlvec(-13 0)\rlvec(0 -2)\move(2 0)\rlvec(0 2)\move(4
4)\rlvec(0.4 0)

\move(0 0)\rlvec(2 2)\move(2 0)\rlvec(2 2)\move(11 0)\rlvec(2
2)\move(9 0)\rlvec(2 2)\move(7 0)\rlvec(2 2)\move(11 0)\rlvec(0
2)\move(9 0)\rlvec(0 2) \move(7 0)\rlvec(0 2)

\move(0 0) \clvec (0 -0.5)(3 -1)(4 -1)

\htext(6.5 -1){$i$ columns}

\move(13 0)\clvec (13 -0.5)(10 -1)(9 -1)

\htext(12 12){$i\!\!-\!\!1$} \rtext td:45 (5.7
5.5){$\cdots$}\htext(1.5 0.5){$0$}\htext(2.5 1.5){$0$}\htext(3.5
0.5){$1$}
\end{texdraw}}
\savebox{\tmpfigb}{\begin{texdraw}\fontsize{6}{6}\selectfont
\drawdim em \setunitscale 0.9 \nc{\dtri}{ \bsegment \lvec(-2
0)\lvec(0 2)\lvec(0 0)\ifill f:0.7 \esegment }

\setgray 0.6 \move(-2 0)\rlvec(15 0)\rlvec(0 25)\move(13
25)\rlvec(-2 0)\rlvec(0 -18)\move(13 23)\rlvec(-4 0)\rlvec(0
-16)\move(13 21)\rlvec(-6 0)\rlvec(0 -14)\move(4 17)\rlvec(-2
0)\rlvec(0 -4)\move(2 15)\rlvec(-2 0)\rlvec(0 -15)

\move(0 0)\rlvec(2 2)\move(11 0)\rlvec(2 2)\move(9 0)\rlvec(2
2)\move(7 0)\rlvec(2 2)\move(11 0)\rlvec(0 2)\move(9 0)\rlvec(0 2)
\move(7 0)\rlvec(0 2)\move(0 2)\rlvec(2 0)

\move(13 25)\rlvec(-2 -2)\htext(12.5 23.5){$1$}

\move(0 0)\setgray 0 \clvec (0 -0.5)(3 -1)(4 -1)

\htext(6.5 -1){$i$ columns}

\move(13 0)\clvec (13 -0.5)(10 -1)(9 -1)

\htext(0.5 1.5){$1$}\htext(1.5 0.5){$0$}\htext(2.5
1.5){$0$}\htext(3.5 0.5){$1$}\htext(1 14){$i$}

\htext(12 12){$i\!\!-\!\!1$} \rtext td:45 (5.7 5.5){$\cdots$}
\rtext td:45 (5.5 18.5){$\cdots$}

\move(0 0) \setgray 0 \linewd 0.12 \rlvec(13 0)\rlvec(0
13)\rlvec(-2 0)\rlvec(0 -6)\move(13 11)\rlvec(-4 0)\rlvec(0
-4)\move(13 9)\rlvec(-6 0)\rlvec(0 -2)\move(4 4)\rlvec(0
-4)\move(2 0)\rlvec(2 2)\move(13 2)\rlvec(-11 0)\rlvec(-2
-2)\move(2 0)\rlvec(0 2)\move(4 4)\rlvec(0.4 0)

\move(0 0)\rlvec(2 2)\move(11 0)\rlvec(2 2)\move(9 0)\rlvec(2
2)\move(7 0)\rlvec(2 2)\move(11 0)\rlvec(0 2)\move(9 0)\rlvec(0 2)
\move(7 0)\rlvec(0 2)
\end{texdraw}}
\savebox{\tmpfigc}{\begin{texdraw} \fontsize{6}{6}\selectfont
\drawdim em \setunitscale 0.9 \textref h:C v:C

\nc{\dtri}{ \bsegment \lvec(-2 0)\lvec(-2 1)\lvec(0 1)\lvec(0
0)\ifill f:0.7 \esegment }

\move(2 0)\lvec(2 1)\move(4 0)\lvec(4 1)\move(6 0)\lvec(6 1)
\move(8 0)\lvec(8 1) \move(10 0)\lvec(10 1)\move(12 0)\lvec(12 1)
\move(14 0)\lvec(14 1)\move(0 0)\rlvec(14 0)\move(0 1)\rlvec(14 0)

\htext(3 0.5){$n$} \htext(5 0.5){$n$} \htext(11 0.5){$n$}
\htext(13 0.5){$n$}\htext(1 0.5){$\cdots$}

\end{texdraw}}
\savebox{\tmpfigd}{\begin{texdraw} \fontsize{6}{6}\selectfont
\drawdim em \setunitscale 0.9 \textref h:C v:C

\nc{\dtri}{ \bsegment \lvec(-2 0)\lvec(-2 1)\lvec(0 1)\lvec(0
0)\ifill f:0.7 \esegment }

\bsegment

\setgray 0.6 \move(2 0)\lvec(2 1)\lvec(6 1) \move(4 0)\lvec(4
2)\lvec(6 2)\lvec(6 0) \move(10 0)\lvec(10 6)\lvec(14 6) \move(12
0)\lvec(12 8) \move(8 0)\lvec(8 1)\lvec(10 1) \move(12 8)\lvec(14
8)\lvec(14 0) \move(8 4)\lvec(12 4) \move(12 6)\lvec(14 8)
\move(10 2)\lvec(14 2)\move(14 0)\rlvec(0 8)

\htext(3 0.5){$n$} \htext(5 0.5){$n$} \htext(5 1.5){$n$} \htext(11
0.5){$n$} \htext(13 0.5){$n$} \htext(11 1.5){$n$} \htext(13
1.5){$n$} \htext(11 5){2} \htext(12.5 7.5){1} \htext(1
0.5){$\cdots$} \vtext(11 3){$\cdots$}

\rtext td:45 (7 3){$\cdots$}

\move(4 0)\setgray 0 \clvec (4 -0.5)(6 -1)(6.5 -1)

\htext(9 -1){$n$ columns}

\move(14 0)\clvec (14 -0.5)(12.2 -1)(11.7 -1)

\setgray 0 \linewd 0.12 \move(0 0)\lvec(14 0)\lvec(14 1)\lvec(0 1)
\esegment
\end{texdraw}}

\vskip 3mm Given a dominant integral weight $\la$ for $\frak g$,
we describe an algorithm of constructing the highest weight vector
$H_{\la}$ and lowest weight vector $L_{\la}$ inside $\Y(\Lambda)$.
For our convenience, we will focus on the case of $\frak g=B_n$
because this case contains all the characteristics of the
remaining cases. If $\la=\omega_i$ ($i=1,\cdots,n$), let
$H_{\omega_i}$ denote the following Young wall:
$$H_{\omega_i}=\raisebox{-0.4\height}{\usebox{\tmpfiga}}\,.$$
Then it is easy to verify that $\tilde{e_j}H_{\omega_i}=0$ for all
$j=1,\cdots,n$. That is, $H_{\omega_i}$ is a highest weight vector
of weight $\omega_i$. For the lowest weight vector, we denote by
$L_{\omega_i}$ the Young wall given below
$$L_{\omega_i}=\raisebox{-0.4\height}{\usebox{\tmpfigb}}\,.$$
Here, $H_{\omega_i}$ is denoted by the dark and bold-faced lines.
Note that $\tilde{f_j}L_{\omega_i}=0$ for all $j=1,\cdots,n$. Thus
$L_{\omega_i}$ is a lowest weight vector of weight $-\omega_i$. In
Theorem \ref{thm:An}, Theorem \ref{thm:Cn}, Theorem \ref{thm:Bn}
and Theorem \ref{thm:Dn}, we will show that $L_{\omega_i}$ is in
fact the lowest weight vector for the crystal graph $B(\omega_i)$;
i.e., $L_{\omega_i}$ and $H_{\omega_i}$ are connected by Kashiwara
operators.

If $\la=\la_n$, then the highest weight vector $H_{\la_n}$ and the
lowest weight vector $L_{\la_n}$ for $B(\la_n)$ are given by
$$H_{\la_n}=\raisebox{-0.1\height}{\usebox{\tmpfigc}}$$ and
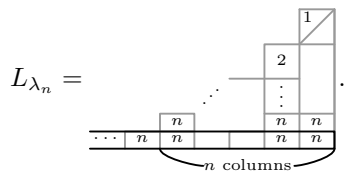
$$L_{\la_n}=\raisebox{-0.5\height}{\usebox{\tmpfigd}}\,.$$
Here, $H_{\la_n}$ is denoted by the dark and bold-faced lines.

\savebox{\tmpfiga}{\begin{texdraw}\fontsize{6}{6}\selectfont
\drawdim em \setunitscale 0.9 \nc{\dtri}{ \bsegment \lvec(-2
0)\lvec(0 2)\lvec(0 0)\ifill f:0.7 \esegment }

\rlvec(10 0)\rlvec(0 25)\rlvec(-1.6 0)\rlvec(0 -1.6)\rlvec(-1.6
0)\rlvec(0 -0.6)

\move(1.6 18.6)\rlvec(0 1.6)\rlvec(1.6 0) \move(1.6
18.6)\rlvec(-1.6 0) \rlvec(0 -18.6) \move(0 17)\rlvec(10 0)

\move(0 1.6)\rlvec(10 0)\move(1.6 0)\rlvec(0 1.6)\rlvec(-1.6 -1.6)
\move(3.2 0)\rlvec(0 1.6)\rlvec(-1.6 -1.6)\move(8.4 0)\rlvec(0
1.6)\move(8.4 0)\rlvec(1.6 1.6)

\htext(7 19){$H_{\omega_{i_k}}$}

\move(0 0)\clvec (0 -0.5)(2 -1)(2.5 -1)

\htext(5 -1){$i_k$ columns}

\move(10 0)\clvec (10 -0.5)(8 -1)(7.5 -1)

\move(10 0)\clvec (11 0)(11.5 5.5)(11.5 7.5)

\htext(12 8.5){$(t-k)$-many $\delta$-columns}

\move(10 17)\clvec (11 17)(11.5 11.5)(11.5 9.5)

\end{texdraw}}
\savebox{\tmpfigb}{\begin{texdraw}\fontsize{6}{6}\selectfont
\drawdim em \setunitscale 0.9 \nc{\dtri}{ \bsegment \lvec(-2
0)\lvec(0 2)\lvec(0 0)\ifill f:0.7 \esegment }

\rlvec(10 0)\rlvec(0 25)\rlvec(-1.6 0)\rlvec(0 -1.6)\rlvec(-1.6
0)\rlvec(0 -0.6)

\move(1.6 18.6)\rlvec(0 1.6)\rlvec(1.6 0) \move(1.6
18.6)\rlvec(-1.6 0) \rlvec(0 -18.6) \move(0 17)\rlvec(10 0)

\move(0 1.6)\rlvec(10 0)\move(1.6 0)\rlvec(0 1.6)\rlvec(-1.6 -1.6)
\move(3.2 0)\rlvec(0 1.6)\rlvec(-1.6 -1.6)\move(8.4 0)\rlvec(0
1.6)\move(8.4 0)\rlvec(1.6 1.6)

\move(0 0)\rlvec(0 -5)\rlvec(10 0)\rlvec(0 5) \move(0
-4.2)\rlvec(10 0)\move(1.6 -4.2)\rlvec(0 -0.8)\move(3.2
-4.2)\rlvec(0 -0.8)\move(8.4 -4.2)\rlvec(0 -0.8)

\htext(7 19){$H_{\omega_{i_k}}$}

\move(0 0)\clvec (0 -0.5)(2 -1)(2.5 -1)

\htext(5 -1){$i_k$ columns}

\move(10 0)\clvec (10 -0.5)(8 -1)(7.5 -1)

\move(10 0)\clvec (11 0)(11.5 5.5)(11.5 7.5)

\htext(12 8.5){$(t-k)$-many $\delta$-columns}

\move(10 17)\clvec (11 17)(11.5 11.5)(11.5 9.5)

\move(10 0)\clvec (11 0)(11.5 -1.5)(11.5 -2)

\htext(12 -2.5){$\frac{1}{2}\delta$-column}

\move(10 -5)\clvec (11 -5)(11.5 -3.5)(11.5 -3)
\end{texdraw}}

\vskip 3mm Suppose $\la$ has the form
$\la=\omega_{i_1}+\cdots+\omega_{i_t}$ ($1\le i_1\le \cdots\le
i_t\le n$). For each $k=1,\cdots,t$, let ${\overline
H_{\omega_{i_k}}}$ (resp. ${\overline L_{\omega_{i_k}}}$) denote
the Young wall consisting of $H_{\omega_{i_k}}$ (resp.
$L_{\omega_{i_k}}$) and $i_k\times (t-k)$-many $\delta$-columns.
Here, we place $H_{\omega_{i_k}}$ (resp. $L_{\omega_{i_k}}$) on
top of $\delta$-columns as is shown below.
$${\overline H_{\omega_{i_k}}}=\raisebox{-0.4\height}{\usebox{\tmpfiga}}$$
We define $H_{\la}$ (resp. $L_{\la}$) to be the Young wall
obtained by attaching ${\overline H_{\omega_{i_{k+1}}}}$ (resp.
${\overline L_{\omega_{i_{k+1}}}}$) to the left-hand side of
${\overline H_{\omega_{i_{k}}}}$ (resp. ${\overline
L_{\omega_{i_{k}}}}$) for $k=1,\cdots,t-1$.

\vskip 2mm On the other hand, suppose $\la$ has the form
$\la=\omega_{i_1}+\cdots+\omega_{i_t}+\la_n$ ($1\le i_1\le
\cdots\le i_t\le n$). For each $k=1,\cdots,t$, let ${\overline
H_{\omega_{i_k}}}$ (resp. ${\overline L_{\omega_{i_k}}}$) denote
the Young wall consisting of $H_{\omega_{i_k}}$ (resp.
$L_{\omega_{i_k}}$) and $i_k\times (t-k+\frac{1}{2})$-many
$\delta$-columns.
$${\overline
H_{\omega_{i_k}}}=\raisebox{-0.4\height}{\usebox{\tmpfigb}}$$ We
define $H_{\la}$ (resp. $L_{\la}$) to be the Young wall obtained
by attaching ${\overline H_{\omega_{i_{k+1}}}}$ (resp. ${\overline
L_{\omega_{i_{k+1}}}}$) to the left-hand side of ${\overline
H_{\omega_{i_{k}}}}$ (resp. ${\overline L_{\omega_{i_{k}}}}$) and
$H_{\la_n}$ (resp. $L_{\la_n}$) to the left-hand side of
${\overline H_{\omega_{i_t}}}$ (resp. ${\overline
L_{\omega_{i_t}}}$).

%
\savebox{\tmpfiga}{\begin{texdraw}\fontsize{6}{6}\selectfont
\drawdim em \setunitscale 0.9 \nc{\dtri}{ \bsegment \lvec(-2
0)\lvec(0 2)\lvec(0 0)\ifill f:0.7 \esegment }

\move(12 0)\dtri \move(10 0)\dtri

\move(15 7)\avec(12 7) \move(15 1)\avec(12 1) \htext(16
7){$L_{\la}$} \htext(16 1){$H_{\la}$}

\setgray 0.6 \move(6 0)\lvec(12 0)\rlvec(0 10)\move(12
10)\rlvec(-2 0)\rlvec(0 -10)\move(12 8)\rlvec(-2 0)\move(12
6)\rlvec(-2 0)\move(12 5)\rlvec(-2 0)\move(12 4)\rlvec(-2 0)
\move(12 2)\rlvec(-2 0)\move(10 2)\rlvec(-2 -2)\move(8 0)\rlvec(0
2)\rlvec(2 0)

\move(12 10)\rlvec(-2 -2)\move(12 2)\rlvec(-2 -2)

\move(12 2)\linewd 0.15 \setgray 0 \rlvec(-2 -2)\rlvec(2 0)
\rlvec(0 2)

\htext(10.5 9.5){$1$} \htext(10.5 1.5){$1$}\htext(11.5
0.5){$0$}\htext(9.5 0.5){$1$} \htext(11 7){$2$} \htext(11
5.5){$3$}\htext(11 4.5){$3$} \htext(11 3){$2$}

\end{texdraw}}
\savebox{\tmpfigb}{\begin{texdraw}\fontsize{6}{6}\selectfont
\drawdim em \setunitscale 0.9 \nc{\dtri}{ \bsegment \lvec(-2
0)\lvec(0 2)\lvec(0 0)\ifill f:0.7 \esegment }

\move(8 0)\dtri\move(10 0)\dtri\move(12 0)\dtri

\move(15 7)\avec(12 7) \move(15 1)\avec(12 1) \htext(16
7){$L_{\la}$} \htext(16 1){$H_{\la}$}

\setgray 0.6 \move(5 0)\lvec(12 0)\rlvec(0 10)\move(12
10)\rlvec(-2 0)\rlvec(0 -10)\move(12 8)\rlvec(-4 0)\rlvec(0 -8)
\move(12 6)\rlvec(-4 0)\move(12 5)\rlvec(-4 0) \move(12
4)\rlvec(-4 0) \move(12 2)\rlvec(-4 0)\move(6 0)\rlvec(0
2)\rlvec(2 0)

\move(12 10)\rlvec(-2 -2)\move(12 2)\rlvec(-2 -2)\move(10
2)\rlvec(-2 -2)\move(8 2)\rlvec(-2 -2)

\move(12 2)\linewd 0.15 \setgray 0 \rlvec(-2 0)\rlvec(0
-2)\rlvec(2 0) \rlvec(0 2)\move(10 0)\rlvec(-2 0)\rlvec(2 2)

\htext(11.5 8.5){$1$} \htext(7.5 0.5){$1$}\htext(8.5
1.5){$1$}\htext(9.5 0.5){$0$}\htext(10.5 1.5){$0$}\htext(11.5
0.5){$1$}

\htext(11 7){$2$} \htext(9 7){$2$}\htext(11 5.5){$3$}\htext(9
5.5){$3$}\htext(11 4.5){$3$}\htext(9 4.5){$3$} \htext(11
3){$2$}\htext(9 3){$2$}
\end{texdraw}}
\savebox{\tmpfigc}{\begin{texdraw}\fontsize{6}{6}\selectfont
\drawdim em \setunitscale 0.9 \nc{\dtri}{ \bsegment \lvec(-2
0)\lvec(-2 1)\lvec(0 1)\lvec(0 0)\ifill f:0.7 \esegment }

\move(0 0)\dtri\move(2 0)\dtri\move(4 0)\dtri\move(6
0)\dtri\move(8 0)\dtri

\move(11 11)\avec(8 11) \move(11 3)\avec(8 3) \htext(12
11){$L_{\la}$} \htext(12 3){$H_{\la}$}

\setgray 0.6 \move(-2 0)\lvec(8 0)\rlvec(0 14)\move(8 14)\rlvec(-2
0)\rlvec(0 -14)\move(8 12)\rlvec(-2 0) \move(8 10)\rlvec(-2
0)\move(8 9)\rlvec(-2 0)\move(8 8)\rlvec(-2 0)\move(8 6)\rlvec(-4
0)\rlvec(0 -6)\move(8 4)\rlvec(-6 0)\rlvec(0 -4)\move(8
2)\rlvec(-8 0)\rlvec(0 -2)\move(8 1)\rlvec(-12 0)

\move(8 14)\rlvec(-2 -2)\move(8 6)\rlvec(-2 -2)\move(6 6)\rlvec(-2
-2)

\move(-4 1)\linewd 0.15 \setgray 0 \rlvec(10 0)\rlvec(0 5)\rlvec(2
0) \rlvec(0 -6)\rlvec(-12 0)\move(6 4)\rlvec(2 2)

\htext(7.5 12.5){$1$} \htext(7 11){$2$}\htext(7 9.5){$3$}\htext(7
8.5){$3$}\htext(7 7){$2$}

\htext(6.5 5.5){$0$}\htext(7.5 4.5){$1$}\htext(4.5 5.5){$1$}

\htext(7 3){$2$}\htext(5 3){$2$}\htext(3 3){$2$}\htext(7 1.5){$3$}
\htext(5 1.5){$3$} \htext(3 1.5){$3$}\htext(1 1.5){$3$}\htext(7
0.5){$3$}\htext(5 0.5){$3$}\htext(3 0.5){$3$}\htext(1
0.5){$3$}\htext(-1 0.5){$3$}
\end{texdraw}}
\savebox{\tmpfigd}{\begin{texdraw}\fontsize{6}{6}\selectfont
\drawdim em \setunitscale 0.9 \nc{\dtri}{ \bsegment \lvec(-2
0)\lvec(-2 1)\lvec(0 1)\lvec(0 0)\ifill f:0.7 \esegment }

\move(2 0)\dtri\move(4 0)\dtri\move(6 0)\dtri\move(8
0)\dtri\move(10 0)\dtri\move(0 0)\dtri

\move(13 11)\avec(10 11) \move(13 3)\avec(10 3) \htext(14
11){$L_{\la}$} \htext(14 3){$H_{\la}$}

\setgray 0.6 \move(-2 0)\lvec(10 0)\rlvec(0 14)\move(10
14)\rlvec(-2 0)\rlvec(0 -14)\move(10 12)\rlvec(-4 0)\rlvec(0 -12)
\move(10 10)\rlvec(-4 0)\move(10 9)\rlvec(-4 0)\move(10
8)\rlvec(-4 0)\move(10 6)\rlvec(-6 0)\rlvec(0 -6)\move(10
4)\rlvec(-8 0)\rlvec(0 -4)\move(10 2)\rlvec(-10 0)\rlvec(0
-2)\move(10 1)\rlvec(-14 0)

\move(10 14)\rlvec(-2 -2)\move(10 6)\rlvec(-2 -2)\move(8
6)\rlvec(-2 -2)\move(6 6)\rlvec(-2 -2)

\move(-4 1)\linewd 0.15 \setgray 0 \rlvec(10 0)\rlvec(0 3)\rlvec(2
2) \rlvec(2 0)\rlvec(0 -6)\rlvec(-14 0)

\htext(9.5 12.5){$1$} \htext(9 11){$2$}\htext(7 11){$2$} \htext(9
9.5){$3$}\htext(7 9.5){$3$}\htext(9 8.5){$3$}\htext(7 8.5){$3$}
\htext(9 7){$2$}\htext(7 7){$2$} \htext(9.5 4.5){$1$}\htext(8.5
5.5){$0$}\htext(7.5 4.5){$0$}\htext(6.5 5.5){$1$}\htext(5.5
4.5){$1$} \htext(9 3){$2$}\htext(7 3){$2$}\htext(5 3){$2$}\htext(3
3){$2$} \htext(9 1.5){$3$}\htext(7 1.5){$3$}\htext(5 1.5){$3$}
\htext(3 1.5){$3$}\htext(1 1.5){$3$} \htext(9 0.5){$3$}\htext(7
0.5){$3$}\htext(5 0.5){$3$} \htext(3 0.5){$3$}\htext(1
0.5){$3$}\htext(-1 0.5){$3$}
\end{texdraw}}

\vskip 3mm \begin{example}\, In this example, we will give
descriptions of $H_{\la}$ and $L_{\la}$ for various choices of
dominant integral weights $\la$ for $\frak g=B_3$. The highest
weight vector $H_{\la}$ will be denoted by the dark, bold-faced
lines and the lowest weight vector $L_{\la}$ will be denoted by
the bright, dotted lines.

\vskip 3mm (a) If $\lambda=\omega_1$, we choose $\La = \La_1$ and
if $\la = \omega_2$, we choose $\Lambda=\Lambda_0$. The vectors
$H_{\la}$ and $L_{\la}$ are given by

$$\raisebox{-0.5\height}{\usebox{\tmpfiga}}\,\,\qquad
\qquad\raisebox{-0.5\height}{\usebox{\tmpfigb}}$$

\vskip 2mm (b) If $\lambda=\omega_3$, we choose $\La = \La_1$, and
if $\la = \omega_1+\omega_3$, we choose $\Lambda=\Lambda_0$. The
vectors $H_{\la}$ and $L_{\la}$ are given by

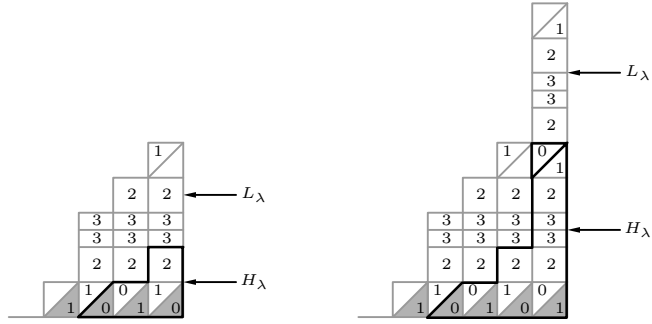
\begin{center}
\begin{texdraw}\fontsize{6}{6}\selectfont \drawdim em
\setunitscale 0.9 \nc{\dtri}{ \bsegment \lvec(-2 0)\lvec(0
2)\lvec(0 0)\ifill f:0.7 \esegment }

\bsegment \move(0 0)\dtri\move(2 0)\dtri\move(4 0)\dtri\move(6
0)\dtri

\move(9 7)\avec(6 7) \move(9 2)\avec(6 2) \htext(10 7){$L_{\la}$}
\htext(10 2){$H_{\la}$}

\setgray 0.6 \move(-4 0)\lvec(6 0)\rlvec(0 10)\move(6 10)\rlvec(-2
0)\rlvec(0 -10)\move(6 8)\rlvec(-4 0)\rlvec(0 -8)\move(6
6)\rlvec(-6 0)\rlvec(0 -6)\move(6 5)\rlvec(-6 0)\move(6
4)\rlvec(-6 0)\move(6 2)\rlvec(-6 0)\rlvec(-2 -2)\move(-2
0)\rlvec(0 2)\rlvec(2 0)

\move(6 10)\rlvec(-2 -2)\move(6 2)\rlvec(-2 -2)\move(4 2)\rlvec(-2
-2)\move(2 2)\rlvec(-2 -2)

\move(6 4)\linewd 0.15 \setgray 0 \rlvec(-2 0)\rlvec(0 -2)
\rlvec(-2 0)\rlvec(-2 -2) \rlvec(6 0)\rlvec(0 4)

\htext(4.5 9.5){$1$} \htext(5 7){$2$}\htext(3 7){$2$}\htext(5
5.5){$3$}\htext(3 5.5){$3$}\htext(1 5.5){$3$}\htext(5
4.5){$3$}\htext(3 4.5){$3$}\htext(1 4.5){$3$} \htext(5
3){$2$}\htext(3 3){$2$}\htext(1 3){$2$}\htext(5.5
0.5){$0$}\htext(4.5 1.5){$1$}\htext(3.5 0.5){$1$}\htext(2.5
1.5){$0$} \htext(1.5 0.5){$0$} \htext(0.5 1.5){$1$}\htext(-0.5
0.5){$1$}

\esegment

\move(22 0)\bsegment \move(-2 0)\dtri\move(0 0)\dtri\move(2
0)\dtri\move(4 0)\dtri\move(6 0)\dtri

\move(9 14)\avec(6 14) \move(9 5)\avec(6 5) \htext(10
14){$L_{\la}$} \htext(10 5){$H_{\la}$}

\setgray 0.6 \move(-6 0)\lvec(6 0)\rlvec(0 18)\move(6 18)\rlvec(-2
0)\rlvec(0 -18)\move(6 16)\rlvec(-2 0) \move(6 14)\rlvec(-2
0)\move(6 13)\rlvec(-2 0)\move(6 12)\rlvec(-2 0)\move(6
10)\rlvec(-4 0)\rlvec(0 -10) \move(6 8)\rlvec(-6 0)\rlvec(0
-8)\move(6 6)\rlvec(-6 0)\move(6 5)\rlvec(-6 0)\move(6 4)\rlvec(-6
0)\move(6 2)\rlvec(-6 0)\move(0 2)\rlvec(-2 -2)\move(-2 0)\rlvec(0
2)\rlvec(2 0)

\move(0 6)\rlvec(-2 0)\rlvec(0 -6)\move(0 5)\rlvec(-2 0)\move(0
4)\rlvec(-2 0)\move(-2 2)\rlvec(-2 0)\rlvec(0 -2)\rlvec(2 2)

\move(6 18)\rlvec(-2 -2)\move(6 10)\rlvec(-2 -2)\move(4
10)\rlvec(-2 -2)\move(6 2)\rlvec(-2 -2) \move(4 2)\rlvec(-2
-2)\move(2 2)\rlvec(-2 -2)

\move(6 10)\linewd 0.15 \setgray 0 \rlvec(-2 -2)\rlvec(0
-4)\rlvec(-2 0)\rlvec(0 -2) \rlvec(-2 0)\rlvec(-2 -2)\rlvec(8
0)\rlvec(0 10)\rlvec(-2 0)\rlvec(0 -2)

\htext(5.5 16.5){$1$} \htext(5 15){$2$}\htext(5 13.5){$3$}\htext(5
12.5){$3$}\htext(5 11){$2$}\htext(4.5 9.5){$0$}\htext(5.5
8.5){$1$}\htext(2.5 9.5){$1$}\htext(5 7){$2$}\htext(3
7){$2$}\htext(1 7){$2$}\htext(5 5.5){$3$}\htext(3 5.5){$3$}
\htext(1 5.5){$3$} \htext(5 4.5){$3$}\htext(3 4.5){$3$}\htext(1
4.5){$3$}\htext(1 3){$2$}\htext(5 3){$2$}\htext(3
3){$2$}\htext(5.5 0.5){$1$}\htext(4.5 1.5){$0$}\htext(3.5
0.5){$0$}\htext(2.5 1.5){$1$}\htext(1.5 0.5){$1$}\htext(0.5
1.5){$0$}\htext(-0.5 0.5){$0$}\htext(-1.5 1.5){$1$}\htext(-2.5
0.5){$1$}\htext(-1 3){$2$}\htext(-1 4.5){$3$}\htext(-1 5.5){$3$}

\esegment
\end{texdraw}
\end{center}

\vskip 2mm (c) If $\lambda=\omega_1+\la_3$ or $\omega_2+\la_3$, we
choose $\Lambda=\Lambda_3$. The vectors $H_{\la}$ and $L_{\la}$
are given by

$$\raisebox{-0.5\height}{\usebox{\tmpfigc}}\,\,\qquad
\qquad\raisebox{-0.5\height}{\usebox{\tmpfigd}}$$

\end{example}

\vskip 5mm We now begin to characterize the crystal graph $B(\la)$
inside $\Y(\La)$. Let $F(\la)$ denote the set of all reduced
proper Young walls lying between $H_{\la}$ and $L_{\la}$. To
describe $B(\la)$ inside $F(\la)$, we need some  additional
conditions. For this purpose, we need to introduce some notations.
Fix a dominant integral weight $\la$ as follows:
\begin{equation}
\la=\left\{\begin{array}{ll} \omega_{i_1}+\cdots+\omega_{i_t}
&\text{if $\frak g=A_n,C_n$},\\
\omega_{i_1}+\cdots+\omega_{i_t}+b\la_n
&\text{if $\frak g=B_n$},\\
\omega_{i_1}+\cdots+\omega_{i_t}+b_1\la_{n-1}+b_2\la_n &\text{if
$\frak g=D_n$}, \end{array}\right.
\end{equation}
where $b=0$ or $1$, $(b_1,b_2)=(1,0)$ or $(0, 1)$.

\vskip 2mm For each $Y \in F(\la)$, we denote by $\overset {\circ}
{Y}_{\omega_{i_k}}$ ($k=1,\cdots,t$) (resp. $\overset {\circ}
{Y}_{\lambda_{n-1}},$ $\overset {\circ}{Y}_{\la_n}$) the part of
$Y$ consisting of the blocks lying above ${\overline
H_{\omega_{i_k}}}$ (resp. $H_{\lambda_{n-1}},$ $H_{\la_n}$) and we
denote by ${\overline Y_{\omega_{i_k}}}$ (resp. ${\overline
Y_{\la_{n-1}}}$, ${\overline Y_{\la_n}}$)  the intersection of $Y$
and ${\overline L_{\omega_{i_k}}}$ (resp. $L_{\la_{n-1}}$,
$L_{\la_n}$) as is shown in the following picture. Moreover, we
denote by ${\overline Y_{\omega_{i_k}+\omega_{i_{k+1}}}}$ (resp.
${\overline Y_{\omega_{i_t}+\la_{n-1}}}$, ${\overline
Y_{\omega_{i_t}+\la_n}}$) the union of ${\overline
Y_{\omega_{i_k}}}$ (resp. ${\overline Y_{\omega_{i_t}}}$) and
${\overline Y_{\omega_{i_{k+1}}}}$ (resp. ${\overline
Y_{\la_{n-1}}}$ or ${\overline Y_{\la_n}}$).

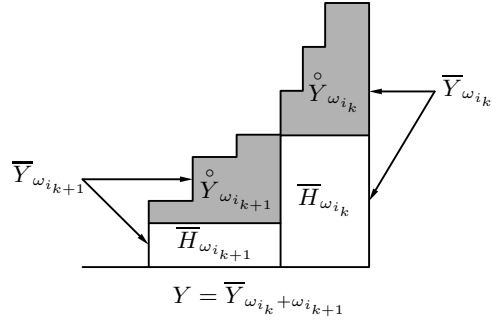
\begin{center}
\begin{texdraw}\fontsize{9}{9}\selectfont \drawdim em
\setunitscale 0.9 \nc{\dtri}{ \bsegment \move(0 0)\rlvec(-2
0)\rlvec(0 -2)\rlvec(-1 0)\rlvec(0 -2)\rlvec(-1 0)\rlvec(0
-2)\rlvec(4 0)\rlvec(0 6)\ifill f:0.7 \esegment }

\nc{\dtrii}{ \bsegment \move(0 0)\rlvec(-2 0)\rlvec(0 -1)\rlvec(-2
0)\rlvec(0 -2)\rlvec(-2 0)\rlvec(0 -1)\rlvec(6 0)\rlvec(0 4)\ifill
f:0.7 \esegment }

\move(10 12)\dtri \move(6 6)\dtrii

\move(-3 0)\lvec(10 0)\rlvec(0 12)\move(10 12)\rlvec(-2 0)\rlvec(0
-2)\rlvec(-1 0)\rlvec(0 -2)\rlvec(-1 0)\rlvec(0 -2)\rlvec(-2
0)\rlvec(0 -1)\rlvec(-2 0)\rlvec(0 -2)\rlvec(-2 0)\rlvec(0 -3)

\move(0 2)\rlvec(6 0)\move(6 0)\rlvec(0 6)\rlvec(4 0)

\move(-3 4)\avec(0 1)\move(-3 4)\avec(2 4)

\move(13 8)\avec(10 3)\move(13 8)\avec(10 8)

\htext(-4.5 4){${\overline Y_{\omega_{i_{k+1}}}}$}

\htext(14.5 8){${\overline Y_{\omega_{i_k}}}$}

\htext(8 3){${\overline H_{\omega_{i_k}}}$}\htext(3 1){${\overline
H_{\omega_{i_{k+1}}}}$}\htext(4 3.5){$\overset {\circ}
{Y}_{\omega_{i_{k+1}}}$}\htext(8.4 8){$\overset
{\circ}{Y}_{\omega_{i_k}}$}

\htext(5 -1.5){$Y={\overline Y_{\omega_{i_k}+\omega_{i_{k+1}}}}$}
\end{texdraw}
\end{center}

\vskip 3mm Now, consider ${\overline
Y_{\omega_{i_k}+\omega_{i_{k+1}}}}$, ${\overline
Y_{\omega_{i_t}+\la_{n-1}}}$ and ${\overline
Y_{\omega_{i_t}+\la_n}}$ of $Y$. Then we define
\begin{equation}
\begin{aligned}
& Y^{\omega_{i_k}} = \overset {\circ}{Y}_{\omega_{i_k}} \cap
{\overline L}_{\omega_{i_{k+1}}} \quad
\text{reading from top to bottom}, \\
& Y^{\omega_{i_{k+1}}} = \overset {\circ}{Y}_{\omega_{i_{k+1}}}
\cap {\overline L}_{\omega_{i_k}} \quad \text{reading from right
to left}\,\,\, \text{in}\,\,\,{\overline
Y_{\omega_{i_k}+\omega_{i_{k+1}}}}.
\end{aligned}
\end{equation}
Similarly, we define
\begin{equation}
\aligned &Y^{\omega_{i_t}}=\overset {\circ}{Y}_{\omega_{i_t}}\cap
L_{\la_{n-1}},\quad Y^{\la_{n-1}} = \overset
{\circ}{Y}_{\la_{n-1}} \cap {\overline L}_{\omega_{i_t}}\quad
\text{in}\,\,\,{\overline
Y_{\omega_{i_t}+\la_{n-1}}},\\
&Y^{\omega_{i_t}}=\overset {\circ}{Y}_{\omega_{i_t}}\cap
L_{\la_{n}},\quad Y^{\la_{n}} = \overset {\circ}{Y}_{\la_{n}} \cap
{\overline L}_{\omega_{i_t}}\quad \text{in}\,\,\,{\overline
Y_{\omega_{i_t}+\la_{n}}}.
\endaligned
\end{equation}

\savebox{\tmpfiga}{\begin{texdraw}\fontsize{6}{6}\selectfont
\drawdim em \setunitscale 0.9 \nc{\dtri}{ \bsegment \lvec(-2
0)\lvec(-2 2)\lvec(0 2)\lvec(0 0)\ifill f:0.7 \esegment }

\bsegment

\setgray 0.6 \move(0 0)\lvec(0 10)\move(0 10)\rlvec(-2 0)\rlvec(0
-10) \move(0 8) \rlvec(-4 0)\rlvec(0 -8)\move(0 6)\rlvec(-4
0)\move(0 4)\rlvec(-6 0)\rlvec(0 -4)\move(0 2)\rlvec(-10
0)\rlvec(0 -2)\move(-8 0)\rlvec(0 2) \move(-12 0)\lvec(0
0)\move(-6 4)\rlvec(-2 0)\rlvec(0 -2)

\move(0 0)\linewd 0.15 \setgray 0 \lvec(0 4)\rlvec(-4 0)\rlvec(0
-4)\rlvec(4 0)

\htext(-1 9){$4$} \htext(-1 7){$3$} \htext(-3 7){$2$} \htext(-1
5){$2$} \htext(-3 5){$1$} \htext(-1 3){$1$} \htext(-3
3){$0$}\htext(-5 3){$4$} \htext(-1 1){$0$} \htext(-3
1){$4$}\htext(-5 1){$3$} \htext(-7 1){$2$} \htext(-7
3){$3$}\htext(-9 1){$1$} \esegment
\end{texdraw}}
\savebox{\tmpfigb}{\begin{texdraw}\fontsize{6}{6}\selectfont
\drawdim em \setunitscale 0.9 \nc{\dtri}{ \bsegment \lvec(-2
0)\lvec(-2 2)\lvec(0 2)\lvec(0 0)\ifill f:0.7 \esegment }

\rlvec(4 0)\rlvec(0 6)\rlvec(-2 0)\rlvec(0 -6)\move(0 0)\rlvec(0
4)\rlvec(4 0)\move(0 2)\rlvec(4 0)

\htext(1 1){$1$} \htext(1 3){$2$} \htext(3 1){$2$} \htext(3
3){$3$} \htext(3 5){$4$}
\end{texdraw}}
\savebox{\tmpfigc}{\begin{texdraw}\fontsize{6}{6}\selectfont
\drawdim em \setunitscale 0.9 \nc{\dtri}{ \bsegment \lvec(-2
0)\lvec(-2 2)\lvec(0 2)\lvec(0 0)\ifill f:0.7 \esegment }

\rlvec(6 0)\rlvec(0 4)\rlvec(-4 0)\rlvec(0 -4)\move(0 0)\rlvec(0
2)\rlvec(6 0)\move(4 0)\rlvec(0 4)

\htext(1 1){$1$} \htext(3 1){$2$} \htext(3 3){$3$} \htext(5
1){$3$} \htext(5 3){$4$}
\end{texdraw}}
\savebox{\tmpfigg}{\begin{texdraw}\fontsize{6}{6}\selectfont
\drawdim em \setunitscale 0.9 \nc{\dtri}{ \bsegment \lvec(-2
0)\lvec(-2 2)\lvec(0 2)\lvec(0 0)\ifill f:0.7 \esegment }

\rlvec(4 0)\rlvec(0 4)\rlvec(-2 0)\rlvec(0 -4)\move(0 0)\rlvec(0
2)\rlvec(4 0)

\htext(1 1){$2$} \htext(3 1){$3$} \htext(3 3){$4$}
\end{texdraw}}
\savebox{\tmpfigh}{\begin{texdraw}\fontsize{6}{6}\selectfont
\drawdim em \setunitscale 0.9 \nc{\dtri}{ \bsegment \lvec(-2
0)\lvec(-2 2)\lvec(0 2)\lvec(0 0)\ifill f:0.7 \esegment }

\rlvec(0 4)\rlvec(4 0)\rlvec(0 -4)\rlvec(-4 0)\move(0 2)\rlvec(4
0)\move(2 0)\rlvec(0 4)

\htext(1 1){$2$} \htext(1 3){$3$} \htext(3 3){$4$} \htext(3
1){$3$}
\end{texdraw}}
\savebox{\tmpfigk}{\begin{texdraw}\fontsize{6}{6}\selectfont
\drawdim em \setunitscale 0.9 \nc{\dtri}{ \bsegment \lvec(-2
0)\lvec(-2 2)\lvec(0 2)\lvec(0 0)\ifill f:0.7 \esegment }

\lvec(0 10)\rlvec(-2 0)\rlvec(0 -10) \move(0 8) \rlvec(-4
0)\rlvec(0 -8)\move(0 6)\rlvec(-4 0)\move(0 4)\rlvec(-4 0)\move(0
2)\rlvec(-4 0)\move(-4 0)\lvec(0 0)

\htext(-1 9){$4$} \htext(-1 7){$3$} \htext(-3 7){$2$} \htext(-1
5){$2$} \htext(-3 5){$1$} \htext(-1 3){$1$} \htext(-3
3){$0$}

\htext(-1 1){$0$} \htext(-3 1){$4$}
\end{texdraw}}

\begin{example}\label{exm:Y_i}
If $\frak g=A_4$, $\la=\omega_2+\omega_3$, and
$$Y=\raisebox{-0.5\height}{\usebox{\tmpfiga}}\in F(\la),$$
then we have
$$\overset {\circ}{Y}_{\omega_2}=\raisebox{-0.5\height}{\usebox{\tmpfigb}}\,,\quad
\overset {\circ}
{Y}_{\omega_3}=\raisebox{-0.5\height}{\usebox{\tmpfigc}}\,,\quad
{\overline
Y_{\omega_2}}=\raisebox{-0.5\height}{\usebox{\tmpfigk}}\quad
\text{and}\quad {\overline
Y_{\omega_3}}=\raisebox{-0.5\height}{\usebox{\tmpfigc}}\,.$$
Moreover, we have
$$Y^{\omega_2}=\raisebox{-0.5\height}{\usebox{\tmpfigg}}\quad
\text{and}\quad
Y^{\omega_3}=\raisebox{-0.5\height}{\usebox{\tmpfigh}}\,.$$
\end{example}

With these notations, we are ready to give an explicit description
of the crystal graph $B(\la)$ over $\frak g = A_n$.


\vskip 3mm
\begin{thm} \label{thm:An}
Let $\la \in P^{+}$ be a dominant integral weight
and write
$$\la = \omega_{i_1} + \cdots + \omega_{i_t}
\ \ (1 \le i_1 \le \cdots \le i_t \le n).$$ Set
$$ Y(\la)=\{Y\in F(\la)\,|\, Y^{\omega_{i_k}}\subset Y^{\omega_{i_{k+1}}}\,\,
\text{in ${\overline Y}_{\omega_{i_k}+\omega_{i_{k+1}}}$ for
all}\,\, k=1,2,\cdots,t-1\}.
$$ Then there exists an isomorphism of $U_q(A_n)$-crystals
\begin{equation}
Y(\lambda) \stackrel{\sim} \longrightarrow B(\la) \quad
\text{given by} \ \ H_{\la} \longmapsto u_{\la},
\end{equation}
where $u_\la$ is the highest weight vector in $B(\la)$.
\end{thm}

\vskip 3mm
\begin{example}

Let $\frak g=A_4$ and $\la = \omega_2 + \omega_3$.
For each Young wall given below, the shaded part represents
$Y^{\omega_2}$ and $Y^{\omega_3}$, respectively.
Hence, by Theorem \ref{thm:An}, the first Young wall belongs to
$Y(\la)$, but the second one doesn't.

\vskip 3mm

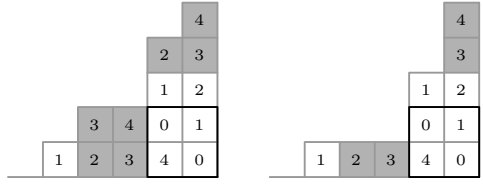
\begin{center}
\begin{texdraw}
\fontsize{6}{6}\selectfont \drawdim em \setunitscale 0.9
\nc{\dtri}{ \bsegment \lvec(-2 0) \lvec(-2 2)\lvec(0 2)\lvec(0
0)\ifill f:0.7 \esegment }

\nc{\dtrii}{ \bsegment \lvec(-2 0) \lvec(-2 -2)\lvec(0 0)\ifill
f:0.7 \esegment }

\bsegment \move(0 8)\dtri\move(0 6)\dtri \move(-4 0)\dtri\move(-6
0)\dtri

\setgray 0.6 \move(0 0)\lvec(0 10)\move(0 10)\rlvec(-2 0)\rlvec(0
-10) \move(0 8) \rlvec(-2 0)\move(0 6)\rlvec(-4 0)\rlvec(0
-6)\move(0 4)\rlvec(-4 0)\move(0 2)\rlvec(-10 0)\rlvec(0
-2)\move(-8 0)\rlvec(0 2)\move(-6 0)\rlvec(0 2) \move(-12
0)\lvec(0 0)

\linewd 0.12 \setgray 0 \lvec(0 4)\rlvec(-4 0)\rlvec(0 -4)\rlvec(4
0)

\htext(-1 9){$4$} \htext(-1 7){$3$}\htext(-1 5){$2$} \htext(-3
5){$1$} \htext(-1 3){$1$} \htext(-3 3){$0$} \htext(-1 1){$0$}
\htext(-3 1){$4$}\htext(-5 1){$3$} \htext(-7 1){$2$} \htext(-9
1){$1$} \esegment

\move(-15 0)\bsegment \move(0 6)\dtri \move(-2 6)\dtri \move(0
8)\dtri \move(-4 0)\dtri \move(-6 0)\dtri \move(-6 2)\dtri
\move(-4 2)\dtri

\setgray 0.6 \move(0 0)\lvec(0 10)\move(0 10)\rlvec(-2 0)\rlvec(0
-10) \move(0 8) \rlvec(-4 0)\rlvec(0 -8)\move(0 6)\rlvec(-4
0)\move(0 4)\rlvec(-8 0)\rlvec(0 -4)\move(-6 4)\rlvec(0 -4)
\move(0 2)\rlvec(-10 0)\rlvec(0 -2)\move(-8 0)\rlvec(0 2)
\move(-12 0)\lvec(0 0)

\linewd 0.12 \setgray 0 \lvec(0 4)\rlvec(-4 0)\rlvec(0 -4)\rlvec(4
0)

\htext(-1 9){$4$} \htext(-1 7){$3$} \htext(-3 7){$2$} \htext(-1
5){$2$} \htext(-3 5){$1$} \htext(-1 3){$1$} \htext(-3
3){$0$}\htext(-5 3){$4$} \htext(-1 1){$0$} \htext(-3
1){$4$}\htext(-5 1){$3$}\htext(-7 3){$3$} \htext(-7 1){$2$}
\htext(-9 1){$1$} \esegment

\end{texdraw}
\end{center}
\end{example}


\vskip 3mm Next, we will consider the case when $\frak g = C_n$ or
$B_n$. Consider ${\overline Y_{\omega_{i_k}}}$ for $k=1,\cdots,t$.
Suppose that ${\overline Y_{\omega_{i_k}}}$ contains a row
consisting of $n$-blocks, which will be called the {\it $n$-row},
as is shown in the following picture.

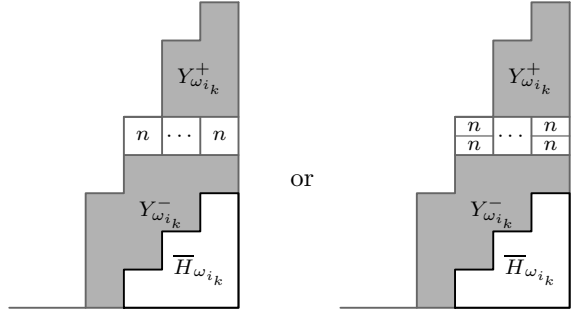
\savebox{\tmpfiga}{\begin{texdraw} \fontsize{8}{8}\selectfont
\drawdim em \setunitscale 1.7 \textref h:C v:C

\bsegment \setgray 0.6 \move(4 5)\lvec(6 5)\lvec(6 8) \lvec(5
8)\lvec(5 7)\lvec(4 7)\lvec(4 5) \ifill f:0.7 \esegment

\bsegment \setgray 0.6 \move(2 0)\lvec(2 3)\lvec(3 3)\lvec(3
4)\lvec(6 4)\lvec(6 3)\lvec(5 3) \lvec(5 2)\lvec(4 2)\lvec(4
1)\lvec(3 1)\lvec(3 0)\lvec(2 0) \ifill f:0.7 \esegment

\bsegment \setgray 0.4 \move(0 0)\lvec(6 0)\lvec(6 8) \move(2
0)\lvec(2 3)\lvec(3 3)\lvec(3 5)\lvec(4 5) \lvec(4 7)\lvec(5
7)\lvec(5 8)\lvec(6 8) \move(3 4)\lvec(6 4) \move(3 5)\lvec(6 5)
\move(4 5)\lvec(4 4) \move(5 5)\lvec(5 4) \esegment

\bsegment \setgray 0 \move(3 0)\lvec(3 1)\lvec(4 1)\lvec(4
2)\lvec(5 2)\lvec(5 3)\lvec(6 3)\lvec(6 0)\lvec(3 0) \esegment

\htext(3.5 4.5){$n$} \htext(5.5 4.5){$n$} \htext(4.5
4.5){$\cdots$} \htext(5 1){${\overline H}_{\omega_{i_k}}$}
\htext(5 6){$Y_{\omega_{i_k}}^+$} \htext(4
2.5){$Y_{\omega_{i_k}}^-$}

\end{texdraw}}
\savebox{\tmpfigb}{\begin{texdraw} \fontsize{8}{8}\selectfont
\drawdim em \setunitscale 1.7 \textref h:C v:C

\setgray 0.6 \move(4 5)\lvec(6 5)\lvec(6 8) \lvec(5 8)\lvec(5
7)\lvec(4 7)\lvec(4 5) \ifill f:0.7

\setgray 0.6 \move(2 0)\lvec(2 3)\lvec(3 3)\lvec(3 4)\lvec(6
4)\lvec(6 3)\lvec(5 3) \lvec(5 2)\lvec(4 2)\lvec(4 1)\lvec(3
1)\lvec(3 0)\lvec(2 0) \ifill f:0.7

\move(3 4.5)\rlvec(1 0)\move(5 4.5)\rlvec(1 0)

\setgray 0.4 \move(0 0)\lvec(6 0)\lvec(6 8) \move(2 0)\lvec(2
3)\lvec(3 3)\lvec(3 5)\lvec(4 5) \lvec(4 7)\lvec(5 7)\lvec(5
8)\lvec(6 8) \move(3 4)\lvec(6 4) \move(3 5)\lvec(6 5) \move(4
5)\lvec(4 4) \move(5 5)\lvec(5 4)

\setgray 0 \move(3 0)\lvec(3 1)\lvec(4 1)\lvec(4 2)\lvec(5
2)\lvec(5 3)\lvec(6 3)\lvec(6 0)\lvec(3 0)

\htext(3.5 4.25){$n$} \htext(5.5 4.25){$n$}\htext(3.5 4.75){$n$}
\htext(5.5 4.75){$n$} \htext(4.5 4.5){$\cdots$} \htext(5
1){${\overline H}_{\omega_{i_k}}$} \htext(5
6){$Y_{\omega_{i_k}}^+$} \htext(4 2.5){$Y_{\omega_{i_k}}^-$}

\end{texdraw}}

\begin{center}
\raisebox{-0.4\height}{\usebox{\tmpfiga}}\qquad or\quad
\raisebox{-0.4\height}{\usebox{\tmpfigb}}

\end{center}

We will denote by $Y_{\omega_{i_k}}^{+}$ (resp.
$Y_{\omega_{i_k}}^{-}$) the part of $Y$ consisting of the blocks
lying above (resp. below) the $n$-row and below ${\overline
L}_{\omega_{i_k}}$ (resp. above ${\overline H}_{\omega_{i_k}}$).
We also denote by $|Y_{\omega_{i_k}}^{-}|$ the wall obtained by
reflecting $Y_{\omega_{i_k}}^{-}$ along the $n$-row and shifting
the blocks to the right as much as possible.

\vskip 3mm
\savebox{\tmpfiga}{\begin{texdraw} \fontsize{6}{6}\selectfont
\drawdim em \setunitscale 0.9 \nc{\dtri}{ \bsegment \lvec(-2 0)
\lvec(-2 2)\lvec(0 2)\lvec(0 0)\ifill f:0.7 \esegment }

\nc{\dtrii}{ \bsegment \lvec(-2 0) \lvec(-2 -2)\lvec(0 0)\ifill
f:0.7 \esegment }

\nc{\dtriii}{ \bsegment \lvec(-2 0) \lvec(0 2)\lvec(0 0)\ifill
f:0.7 \esegment }

\setgray 0.6 \move(8 0)\lvec(8 12)\move(8 12)\rlvec(-2 0)\rlvec(0
-12) \move(8 10)\rlvec(-4 0)\rlvec(0 -10) \move(8 8)\rlvec(-4
0)\move(8 6)\rlvec(-6 0)\rlvec(0 -6)\move(8 4)\rlvec(-8 0)\rlvec(0
-4)\move(8 2)\rlvec(-8 0)\move(-3 0)\lvec(8 0)

\move(2 2)\rlvec(-2 -2)\move(4 2)\rlvec(-2 -2)\move(6 2)\rlvec(-2
-2)\move(8 2)\rlvec(-2 -2)

\linewd 0.15 \setgray 0 \move(8 6)\rlvec(-2 0)\rlvec(0
-2)\rlvec(-2 0)\rlvec(0 -2) \rlvec(-2 0)\rlvec(-2 -2)\rlvec(8
0)\rlvec(0 6)

\htext(7 11){$2$} \htext(5 9){$3$} \htext(7 9){$3$}

\htext(5 7){$4$}\htext(7 7){$4$} \htext(3 5){$3$}\htext(5
5){$3$}\htext(7 5){$3$} \htext(1 3){$2$}\htext(3 3){$2$}\htext(5
3){$2$}\htext(7 3){$2$}

\htext(0.5 1.5){$1$} \htext(1.5 0.5){$0$} \htext(2.5 1.5){$0$}
\htext(3.5 0.5){$1$} \htext(4.5 1.5){$1$} \htext(5.5 0.5){$0$}
\htext(7.5 0.5){$1$} \htext(6.5 1.5){$0$}
\end{texdraw}}
\savebox{\tmpfigb}{\begin{texdraw} \fontsize{6}{6}\selectfont
\drawdim em \setunitscale 0.9

\rlvec(4 0)\rlvec(0 4)\rlvec(-2 0)\rlvec(0 -4) \move(0 0)\rlvec(0
2)\rlvec(4 0)

\htext(1 1){$3$} \htext(3 1){$3$} \htext(3 3){$2$}

\end{texdraw}}
\savebox{\tmpfigc}{\begin{texdraw} \fontsize{6}{6}\selectfont
\drawdim em \setunitscale 0.9

\rlvec(0 4)\rlvec(6 0)\rlvec(0 2)\rlvec(-4 0)\rlvec(0 -4)\rlvec(-2
-2) \move(0 2)\rlvec(4 0)\rlvec(0 4)\move(0 0)\rlvec(2 0)\rlvec(0
2)

\htext(0.5 1.5){$1$}\htext(1 3){$2$} \htext(3 3){$2$} \htext(3
5){$3$}\htext(5 5){$3$}

\end{texdraw}}
\savebox{\tmpfigd}{\begin{texdraw} \fontsize{6}{6}\selectfont
\drawdim em \setunitscale 0.9

\rlvec(4 0)\rlvec(0 6)\rlvec(-2 -2)\rlvec(0 -4)\move(0 0)\rlvec(0
4) \rlvec(4 0)\move(0 2)\rlvec(4 0)\move(2 4)\rlvec(0 2)\rlvec(2
0)

\htext(3.5 4.5){$1$}\htext(1 1){$3$} \htext(3 1){$3$} \htext(3
3){$2$}\htext(1 3){$2$}

\end{texdraw}}

\begin{example}
If $\frak g=C_4$, $\la=\omega_4$, and
$$Y= \raisebox{-0.4\height}{\usebox{\tmpfiga}}\in F(\la),$$
then we have
$$Y_{\omega_4}^+=\raisebox{-0.4\height}{\usebox{\tmpfigb}}\,,\quad
Y_{\omega_4}^-=\raisebox{-0.4\height}{\usebox{\tmpfigc}}\quad
\text{and}\quad
|Y_{\omega_4}^-|=\raisebox{-0.4\height}{\usebox{\tmpfigd}}\,.$$

\end{example}

\savebox{\tmpfiga} {\begin{texdraw} \fontsize{7}{7}\selectfont
\drawdim em \setunitscale 1.7 \textref h:C v:C

\bsegment \setgray 0.6 \move(4.5 9.5)\lvec(6.5 9.5)\lvec(4.5
7.5)\lvec(4.5 9.5)\ifill f:0.7 \esegment

\bsegment \setgray 0.6 \move(1.5 3.5)\lvec(3.5 3.5)\lvec(3.5
5.5)\lvec(3.5 5.5)\ifill f:0.7 \esegment

\bsegment \move(-1 0)\lvec(7 0)\lvec(7 13) \move(0 0)\lvec(0
4)\lvec(1 4)\lvec(1 5)\lvec(2 5) \lvec(2 6)\lvec(3 6)\lvec(3
7)\lvec(4 7)\lvec(4 11)\lvec(5 11) \lvec(5 12)\lvec(6 12)\lvec(6
13)\lvec(7 13) \move(4 9)\lvec(7 9) \move(4 10)\lvec(7 10) \move(5
9)\lvec(5 10) \move(6 9)\lvec(6 10) \move(0 3)\lvec(7 3) \move(0
4)\lvec(7 4) \move(1 3)\lvec(1 4) \move(6 3)\lvec(6 4) \esegment

\htext(4.5 9.5){$n$} \htext(6.5 9.5){$n$} \htext(5.5
9.5){$\cdots$} \htext(0.5 3.5){$n$} \htext(6.5 3.5){$n$} \htext(4
3.5){$\cdots$} \htext(5
8.5){$L_{(\omega_{i_k},\omega_{i_{k+1}})}^-$} \htext(3
4.5){$L_{(\omega_{i_k},\omega_{i_{k+1}})}^+$}

\move(8.5 9.5)\avec(6.8 9.5) \htext(9 9.5){$b_R$} \move(2.5
9.5)\avec(4.2 9.5) \htext(2 9.5){$b_L$} \move(2.8 7.5)\avec(4.5
7.5) \htext(2.3 7.5){$b$} \move(1.6 5.5)\avec(3.5 5.5) \htext(1.1
5.5){$b'$} \move(3.5 2)\avec(3.5 3.5) \htext(3.5 1.5){$b'_R$}
\move(1.5 2)\avec(1.5 3.5) \htext(1.5 1.5){$b'_L$}

\end{texdraw}}
\savebox{\tmpfigb} {\begin{texdraw} \fontsize{7}{7}\selectfont
\drawdim em \setunitscale 1.7 \textref h:C v:C

\bsegment \setgray 0.6 \move(4.5 9.5)\lvec(6.5 9.5)\lvec(4.5
7.5)\lvec(4.5 9.5)\ifill f:0.7 \esegment

\bsegment \setgray 0.6 \move(1.5 3.5)\lvec(3.5 3.5)\lvec(3.5
5.5)\lvec(3.5 5.5)\ifill f:0.7 \esegment

\bsegment \move(-1 3)\rlvec(8 0)\rlvec(0 10) \move(0 3)\lvec(0
4)\move(0 4)\lvec(1 4)\lvec(1 5)\lvec(2 5) \lvec(2 6)\lvec(3
6)\lvec(3 7)\lvec(4 7)\lvec(4 11)\lvec(5 11) \lvec(5 12)\lvec(6
12)\lvec(6 13)\lvec(7 13) \move(4 9)\lvec(7 9) \move(4 10)\lvec(7
10) \move(5 9)\lvec(5 10) \move(6 9)\lvec(6 10) \move(0 3)\lvec(7
3) \move(0 4)\lvec(7 4) \move(1 3)\lvec(1 4) \move(6 3)\lvec(6 4)
\move(0 3.5)\rlvec(7 0)\move(4 9.5)\rlvec(3 0) \esegment

\htext(4.5 9.25){$n$}\htext(4.5 9.75){$n$}
 \htext(6.5 9.25){$n$}\htext(6.5 9.75){$n$} \htext(5.5
9.25){$\cdots$}\htext(5.5 9.75){$\cdots$}

\htext(0.5 3.25){$n$}\htext(0.5 3.75){$n$} \htext(6.5
3.25){$n$}\htext(6.5 3.75){$n$} \htext(4 3.25){$\cdots$} \htext(4
3.75){$\cdots$} \htext(5 8.5){$L_{(\omega_{i_t},\la_n)}^-$}
\htext(3 4.5){$L_{(\omega_{i_t},\la_n)}^+$}

\move(8.5 9.5)\avec(6.8 9.5) \htext(9 9.5){$b_R$} \move(2.5
9.5)\avec(4.2 9.5) \htext(2 9.5){$b_L$} \move(2.8 7.5)\avec(4.5
7.5) \htext(2.3 7.5){$b$} \move(1.6 5.5)\avec(3.5 5.5) \htext(1.1
5.5){$b'$} \move(3.5 2)\avec(3.5 3.5) \htext(3.5 1.5){$b'_R$}
\move(1.5 2)\avec(1.5 3.5) \htext(1.5 1.5){$b'_L$}

\end{texdraw}}

\vskip 5mm Now, consider ${\overline
L}_{\omega_{i_k}+\omega_{i_{k+1}}}$ ($k=1,\cdots,t-1$) and
${\overline L}_{\omega_{i_t}+\la_n}$.

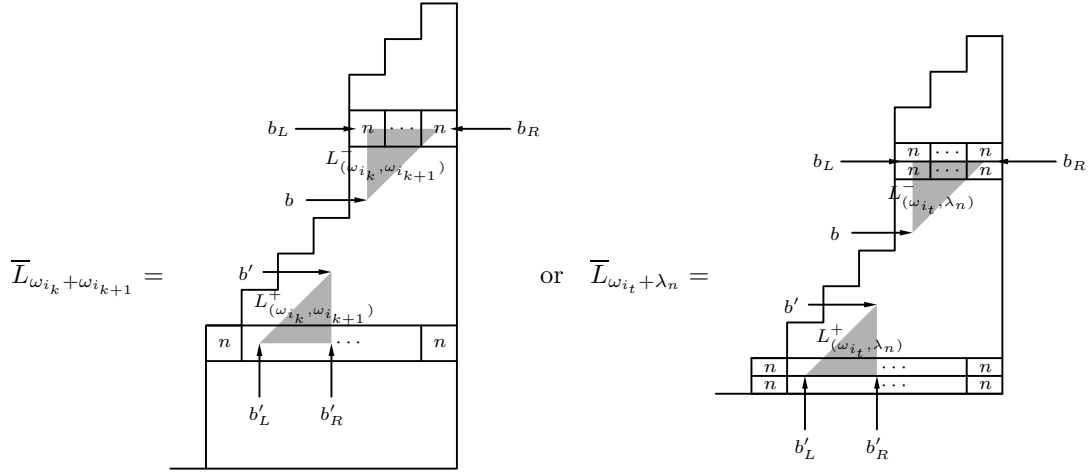
$${\overline L}_{\omega_{i_k}+\omega_{i_{k+1}}}=
\raisebox{-0.4\height}{\usebox{\tmpfiga}}\,\,\,\text{or}\,\,\,\,\,\,
{\overline
L}_{\omega_{i_t}+\la_n}=\raisebox{-0.4\height}{\usebox{\tmpfigb}}$$

\vskip 3mm As we can see from the above picture, ${\overline
L}_{\omega_{i_k}+\omega_{i_{k+1}}}$ and ${\overline
L}_{\omega_{i_t}+\la_{n}}$ contains two $n$-rows above
${\overline H}_{\omega_{i_k}+\omega_{i_{k+1}}}$. Note that there
are $i$-many blocks in the upper $n$-row. Let us denote by $b_{L}$
(resp. $b_{R}$) the left-most (resp. right-most) block in the
upper $n$-row. Then the blocks $b_{L}$, $b_{R}$ and the block $b$
lying in the $(i-1)$-th row below $b_{L}$ form a {\it right
isosceles triangle}. We denote by
$L_{(\omega_{i_k},\omega_{i_{k+1}})}^{-}$ (resp.
$L_{(\omega_{i_t},\la_{n})}^{-}$) the part of ${\overline
L}_{\omega_{i_k}+\omega_{i_{k+1}}}$ (resp. ${\overline
L}_{\omega_{i_t}+\la_{n}}$) constituting this right isosceles
triangle.

Similarly, let $b_{R}'$ be the right-most block in the lower
$n$-row outside the highest weight vector ${\overline
H}_{\omega_{i_k}+\omega_{i_{k+1}}}$ and ${\overline
H}_{\omega_{i_t}+\la_{n}}$ and let $b_{L}'$ be the $n$-block lying
in the $(i-1)$-th column to the left of $b_{R}'$. Then $b_{R}'$,
$b_{L}'$ and the block $b'$ lying in the $(i-1)$-th row above
$b_{R}'$ form another {\it right isosceles triangle}. We denote
by $L_{(\omega_{i_k},\omega_{i_{k+1}})}^{+}$ (resp.
$L_{(\omega_{i_t},\la_{n})}^{+}$) the part of ${\overline
L}_{\omega_{i_k}+\omega_{i_{k+1}}}$ (resp. ${\overline
L}_{\omega_{i_t}+\la_{n}}$) constituting this right isosceles
triangle. Note that $L_{(\omega_{i_k},\omega_{i_{k+1}})}^-$ (resp.
$L_{(\omega_{i_t},\la_{n})}^{-}$) and
$L_{(\omega_{i_k},\omega_{i_{k+1}})}^+$ (resp.
$L_{(\omega_{i_t},\la_{n})}^{+}$) are of the same size with each
base of length $i$. Now, for each $Y \in F(\lambda)$, set
$$
\aligned &Y_{(\omega_{i_k},\omega_{i_{k+1}})}^{-} = Y \cap
L_{(\omega_{i_k},\omega_{i_{k+1}})}^{-}, \quad
Y_{(\omega_{i_k},\omega_{i_{k+1}})}^{+} = Y \cap
L_{(\omega_{i_k},\omega_{i_{k+1}})}^{+}\\
&Y_{(\omega_{i_t},\la_n)}^{-} = Y \cap
L_{(\omega_{i_t},\la_n)}^{-}, \quad Y_{(\omega_{i_t},\la_n)}^{+} =
Y \cap L_{(\omega_{i_t},\la_n)}^{+}
\endaligned
$$
and denote by $|Y_{(\omega_{i_k},\omega_{i_{k+1}})}^{-}|$ (resp.
$|Y_{(\omega_{i_t},\la_n)}^{-}|$) the wall obtained by reflecting
$Y_{(\omega_{i_k},\omega_{i_{k+1}})}^{-}$ (resp.
$Y_{(\omega_{i_t},\la_n)}^{-}$) with respect to the upper $n$-row
and shifting the blocks to the right as much as possible.

\savebox{\tmpfiga}{\begin{texdraw} \fontsize{6}{6}\selectfont
\drawdim em \setunitscale 0.9 \nc{\dtri}{ \bsegment \lvec(-2 0)
\lvec(-2 2)\lvec(0 2)\lvec(0 0)\ifill f:0.75 \esegment }

\nc{\dtrii}{ \bsegment \lvec(-2 0) \lvec(-2 1)\lvec(0 1)\lvec(0
0)\ifill f:0.75 \esegment }

\bsegment \move(16 16)\dtrii \move(14 16)\dtrii \move(12 16)\dtrii
\move(14 14)\dtri \move(12 14)\dtri \move(12 12)\dtri \move(10
4)\dtri \move(10 2)\dtri \move(8 2)\dtri \move(10 1)\dtrii \move(8
1)\dtrii \move(6 1)\dtrii

\setgray 0.6 \move(16 0) \lvec(16 22) \move(16 22)\rlvec(-2
0)\rlvec(0 -22) \move(16 20) \rlvec(-4 0)\rlvec(0 -20)

\move(16 18)\rlvec(-4 0)\move(16 17)\rlvec(-4 0) \move(16
16)\rlvec(-6 0)\rlvec(0 -16)\move(16 14)\rlvec(-6 0)\move(16
12)\rlvec(-6 0)\move(16 10)\rlvec(-8 0)\rlvec(0 -10)

\move(16 8)\rlvec(-10 0)\rlvec(0 -8)\move(16 6)\rlvec(-12
0)\rlvec(0 -6)\move(16 4)\rlvec(-14 0)\rlvec(0 -4) \move(16
2)\rlvec(-16 0)\rlvec(0 -2)\move(16 1)\rlvec(-16 0)\move(-3
0)\lvec(16 0)

\move(16 10)\rlvec(-2 -2)\move(14 10)\rlvec(-2 -2)\move(12
10)\rlvec(-2 -2)\move(10 10)\rlvec(-2 -2)

\move(16 12)\linewd 0.15 \setgray 0 \rlvec(-2 0)\rlvec(0
-2)\rlvec(-4 0)\rlvec(0 -9)\rlvec(-10 0)\rlvec(0 -1)\rlvec(16
0)\rlvec(0 12)\move(12 10)\rlvec(-2 -2)

\htext(15 21){$3$}\htext(15 19){$4$} \htext(13 19){$4$}\htext(15
17.5){$5$} \htext(13 17.5){$5$} \htext(15 16.5){$5$} \htext(13
16.5){$5$} \htext(15 15){$4$} \htext(13 15){$4$} \htext(11
15){$4$} \htext(15 13){$3$} \htext(13 13){$3$} \htext(11 13){$3$}
\htext(15 11){$2$} \htext(13 11){$2$} \htext(11 11){$2$}

\htext(14.5 9.5){$0$}\htext(15.5 8.5){$1$} \htext(13.5
8.5){$0$}\htext(12.5 9.5){$1$} \htext(11.5 8.5){$1$}\htext(10.5
9.5){$0$}\htext(8.5 9.5){$1$}

\htext(15 7){$2$} \htext(13 7){$2$} \htext(11 7){$2$}\htext(9
7){$2$} \htext(7 7){$2$} \htext(15 5){$3$} \htext(13 5){$3$}
\htext(11 5){$3$}\htext(9 5){$3$} \htext(7 5){$3$}\htext(5 5){$3$}
\htext(15 3){$4$} \htext(13 3){$4$} \htext(11 3){$4$}\htext(9
3){$4$} \htext(7 3){$4$}\htext(5 3){$4$}\htext(3 3){$4$}

\htext(15 1.5){$5$} \htext(13 1.5){$5$} \htext(11
1.5){$5$}\htext(9 1.5){$5$} \htext(7 1.5){$5$}\htext(5
1.5){$5$}\htext(3 1.5){$5$}\htext(1 1.5){$5$} \htext(15 0.5){$5$}
\htext(13 0.5){$5$} \htext(11 0.5){$5$}\htext(9 0.5){$5$} \htext(7
0.5){$5$}\htext(5 0.5){$5$}\htext(3 0.5){$5$}\htext(1 0.5){$5$}

\esegment
\end{texdraw}}
\savebox{\tmpfigb}{\begin{texdraw} \fontsize{6}{6}\selectfont
\drawdim em \setunitscale 0.9

\rlvec(0 4)\rlvec(6 0)\rlvec(0 1)\rlvec(-4 0)\rlvec(0 -5)\rlvec(-2
0)\move(0 2)\rlvec(4 0)\rlvec(0 3)

\htext(1 1){$3$}\htext(1 3){$4$}\htext(3 3){$4$}\htext(3
4.5){$5$}\htext(5 4.5){$5$}
\end{texdraw}}
\savebox{\tmpfigc}{\begin{texdraw} \fontsize{6}{6}\selectfont
\drawdim em \setunitscale 0.9

\rlvec(6 0)\rlvec(0 5)\rlvec(-2 0)\rlvec(0 -5)\move(0 0)\rlvec(0
1) \rlvec(6 0)\move(2 0)\rlvec(0 3)\rlvec(4 0)

\htext(1 0.5){$5$}\htext(3 0.5){$5$}\htext(5 0.5){$5$}\htext(3
2){$4$}\htext(5 2){$4$}\htext(5 4){$3$}
\end{texdraw}}
\savebox{\tmpfigd}{\begin{texdraw} \fontsize{6}{6}\selectfont
\drawdim em \setunitscale 0.9

\rlvec(4 0)\rlvec(0 5)\rlvec(-2 0)\rlvec(0 -5)\move(0 0)\rlvec(0
3) \rlvec(4 0)\move(0 1)\rlvec(4 0)

\htext(1 0.5){$5$}\htext(1 2){$4$}\htext(3 0.5){$5$}\htext(3
2){$4$}\htext(3 4){$3$}
\end{texdraw}}

\vskip 3mm
\begin{example}
If $\frak g=B_5$, $\la=\omega_3+\la_5$ and
$$Y= \raisebox{-0.4\height}{\usebox{\tmpfiga}}\in F(\la),$$
then we have
$$Y_{(\omega_3,\lambda_5)}^-=\raisebox{-0.4\height}{\usebox{\tmpfigb}}\,, \quad
Y_{(\omega_3,\lambda_5)}^+=\raisebox{-0.4\height}{\usebox{\tmpfigc}}\quad
\text{and}\quad
|Y_{(\omega_3,\lambda_5)}^-|=\raisebox{-0.4\height}{\usebox{\tmpfigd}}\,
.$$ Here, the shaded parts represent $L_{(\omega_3,\lambda_5)}^-$
and $L_{(\omega_3,\lambda_5)}^+$.

\end{example}

\vskip 5mm

For $a=1, \cdots, n-1$, consider the ${\overline
Y}_{\omega_{i_k}+\omega_{i_{k+1}}}$  and ${\overline
Y}_{\omega_{i_t}+\la_n}$ of $Y\in F(\la)$ having the following
configuration :

\vskip 3mm
\savebox{\tmpfiga}{\begin{texdraw} \fontsize{7.5}{7.5}\selectfont
\drawdim em \setunitscale 0.9 \nc{\dtri}{ \bsegment \lvec(-2 0)
\lvec(-2 2)\lvec(0 2)\lvec(0 0)\ifill f:0.7 \esegment }

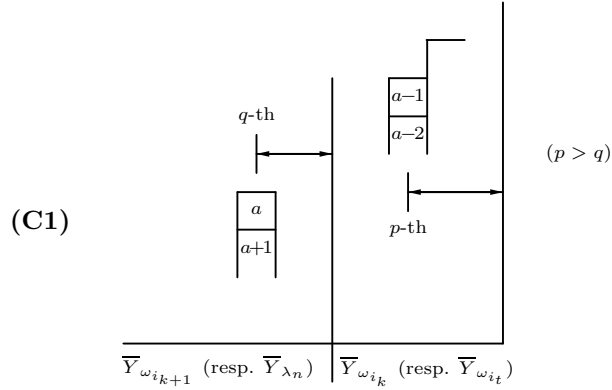
\bsegment \setgray 0 \move(0 0) \lvec(0 18) \move(-2 16 )\rlvec(-2
0)\move(-4 14)\rlvec(-2 0)\move(-4 12)\rlvec(-2 0) \move(-12
8)\rlvec(-2 0) \move(-12 6)\rlvec(-2 0) \move(-4 10)\rlvec(0
6)\move(-6 10)\rlvec(0 4) \move(-12 3.5)\rlvec(0 4.5)\move(-14
3.5)\rlvec(0 4.5) \move(0 0)\rlvec(-20 0) \move(-9 -2)\rlvec(0 16)

\move(-5 8)\ravec(5 0) \move(0 8)\ravec(-5 0) \move(-5 9)\rlvec(0
-2)\htext(-5 6){$p$-th}

\move(-13 10)\ravec(4 0) \move(-9 10)\ravec(-4 0) \move(-13
9)\rlvec(0 2)\htext(-13 12){$q$-th}

\htext(-5 13){$a\!\!-\!\!1$}\htext(-5 11){$a\!\!-\!\!2$}
\htext(-13 7){$a$} \htext(-13 5){$a\!\!+\!\!1$}

\htext(-15 -1.5){${\overline Y_{\omega_{i_{k+1}}}}$ (resp.
${\overline Y_{\la_n}}$)}\htext(-4 -1.5){${\overline
Y_{\omega_{i_k}}}$ (resp. ${\overline Y_{\omega_{i_t}}}$)}
\htext(4 10){($p>q$)}

\esegment
\end{texdraw}}

\begin{center}
{\bf (C1)}\qquad $\raisebox{-0.4\height}{\usebox{\tmpfiga}}$
\end{center}

\savebox{\tmpfiga}{\begin{texdraw} \fontsize{7.5}{7.5}\selectfont
\drawdim em \setunitscale 0.9

\move(0 -0.5)\lvec(0 3)\rlvec(2 0)\rlvec(0 -3.5)\move(0
1.5)\rlvec(2 0)

\htext(1 0.75){$a\!\!-\!\!2$}\htext(1 2.25){$a\!\!-\!\!1$}
\end{texdraw}}
\savebox{\tmpfigb}{\begin{texdraw} \fontsize{7.5}{7.5}\selectfont
\drawdim em \setunitscale 0.9

\move(0 -0.5)\lvec(0 3)\rlvec(2 0)\rlvec(0 -3.5)\move(0
1.5)\rlvec(2 0)

\htext(1 0.75){$a\!\!+\!\!1$}\htext(1 2.25){$a$}
\end{texdraw}}

\vskip 3mm \noindent That is, the top of the $p$-th column of
${\overline Y_{\omega_{i_k}}}$ (resp. ${\overline
Y_{\omega_{i_t}}}$) from the right is
$\raisebox{-0.4\height}{\usebox{\tmpfiga}}$ and the top of the
$q$-th column of ${\overline Y_{\omega_{i_{k+1}}}}$ (resp.
${\overline Y_{\la_n}}$) from the right is
$\raisebox{-0.4\height}{\usebox{\tmpfigb}}$ with $p > q$.

We define $L_{\omega_{i_k}}^{+}(a;p,q)$ (resp.
$L_{\omega_{i_{k+1}}}^{+}(a;p,q)$) to be the right isosceles
triangle formed by $a$-block in the $q$-th column,
$(a+p-q-1)$-block in the $q$-th column and $(a+p-q-1)$-block in
the $(p-1)$-th column in $Y_{\omega_{i_k}}$ (resp.
$Y_{\omega_{i_{k+1}}}$). Then the wall obtained by reflecting
$L_{\omega_{i_k}}^{+}(a;p,q)$ (resp.
$L_{\omega_{i_{k+1}}}^{+}(a;p,q)$) with respect to the $n$-row
will be denoted by $L_{\omega_{i_k}}^{-}(a;p,q)$ (resp.
$L_{\omega_{i_{k+1}}}^{-}(a;p,q)$). The shaded parts in the
following picture represent $L_{\omega_{i_k}}^{\pm}(a;p,q)$  and
$L_{\omega_{i_{k+1}}}^{\pm}(a;p,q)$.

\vskip 3mm
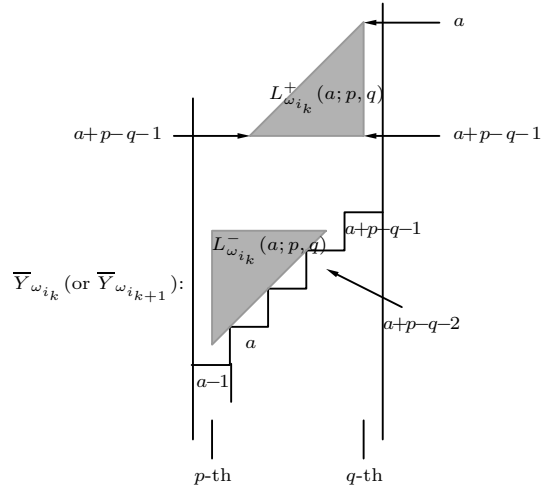
\begin{center}
\begin{texdraw}
\fontsize{7.5}{7.5}\selectfont \drawdim em \setunitscale 0.9

\nc{\dtri}{ \bsegment \lvec(-6 0) \lvec(0 6) \lvec(0 0)\ifill
f:0.7 \esegment }

\nc{\dtrii}{ \bsegment \lvec(-6 0) \lvec(-6 -6) \lvec(0 0)\ifill
f:0.7 \esegment }

\bsegment  \move(-1 16) \dtri \move(-3 11) \dtrii

\setgray 0 \move(0 0) \lvec(0 23)

\move(0 12)\rlvec(-2 0)\rlvec(0 -2)\rlvec(-2 0)\rlvec(0
-2)\rlvec(-2 0)\rlvec(0 -2)\rlvec(-2 0) \rlvec(0 -2) \rlvec(-2
0)\move(-8 4)\rlvec(0 -2) \move(-10 0)\rlvec(0 18)\move(-9
-1)\rlvec(0 2) \move(-1 -1)\rlvec(0 2)

\htext(0 11){$a\!\!+\!\!p\!\!-\!\!q\!\!-\!\!1$}

\move(1 7)\avec(-3 9)\htext(2
6){$a\!\!+\!\!p\!\!-\!\!q\!\!-\!\!2$}


\htext(-7 5){$a$} \htext(-9 3){$a\!\!-\!\!1$}

\htext(-1 -2){$q$-th}\htext(-9 -2){$p$-th} \htext(-15
8){${\overline Y_{\omega_{i_k}}}$(or ${\overline
Y_{\omega_{i_{k+1}}} }$):} \htext(-3
18){$L_{\omega_{i_k}}^{+}(a;p,q)$}\htext(-6
10){$L_{\omega_{i_k}}^{-}(a;p,q)$}

\move(3 22)\avec(-1 22) \move(3 16)\avec(-1 16) \move(-11
16)\avec(-7 16) \htext(4 22){$a$} \htext(6
16){$a\!\!+\!p\!\!-\!q\!\!-\!1$} \htext(-14
16){$a\!\!+\!p\!\!-\!q\!\!-\!1$}

\setgray 0.6 \move(-1 16)\rlvec(-6 0)\rlvec(6 6)\rlvec(0 -6)
\move(-3 11)\rlvec(-6 0)\rlvec(0 -6)\rlvec(6 6)
\esegment
\end{texdraw}
\end{center}
Now, we also define $L^{\pm}_{\omega_{i_t}}(a,p,q)$ in a similar
way, and for each $Y\in F(\la)$, set
\begin{equation}
\begin{aligned}
& Y_{\omega_{i_k}}^{\pm}(a;p,q) = L_{\omega_{i_k}}^{\pm}(a;p,q)
\cap Y, \quad Y_{\omega_{i_{k+1}}}^{\pm}(a;p,q) =
L_{\omega_{i_{k+1}}}^{\pm}(a;p,q) \cap Y,\\
&Y_{\omega_{i_t}}^{\pm}(a;p,q) = L_{\omega_{i_t}}^{\pm}(a;p,q)
\cap Y,
\end{aligned}
\end{equation}
and let $|Y_{\omega_{i_k}}^{-}(a;p,q)|$ be the wall obtained by
reflecting $Y_{\omega_{i_k}}^{-}(a;p,q)$ with respect to the
$n$-row and shifting the blocks to the right as much as possible.


%
\savebox{\tmpfiga}{\begin{texdraw} \fontsize{6}{6}\selectfont
\drawdim em \setunitscale 0.9 \nc{\dtri}{ \bsegment \lvec(-2 0)
\lvec(-2 2)\lvec(0 2)\lvec(0 0)\ifill f:0.8 \esegment }

\nc{\dtrii}{ \bsegment \lvec(-2 0) \lvec(-2 1)\lvec(0 1)\lvec(0
0)\ifill f:0.8 \esegment }

\bsegment \move(14 30)\dtri \move(14 28)\dtri \move(12 28)\dtri
\move(12 20)\dtri \move(10 20)\dtri \move(10 18)\dtri \move(8
14)\dtri \move(8 12)\dtri \move(6 12)\dtri \move(6 4)\dtri \move(4
4)\dtri \move(4 2)\dtri

\setgray 0.6 \move(-3 0)\lvec(14 0) \lvec(14 32) \move(14
30)\rlvec(-2 0)\rlvec(0 -30) \move(14 28) \rlvec(-2 0)\move(14
26)\rlvec(-4 0)\rlvec(0 -26) \move(14 24)\rlvec(-4 0) \move(14
22)\rlvec(-4 0) \move(14 20)\rlvec(-4 0) \move(14 18)\rlvec(-6
0)\rlvec(0 -18) \move(14 16)\rlvec(-8 0)\rlvec(0 -16)\move(14
14)\rlvec(-10 0)\rlvec(0 -14)\move(14 12)\rlvec(-10 0)\move(14
10)\rlvec(-10 0) \move(14 8)\rlvec(-10 0)\move(14 6)\rlvec(-10
0)\move(14 4)\rlvec(-12 0)\rlvec(0 -4) \move(14 2)\rlvec(-14
0)\rlvec(0 -2)

\move(14 18)\rlvec(-2 -2) \move(12 18)\rlvec(-2 -2)\move(10
18)\rlvec(-2 -2)\move(14 2)\rlvec(-2 -2) \move(12 2)\rlvec(-2
-2)\move(10 2)\rlvec(-2 -2) \move(8 2)\rlvec(-2 -2)\move(6
2)\rlvec(-2 -2) \move(4 2)\rlvec(-2 -2)\move(2 2)\rlvec(-2 -2)

\move(14 20)\linewd 0.15 \setgray 0 \rlvec(-2 0)\rlvec(0
-2)\rlvec(-2 0)\rlvec(-2 -2)\rlvec(0 -10)\rlvec(-2 0)\rlvec(0
-2)\rlvec(-2 0)\rlvec(0 -2)\rlvec(-2 0)\rlvec(-2 -2)\lvec(14
0)\lvec(14 20)

\linewd 0.15 \setgray 0.3 \move(8 14)\rlvec(0 14)\move(14
16)\lvec(14 27)\move(2 -1)\rlvec(0 17)

\htext(13 29){$3$} \htext(13 27){$4$} \htext(13 25){$5$} \htext(11
25){$5$}\htext(13 23){$4$} \htext(11 23){$4$} \htext(13 21){$3$}
\htext(11 21){$3$}\htext(13 19){$2$} \htext(11 19){$2$}\htext(13.5
16.5){$0$}\htext(12.5 17.5){$1$} \htext(11.5 16.5){$1$}\htext(10.5
17.5){$0$} \htext(9.5 16.5){$0$}\htext(8.5 17.5){$1$}\htext(13
15){$2$} \htext(11 15){$2$} \htext(9 15){$2$}\htext(7 15){$2$}
\htext(13 13){$3$} \htext(11 13){$3$} \htext(9 13){$3$} \htext(7
13){$3$}\htext(5 13){$3$} \htext(13 11){$4$} \htext(11
11){$4$}\htext(9 11){$4$} \htext(7 11){$4$}\htext(5
11){$4$}\htext(13 9){$5$} \htext(11 9){$5$}\htext(9 9){$5$}
\htext(7 9){$5$}\htext(5 9){$5$} \htext(13 7){$4$} \htext(11
7){$4$}\htext(9 7){$4$} \htext(7 7){$4$}\htext(5 7){$4$}\htext(13
5){$3$}\htext(11 5){$3$} \htext(9 5){$3$}\htext(7 5){$3$} \htext(5
5){$3$} \htext(13 3){$2$}\htext(11 3){$2$} \htext(9
3){$2$}\htext(7 3){$2$} \htext(5 3){$2$}\htext(3 3){$2$}

\htext(13.5 0.5){$0$}\htext(12.5 1.5){$1$} \htext(11.5
0.5){$1$}\htext(10.5 1.5){$0$} \htext(9.5 0.5){$0$}\htext(8.5
1.5){$1$}\htext(7.5 0.5){$1$}\htext(6.5 1.5){$0$} \htext(5.5
0.5){$0$}\htext(4.5 1.5){$1$}\htext(3.5 0.5){$1$}\htext(2.5
1.5){$0$}\htext(1.5 0.5){$0$}\htext(0.5 1.5){$1$}

\esegment
\end{texdraw}}
\savebox{\tmpfigb}{\begin{texdraw} \fontsize{6}{6}\selectfont
\drawdim em \setunitscale 0.9 \nc{\dtri}{ \bsegment \lvec(-2 0)
\lvec(-2 2)\lvec(0 2)\lvec(0 0)\ifill f:0.8 \esegment }

\nc{\dtrii}{ \bsegment \lvec(-2 0) \lvec(-2 1)\lvec(0 1)\lvec(0
0)\ifill f:0.8 \esegment }

\lvec(2 0) \rlvec(0 2)\rlvec(-2 0)\rlvec(0 -2)

\htext(1 1){$3$}
\end{texdraw}}
\savebox{\tmpfigc}{\begin{texdraw} \fontsize{6}{6}\selectfont
\drawdim em \setunitscale 0.9 \nc{\dtri}{ \bsegment \lvec(-2 0)
\lvec(-2 2)\lvec(0 2)\lvec(0 0)\ifill f:0.8 \esegment }

\nc{\dtrii}{ \bsegment \lvec(-2 0) \lvec(-2 1)\lvec(0 1)\lvec(0
0)\ifill f:0.8 \esegment }

\lvec(4 0) \rlvec(0 4)\rlvec(-2 0)\rlvec(0 -4)\move(0 0)\rlvec(0
2)\rlvec(4 0)

\htext(1 1){$3$} \htext(3 1){$3$} \htext(3 3){$2$}
\end{texdraw}}
\savebox{\tmpfigd}{\begin{texdraw} \fontsize{6}{6}\selectfont
\drawdim em \setunitscale 0.9 \nc{\dtri}{ \bsegment \lvec(-2 0)
\lvec(-2 2)\lvec(0 2)\lvec(0 0)\ifill f:0.8 \esegment }

\nc{\dtrii}{ \bsegment \lvec(-2 0) \lvec(-2 1)\lvec(0 1)\lvec(0
0)\ifill f:0.8 \esegment }

\lvec(2 0) \rlvec(0 2)\rlvec(-2 0)\rlvec(0 -2)\move(2 2)\rlvec(2
0) \rlvec(0 2)\rlvec(-2 0)\rlvec(0 -2)

\htext(1 1){$2$} \htext(3 3){$3$}
\end{texdraw}}
\savebox{\tmpfige}{\begin{texdraw} \fontsize{6}{6}\selectfont
\drawdim em \setunitscale 0.9 \nc{\dtri}{ \bsegment \lvec(-2 0)
\lvec(-2 2)\lvec(0 2)\lvec(0 0)\ifill f:0.8 \esegment }

\nc{\dtrii}{ \bsegment \lvec(-2 0) \lvec(-2 1)\lvec(0 1)\lvec(0
0)\ifill f:0.8 \esegment }

\lvec(2 0) \rlvec(0 4)\rlvec(-2 0)\rlvec(0 -4)\move(2 2)\rlvec(-2
0)

\htext(1 1){$3$} \htext(1 3){$2$}
\end{texdraw}}

\vskip 3mm
\begin{example}
If $\frak g=C_5$, $\la=\omega_3+\omega_4$ and
$$Y=\raisebox{-0.4\height}{\usebox{\tmpfiga}}\in F(\la),$$
then we have
$$Y_{\omega_3}^\pm(2;3,1)=|Y_{\omega_3}^-(2;3,1)|=\raisebox{-0.4\height}{\usebox{\tmpfigb}}\,,$$
$$Y_{\omega_4}^+(2;3,1)=\raisebox{-0.4\height}{\usebox{\tmpfigc}}\,,
\quad
Y_{\omega_4}^-(2;3,1)=\raisebox{-0.4\height}{\usebox{\tmpfigd}}$$
and
$$ |Y_{\omega_4}^-(2;3,1)|=\raisebox{-0.4\height}{\usebox{\tmpfige}}\, .$$
Here, the shaded parts represent $L_{\omega_i}^{\pm}(2;3,1)$ for
$i=3,4$.

\end{example}

\vskip 3mm Now, we are ready to give an explicit description of
the crystal graph $B(\la)$ over $\frak g = C_n$ and $B_n$.

\vskip 3mm
\begin{thm} \label{thm:Cn}
Let $\la \in P^{+}$ be a dominant integral weight for $\frak g =
C_n$, and write
$$\la = \omega_{i_1} + \cdots + \omega_{i_t}
\qquad (1 \le i_1 \le \cdots \le i_t \le n).$$ We define $Y(\la)$
to be the set of all reduced proper Young walls in $F(\la)$
satisfying the following conditions {\rm :}

{\bf (Y1)} For each $k = 1, \cdots, t$, we have
$Y_{\omega_{i_k}}^{+} \subset |Y_{\omega_{i_k}}^{-}|$.

{\bf (Y2)} For each $k=1, \cdots, t-1$, we have $Y^{\omega_{i_k}}
\subset Y^{\omega_{i_{k+1}}}$ in ${\overline
Y}_{\omega_{i_k}+\omega_{i_{k+1}}}$.

{\bf (Y3)} For each $k=1, \cdots, t-1$, we have
$|Y_{(\omega_{i_k},\omega_{i_{k+1}})}^{-}| \subset
Y_{(\omega_{i_k},\omega_{i_{k+1}})}^{+}$.

{\bf (Y4)} For each $k=1, \cdots, t-1$, if ${\overline
Y_{\omega_{i_k}+\omega_{i_{k+1}}}}$ satisfies $(\bf{C1})$, then we
have
$$Y_{\omega_{i_k}}^{+}(a;p,
q) \subset |Y_{\omega_{i_k}}^{-}(a;p,q)|, \quad
Y_{\omega_{i_{k+1}}}^{+}(a;p,q) \subset
|Y_{\omega_{i_{k+1}}}^{-}(a;p,q)|.$$

Then there is an isomorphism of $U_q(C_n)$-crystals
\begin{equation}
Y(\lambda) \stackrel{\sim} \longrightarrow B(\la) \quad
\text{given by} \ \ H_{\la} \longmapsto u_{\la},
\end{equation}
where $u_\la$ is the highest weight vector in $B(\la)$.

\end{thm}

\vskip 3mm

\begin{thm} \label{thm:Bn}
Let $\la \in P^{+}$ be a dominant integral weight for $\frak g =
B_n$, and write
$$\aligned
&\la = \omega_{i_1} + \cdots + \omega_{i_t} \qquad (1 \le i_1 \le
\cdots \le i_t \le n)\quad \text{or}\\
&\la = \omega_{i_1} + \cdots + \omega_{i_t} + \la_n \qquad (1 \le
i_1 \le \cdots \le i_t \le n). \endaligned$$

We define $Y(\la)$ to be the set of all reduced proper Young walls
in $F(\la)$ satisfying the following conditions
{\rm :}

{\bf (Y1)} For each $k = 1, \cdots, t$, we have
$Y_{\omega_{i_k}}^{+} \subset |Y_{\omega_{i_k}}^{-}|$.

{\bf (Y2)} For each $k=1, \cdots, t-1$, we have
$$\text{$Y^{\omega_{i_k}}
\subset Y^{\omega_{i_{k+1}}}$ in ${\overline
Y}_{\omega_{i_k}+\omega_{i_{k+1}}}$ and \,\,\,$Y^{\omega_{i_t}}
\subset Y^{\la_n}$ in ${\overline Y}_{\omega_{i_t}+\la_n}$}.$$

{\bf (Y3)} For each $k=1, \cdots, t-1$, we have
$$|Y_{(\omega_{i_k},\omega_{i_{k+1}})}^{-}| \subset
Y_{(\omega_{i_k},\omega_{i_{k+1}})}^{+}, \quad
|Y_{(\omega_{i_t},\la_n)}^{-}| \subset
Y_{(\omega_{i_t},\la_n)}^{+}.$$

{\bf (Y4)} For each $k=1, \cdots, t-1$, if ${\overline
Y_{\omega_{i_k}+\omega_{i_{k+1}}}}$ or ${\overline
Y_{\omega_{i_t}+\la_n}}$ satisfies $(\bf{C1})$, then we have
$$\aligned
&Y_{\omega_{i_k}}^{+}(a;p,q) \subset
|Y_{\omega_{i_k}}^{-}(a;p,q)|, \quad
Y_{\omega_{i_{k+1}}}^{+}(a;p,q) \subset
|Y_{\omega_{i_{k+1}}}^{-}(a;p,q)|,\\
&Y_{\omega_{i_t}}^{+}(a;p,q) \subset
|Y_{\omega_{i_t}}^{-}(a;p,q)|.
\endaligned$$

Then there is an isomorphism of crystal graphs for
$U_q(B_n)$-modules
\begin{equation}
Y(\lambda) \stackrel{\sim} \longrightarrow B(\la) \quad
\text{given by} \ \ H_{\la} \longmapsto u_{\la},
\end{equation}
where $u_\la$ is the highest weight vector in $B(\la)$.

\end{thm}

\vskip 2mm
\begin{remark}  If $\la=\la_n$, then $Y(\la_n)=F(\la_n)$, the set of all
reduced proper Young walls lying between $H_{\la_n}$ and
$L_{\la_n}$.
\end{remark}

\vskip 3mm
\begin{example}
Let $\frak g = C_3$ and $\la = \omega_2 + \omega_3$. Then in the
following picture, the first Young wall belongs to $Y(\la)$ but
the second one and the third one don't. The second one does not
satisfy {\bf (Y1)} and the third one does not satisfy {\bf (Y2)}.

\vskip 3mm

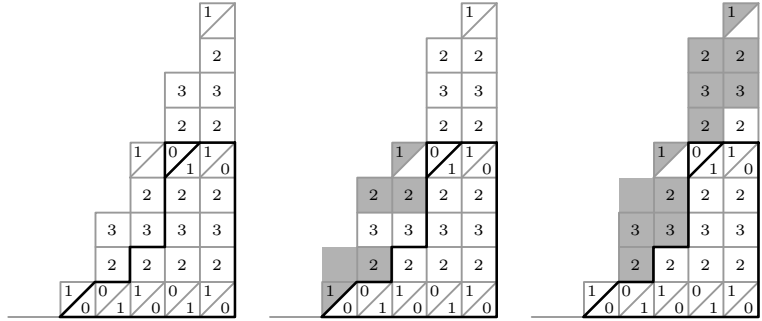
\begin{center}
\begin{texdraw}
\fontsize{6}{6}\selectfont \drawdim em \setunitscale 0.9
\nc{\dtri}{ \bsegment \lvec(-2 0) \lvec(-2 2)\lvec(0 2)\lvec(0
0)\ifill f:0.7 \esegment }

\nc{\dtrii}{ \bsegment \lvec(-2 0) \lvec(-2 1)\lvec(0 1)\lvec(0
0)\ifill f:0.7 \esegment }

\nc{\dtriii}{ \bsegment \lvec(-2 0) \lvec(-2 -2)\lvec(0 0)\ifill
f:0.7 \esegment }

\bsegment

\setgray 0.6 \move(-3 0) \lvec(10 0)\lvec(10 18) \move(10
18)\rlvec(-2 0)\rlvec(0 -18)\move(10 16)\rlvec(-2 0)\move(10
14)\rlvec(-4 0)\rlvec(0 -14)\move(10 12)\rlvec(-4 0)\move(10
10)\rlvec(-6 0)\rlvec(0 -10)\move(10 8)\rlvec(-6 0)\move(10
6)\rlvec(-8 0)\rlvec(0 -6) \move(10 4) \rlvec(-8 0)\move(10
2)\rlvec(-10 0) \rlvec(0 -2)

\move(10 18)\rlvec(-2 -2) \move(10 10)\rlvec(-2 -2)\move(8
10)\rlvec(-2 -2)\move(6 10)\rlvec(-2 -2) \move(10 2)\rlvec(-2
-2)\move(8 2)\rlvec(-2 -2)\move(6 2)\rlvec(-2 -2)\move(4
2)\rlvec(-2 -2) \move(2 2)\rlvec(-2 -2)

\move(10 10)\linewd 0.15 \setgray 0 \rlvec(-4 0)\rlvec(0
-6)\rlvec(-2 0) \rlvec(0 -2)\rlvec(-2 0)\rlvec(-2 -2)\lvec(10
0)\lvec(10 10)\move(8 10)\rlvec(-2 -2)

\htext(8.5 17.5){$1$}\htext(9 15){$2$}\htext(9 13){$3$} \htext(7
13){$3$} \htext(9 11){$2$} \htext(7 11){$2$} \htext(9.5 8.5){$0$}
\htext(8.5 9.5){$1$} \htext(7.5 8.5){$1$} \htext(6.5 9.5){$0$}
\htext(4.5 9.5){$1$} \htext(9 7){$2$} \htext(7 7){$2$} \htext(5
7){$2$} \htext(9 5){$3$} \htext(7 5){$3$} \htext(5 5){$3$}
\htext(3 5){$3$}\htext(9 3){$2$} \htext(7 3){$2$} \htext(5 3){$2$}
\htext(3 3){$2$}

\htext(9.5 0.5){$0$}\htext(8.5 1.5){$1$} \htext(7.5
0.5){$1$}\htext(6.5 1.5){$0$} \htext(5.5 0.5){$0$}\htext(4.5
1.5){$1$}\htext(3.5 0.5){$1$} \htext(2.5 1.5){$0$} \htext(1.5
0.5){$0$} \htext(0.5 1.5){$1$}

\esegment

\move(15 0) \bsegment \move(4 6)\dtri \move(6 6)\dtri \move(6 10
)\dtriii

\move(4 2)\dtri \move(2 2)\dtri \move(2 2)\dtriii

\setgray 0.6 \move(-3 0) \lvec(10 0)\lvec(10 18) \move(10
18)\rlvec(-2 0)\rlvec(0 -18)\move(10 16)\rlvec(-4 0)\rlvec(0 -16)
\move(10 14)\rlvec(-4 0)\move(10 12)\rlvec(-4 0)\move(10
10)\rlvec(-6 0)\rlvec(0 -10)\move(10 8)\rlvec(-8 0)\rlvec(0 -8)
\move(10 6)\rlvec(-8 0) \move(10 4) \rlvec(-8 0)\move(10
2)\rlvec(-10 0) \rlvec(0 -2)

\move(10 18)\rlvec(-2 -2) \move(10 10)\rlvec(-2 -2)\move(8
10)\rlvec(-2 -2)\move(6 10)\rlvec(-2 -2) \move(10 2)\rlvec(-2
-2)\move(8 2)\rlvec(-2 -2)\move(6 2)\rlvec(-2 -2)\move(4
2)\rlvec(-2 -2) \move(2 2)\rlvec(-2 -2)

\move(10 10)\linewd 0.15 \setgray 0 \rlvec(-4 0)\rlvec(0
-6)\rlvec(-2 0) \rlvec(0 -2)\rlvec(-2 0)\rlvec(-2 -2)\lvec(10
0)\lvec(10 10)\move(8 10)\rlvec(-2 -2)

\htext(8.5 17.5){$1$}\htext(9 15){$2$}\htext(7 15){$2$}\htext(9
13){$3$} \htext(7 13){$3$} \htext(9 11){$2$} \htext(7 11){$2$}
\htext(9.5 8.5){$0$} \htext(8.5 9.5){$1$} \htext(7.5 8.5){$1$}
\htext(6.5 9.5){$0$} \htext(4.5 9.5){$1$} \htext(9 7){$2$}
\htext(7 7){$2$} \htext(5 7){$2$}\htext(3 7){$2$} \htext(9 5){$3$}
\htext(7 5){$3$} \htext(5 5){$3$} \htext(3 5){$3$}\htext(9 3){$2$}
\htext(7 3){$2$} \htext(5 3){$2$} \htext(3 3){$2$}

\htext(9.5 0.5){$0$}\htext(8.5 1.5){$1$} \htext(7.5
0.5){$1$}\htext(6.5 1.5){$0$} \htext(5.5 0.5){$0$}\htext(4.5
1.5){$1$}\htext(3.5 0.5){$1$} \htext(2.5 1.5){$0$} \htext(1.5
0.5){$0$} \htext(0.5 1.5){$1$}

\esegment

\move(30 0) \bsegment \move(8 10)\dtri \move(8 14)\dtri \move(8
12)\dtri  \move(10 12)\dtri \move(10 14)\dtri  \move(10 18)\dtriii

\move(4 2)\dtri \move(4 4)\dtri \move(4 6)\dtri \move(6 4)\dtri
\move(6 6)\dtri \move(6 10)\dtriii

\setgray 0.6 \move(-3 0) \lvec(10 0)\lvec(10 18) \move(10
18)\rlvec(-2 0)\rlvec(0 -18)\move(10 16)\rlvec(-4 0)\rlvec(0 -16)
\move(10 14)\rlvec(-4 0)\move(10 12)\rlvec(-4 0)\move(10
10)\rlvec(-6 0)\rlvec(0 -10)\move(10 8)\rlvec(-6 0)\move(10
6)\rlvec(-8 0)\rlvec(0 -6) \move(10 4) \rlvec(-8 0)\move(10
2)\rlvec(-10 0) \rlvec(0 -2)

\move(10 18)\rlvec(-2 -2) \move(10 10)\rlvec(-2 -2)\move(8
10)\rlvec(-2 -2)\move(6 10)\rlvec(-2 -2) \move(10 2)\rlvec(-2
-2)\move(8 2)\rlvec(-2 -2)\move(6 2)\rlvec(-2 -2)\move(4
2)\rlvec(-2 -2) \move(2 2)\rlvec(-2 -2)

\move(10 10)\linewd 0.15 \setgray 0 \rlvec(-4 0)\rlvec(0
-6)\rlvec(-2 0) \rlvec(0 -2)\rlvec(-2 0)\rlvec(-2 -2)\lvec(10
0)\lvec(10 10)\move(8 10)\rlvec(-2 -2)

\htext(8.5 17.5){$1$}\htext(9 15){$2$}\htext(7 15){$2$}\htext(9
13){$3$} \htext(7 13){$3$} \htext(9 11){$2$} \htext(7 11){$2$}
\htext(9.5 8.5){$0$} \htext(8.5 9.5){$1$} \htext(7.5 8.5){$1$}
\htext(6.5 9.5){$0$} \htext(4.5 9.5){$1$} \htext(9 7){$2$}
\htext(7 7){$2$} \htext(5 7){$2$} \htext(9 5){$3$} \htext(7
5){$3$} \htext(5 5){$3$} \htext(3 5){$3$}\htext(9 3){$2$} \htext(7
3){$2$} \htext(5 3){$2$} \htext(3 3){$2$}

\htext(9.5 0.5){$0$}\htext(8.5 1.5){$1$} \htext(7.5
0.5){$1$}\htext(6.5 1.5){$0$} \htext(5.5 0.5){$0$}\htext(4.5
1.5){$1$}\htext(3.5 0.5){$1$} \htext(2.5 1.5){$0$} \htext(1.5
0.5){$0$} \htext(0.5 1.5){$1$}

\esegment

\end{texdraw}
\end{center}

\noindent Here, in the second Young wall, the shaded parts
represent $L_{\omega_3}^{\pm}$, and in the third Young wall, the
shaded parts represent $L^{\omega_2}$ and $L^{\omega_3}$.

\end{example}

\vskip 3mm
\begin{example}
Let $\frak g=B_4$ and $\la = \omega_3 + \la_4$. Then, in the
following picture, the first Young wall belongs to $Y(\la)$, but
the other ones don't. They do not satisfy the conditions {\bf
(Y1)}, {\bf (Y3)} and {\bf (Y4)}, respectively.

\vskip 3mm

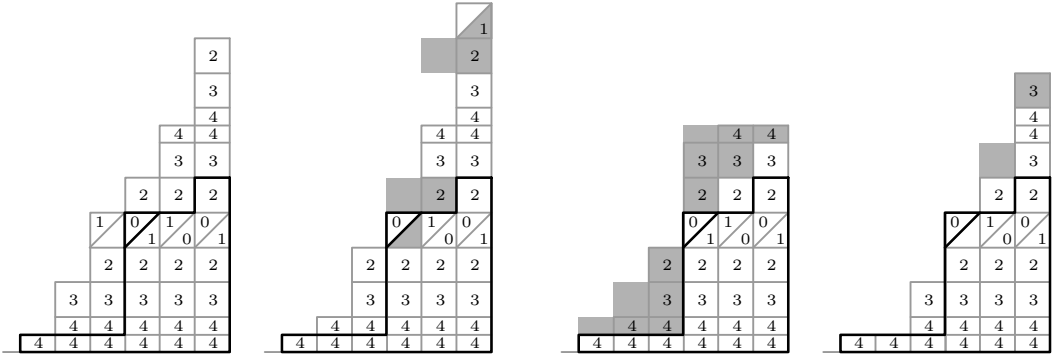
\begin{center}
\begin{texdraw}
\fontsize{6}{6}\selectfont \drawdim em \setunitscale 0.9
\nc{\dtri}{ \bsegment \lvec(-2 0) \lvec(-2 2)\lvec(0 2)\lvec(0
0)\ifill f:0.7 \esegment }

\nc{\dtrii}{ \bsegment \lvec(-2 0) \lvec(-2 1)\lvec(0 1)\lvec(0
0)\ifill f:0.7 \esegment }

\nc{\dtriii}{ \bsegment \lvec(-2 0) \lvec(0 2)\lvec(0 0)\ifill
f:0.7 \esegment }

\bsegment

\setgray 0.6 \move(-1 0) \lvec(12 0)\lvec(12 18) \move(12
18)\rlvec(-2 0)\rlvec(0 -18)\move(12 16)\rlvec(-2 0)\move(12
14)\rlvec(-2 0)\move(12 13)\rlvec(-4 0)\rlvec(0 -13)\move(12
12)\rlvec(-4 0)\move(12 10)\rlvec(-6 0)\rlvec(0 -10)\move(12
8)\rlvec(-8 0)\rlvec(0 -8)\move(12 6)\rlvec(-8 0)\move(12
4)\rlvec(-10 0)\rlvec(0 -4)\move(12 2)\rlvec(-10 0)\move(12
1)\rlvec(-12 0)\rlvec(0 -1)

\move(12 8)\rlvec(-2 -2) \move(10 8)\rlvec(-2 -2)\move(8
8)\rlvec(-2 -2) \move(6 8)\rlvec(-2 -2)

\move(12 10)\linewd 0.15 \setgray 0 \rlvec(-2 0)\rlvec(0 -2)
\rlvec(-4 0) \rlvec(0 -7)\rlvec(-6 0)\rlvec(0 -1)\lvec(12
0)\lvec(12 10)\move(8 8)\rlvec(-2 -2)

\htext(11 17){$2$}\htext(11 15){$3$}\htext(11 13.5){$4$} \htext(11
12.5){$4$}\htext(9 12.5){$4$} \htext(11 11){$3$}\htext(9 11){$3$}
\htext(11 9){$2$}\htext(9 9){$2$}\htext(7 9){$2$} \htext(11.5
6.5){$1$} \htext(10.5 7.5){$0$} \htext(9.5 6.5){$0$} \htext(8.5
7.5){$1$} \htext(7.5 6.5){$1$} \htext(6.5 7.5){$0$} \htext(4.5
7.5){$1$} \htext(11 5){$2$} \htext(9 5){$2$} \htext(7 5){$2$}
\htext(5 5){$2$}\htext(11 3){$3$} \htext(9 3){$3$} \htext(7
3){$3$} \htext(5 3){$3$}\htext(3 3){$3$} \htext(11 1.5){$4$}
\htext(9 1.5){$4$} \htext(7 1.5){$4$} \htext(5 1.5){$4$}\htext(3
1.5){$4$}\htext(11 0.5){$4$} \htext(9 0.5){$4$} \htext(7 0.5){$4$}
\htext(5 0.5){$4$}\htext(3 0.5){$4$}\htext(1 0.5){$4$}

\esegment

\move(15 0) \bsegment \move(12 16)\dtri \move(10 16)\dtri \move(12
18)\dtriii

\move(10 8)\dtri\move(8 8)\dtri\move(8 6)\dtriii

\setgray 0.6 \move(-1 0) \lvec(12 0)\lvec(12 20)\move(12
18)\rlvec(0 2)\rlvec(-2 0)\rlvec(0 -2) \move(12 20)\rlvec(-2
-2)\move(12 18)\rlvec(-2 0)\rlvec(0 -18)\move(12 16)\rlvec(-2
0)\move(12 14)\rlvec(-2 0)\move(12 13)\rlvec(-4 0)\rlvec(0
-13)\move(12 12)\rlvec(-4 0)\move(12 10)\rlvec(-4 0)\move(12
8)\rlvec(-6 0)\rlvec(0 -8)\move(12 6)\rlvec(-8 0)\rlvec(0
-6)\move(12 4)\rlvec(-8 0)\move(12 2)\rlvec(-10 0)\rlvec(0
-2)\move(12 1)\rlvec(-12 0)\rlvec(0 -1)

\move(12 8)\rlvec(-2 -2) \move(10 8)\rlvec(-2 -2)\move(8
8)\rlvec(-2 -2)

\move(12 10)\linewd 0.15 \setgray 0 \rlvec(-2 0)\rlvec(0 -2)
\rlvec(-4 0) \rlvec(0 -7)\rlvec(-6 0)\rlvec(0 -1)\lvec(12
0)\lvec(12 10)\move(8 8)\rlvec(-2 -2)

\htext(11.5 18.5){$1$} \htext(11 17){$2$}\htext(11
15){$3$}\htext(11 13.5){$4$} \htext(11 12.5){$4$}\htext(9
12.5){$4$} \htext(11 11){$3$}\htext(9 11){$3$} \htext(11
9){$2$}\htext(9 9){$2$}\htext(11.5 6.5){$1$} \htext(10.5 7.5){$0$}
\htext(9.5 6.5){$0$} \htext(8.5 7.5){$1$} \htext(6.5 7.5){$0$}
\htext(11 5){$2$} \htext(9 5){$2$} \htext(7 5){$2$} \htext(5
5){$2$}\htext(11 3){$3$} \htext(9 3){$3$} \htext(7 3){$3$}
\htext(5 3){$3$}\htext(11 1.5){$4$} \htext(9 1.5){$4$} \htext(7
1.5){$4$} \htext(5 1.5){$4$}\htext(3 1.5){$4$}\htext(11 0.5){$4$}
\htext(9 0.5){$4$} \htext(7 0.5){$4$} \htext(5 0.5){$4$}\htext(3
0.5){$4$}\htext(1 0.5){$4$}

\esegment

\move(32 0) \bsegment \move(12 12)\dtrii \move(10 12)\dtrii
\move(8 12)\dtrii \move(10 10)\dtri \move(8 8)\dtri \move(8
10)\dtri

\move(6 1)\dtrii \move(4 1)\dtrii \move(2 1)\dtrii \move(6 2)\dtri
\move(4 2)\dtri \move(6 4)\dtri

\setgray 0.6 \move(-1 0) \lvec(12 0)\lvec(12 13) \move(12
13)\rlvec(-4 0)\rlvec(0 -13)\move(12 12)\rlvec(-6 0)\rlvec(0 -12)
\move(12 10)\rlvec(-6 0)\move(12 8)\rlvec(-6 0)\move(12
6)\rlvec(-8 0)\rlvec(0 -6)\move(12 4)\rlvec(-8 0)\move(12
2)\rlvec(-10 0)\rlvec(0 -2)\move(12 1)\rlvec(-12 0)\rlvec(0
-1)\move(10 0)\rlvec(0 13)

\move(12 8)\rlvec(-2 -2) \move(10 8)\rlvec(-2 -2)\move(8
8)\rlvec(-2 -2)

\move(12 10)\linewd 0.15 \setgray 0 \rlvec(-2 0)\rlvec(0 -2)
\rlvec(-4 0) \rlvec(0 -7)\rlvec(-6 0)\rlvec(0 -1)\lvec(12
0)\lvec(12 10)\move(8 8)\rlvec(-2 -2)

\htext(11 12.5){$4$}\htext(9 12.5){$4$} \htext(11 11){$3$}\htext(9
11){$3$} \htext(7 11){$3$} \htext(11 9){$2$}\htext(9
9){$2$}\htext(7 9){$2$} \htext(11.5 6.5){$1$} \htext(10.5
7.5){$0$} \htext(9.5 6.5){$0$} \htext(8.5 7.5){$1$} \htext(7.5
6.5){$1$} \htext(6.5 7.5){$0$} \htext(11 5){$2$} \htext(9 5){$2$}
\htext(7 5){$2$} \htext(5 5){$2$}\htext(11 3){$3$} \htext(9
3){$3$} \htext(7 3){$3$} \htext(5 3){$3$} \htext(11 1.5){$4$}
\htext(9 1.5){$4$} \htext(7 1.5){$4$} \htext(5 1.5){$4$}\htext(3
1.5){$4$}\htext(11 0.5){$4$} \htext(9 0.5){$4$} \htext(7 0.5){$4$}
\htext(5 0.5){$4$}\htext(3 0.5){$4$}\htext(1 0.5){$4$}

\esegment

\move(47 0) \bsegment  \move(12 14)\dtri

\move(10 10)\dtri

\setgray 0.6 \move(-1 0) \lvec(12 0)\lvec(12 16) \move(12
16)\rlvec(-2 0)\rlvec(0 -16)\move(12 14)\rlvec(-2 0)\move(12
13)\rlvec(-2 0)\move(12 12)\rlvec(-2 0)\move(12 10)\rlvec(-4
0)\rlvec(0 -10)\move(12 8)\rlvec(-6 0)\rlvec(0 -8)\move(12
6)\rlvec(-6 0)\rlvec(0 -6)\move(12 4)\rlvec(-8 0)\rlvec(0 -4
)\move(12 2)\rlvec(-8 0)\move(12 1)\rlvec(-12 0)\rlvec(0
-1)\move(2 0)\rlvec(0 1)

\move(12 8)\rlvec(-2 -2) \move(10 8)\rlvec(-2 -2)\move(8
8)\rlvec(-2 -2)

\move(12 10)\linewd 0.15 \setgray 0 \rlvec(-2 0)\rlvec(0 -2)
\rlvec(-4 0) \rlvec(0 -7)\rlvec(-6 0)\rlvec(0 -1)\lvec(12
0)\lvec(12 10)\move(8 8)\rlvec(-2 -2)

\htext(11 15){$3$}\htext(11 13.5){$4$} \htext(11 12.5){$4$}
\htext(11 11){$3$} \htext(11 9){$2$}\htext(9 9){$2$}\htext(11.5
6.5){$1$} \htext(10.5 7.5){$0$} \htext(9.5 6.5){$0$} \htext(8.5
7.5){$1$} \htext(6.5 7.5){$0$} \htext(11 5){$2$} \htext(9 5){$2$}
\htext(7 5){$2$}
\htext(11 3){$3$} \htext(9 3){$3$} \htext(7 3){$3$} \htext(5
3){$3$}\htext(11 1.5){$4$} \htext(9 1.5){$4$} \htext(7 1.5){$4$}
\htext(5 1.5){$4$}\htext(11 0.5){$4$} \htext(9 0.5){$4$} \htext(7
0.5){$4$} \htext(5 0.5){$4$}\htext(3 0.5){$4$}\htext(1 0.5){$4$}

\esegment
\end{texdraw}
\end{center}

\noindent Here, the shaded parts represent $L_{\omega_3}^{\pm}$,
$L_{(\omega_3,\la_4)}^{\pm}$ and $L_{\omega_3}^{\pm}(3;2,1)$ in
the second, third and fourth Young walls, respectively.
\end{example}

%
\savebox{\tmpfiga}{\begin{texdraw} \fontsize{6}{6}\selectfont
\drawdim em \setunitscale 0.9

\setgray 0.6 \move(-2 0)\lvec(10 0)\lvec(10 12)\move(10 12)
\rlvec(-2 0)\rlvec(0 -12)\move(10 10)\rlvec(-2 0) \move(10
8)\rlvec(-4 0)\rlvec(0 -8)\move(10 6)\rlvec(-6 0)\rlvec(0
-6)\move(10 4)\rlvec(-8 0)\rlvec(0 -4)\move(10 2)\rlvec(-10
0)\rlvec(0 -2)

\move(10 8)\rlvec(-2 -2)\move(8 8)\rlvec(-2 -2)

\move(2 2)\rlvec(-2 -2)\move(4 2)\rlvec(-2 -2)\move(6 2)\rlvec(-2
-2)\move(8 2)\rlvec(-2 -2)\move(10 2)\rlvec(-2 -2)

\linewd 0.15 \setgray 0 \move(10 8)\rlvec(-2 -2)\rlvec(-2
0)\rlvec(0 -2) \rlvec(-2 0)\rlvec(0 -2)\rlvec(-2 0)\rlvec(-2
-2)\rlvec(10 0)\rlvec(0 8)

\linewd 0.1 \setgray 0.3 \move(-2 7)\rlvec(14 0)

\htext(9 11){$2$} \htext(9 9){$3$}

\htext(6.5 7.5){$5$}\htext(7.5 6.5){$4$}\htext(8.5
7.5){$4$}\htext(9.5 6.5){$5$}

\htext(5 5){$3$}\htext(7 5){$3$}\htext(9 5){$3$}\htext(3
3){$2$}\htext(5 3){$2$}\htext(7 3){$2$}\htext(9 3){$2$}

\htext(0.5 1.5){$1$} \htext(1.5 0.5){$0$} \htext(2.5 1.5){$0$}
\htext(3.5 0.5){$1$} \htext(4.5 1.5){$1$} \htext(5.5 0.5){$0$}
\htext(7.5 0.5){$1$} \htext(6.5 1.5){$0$} \htext(8.5 1.5){$1$}
\htext(9.5 0.5){$0$}
\end{texdraw}}
\savebox{\tmpfigb}{\begin{texdraw} \fontsize{6}{6}\selectfont
\drawdim em \setunitscale 0.9

\rlvec(0 2)\rlvec(4 0)\rlvec(0 4)\rlvec(-2 0)\rlvec(0 -6)\rlvec(2
2)\move(0 0)\rlvec(2 2)\move(2 4)\rlvec(2 0)\move(0 0)\rlvec(4
0)\rlvec(0 2)

\htext(0.5 1.5){$5$}\htext(2.5 1.5){$4$}\htext(3 3){$3$}\htext(3
5){$2$}
\end{texdraw}}
\savebox{\tmpfigc}{\begin{texdraw} \fontsize{6}{6}\selectfont
\drawdim em \setunitscale 0.9

\rlvec(0 2)\rlvec(4 0)\rlvec(0 4)\rlvec(4 0)\rlvec(0 2)\rlvec(-2
-2)\rlvec(0 -2)\rlvec(-4 0)\rlvec(0 -2)\rlvec(-2 -2)\move(0
0)\rlvec(2 0)\rlvec(0 2)\move(6 6)\rlvec(0 2)\rlvec(2 0)

\htext(0.5 1.5){$1$}\htext(3 3){$2$}\htext(5 5){$3$}\htext(7.5
6.5){$4$}
\end{texdraw}}
\savebox{\tmpfigd}{\begin{texdraw} \fontsize{6}{6}\selectfont
\drawdim em \setunitscale 0.9

\rlvec(0 8)\rlvec(2 0)\rlvec(-2 -2)\rlvec(2 0)\rlvec(0
-4)\rlvec(-2 -2)\move(0 2)\rlvec(2 0)\move(0 4)\rlvec(2 0)\move(2
6)\rlvec(0 2)\move(0 0)\rlvec(2 0)\rlvec(0 2)

\htext(0.5 1.5){$4$}\htext(1 3){$3$}\htext(1 5){$2$}\htext(0.5
7.5){$1$}
\end{texdraw}}

\vskip 3mm Finally, we focus on the case $\frak g = D_n$. If
${\overline Y_{\omega_{i_k}}}$ of $Y\in F(\la)$ contains a row
consisting of $n$-blocks and $(n-1)$-blocks, which will be called
the {\it $(n-1,n)$-row}, then we define the walls
$Y_{\omega_{i_k}}^{\pm}$ and $|Y_{\omega_{i_k}}^{-}|$ as in the
case of $\frak g = C_n$ or $B_n$.

\vskip 3mm
\begin{example} If $\frak g = D_5$, $\la = \omega_6$
and
$$Y= \raisebox{-0.4\height}{\usebox{\tmpfiga}}\in F(\la),$$
then we have
$$Y_{\omega_i}^+=\raisebox{-0.4\height}{\usebox{\tmpfigb}}\,, \quad
Y_{\omega_i}^-=\raisebox{-0.4\height}{\usebox{\tmpfigc}}\quad
\text{and}\quad
|Y_{\omega_i}^-|=\raisebox{-0.4\height}{\usebox{\tmpfigd}}\,.$$
\end{example}

Consider ${\overline L}_{\omega_{i_k} + \omega_{i_{k+1}}}$,
${\overline L}_{\omega_{i_t}+\la_n}$ or ${\overline
L}_{\omega_{i_t}+\la_{n-1}}$ of $Y\in F(\la)$.

\vskip 3mm
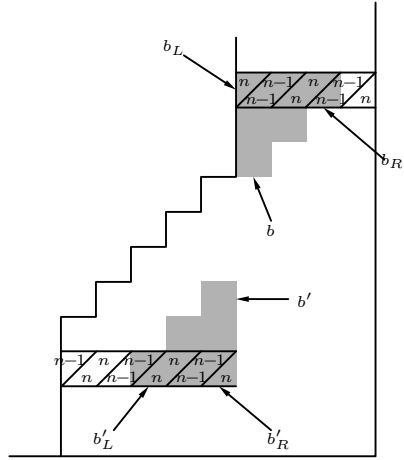
\begin{center}
\begin{texdraw}
\fontsize{6}{6}\selectfont \drawdim em \setunitscale 0.9
\nc{\dtri}{ \bsegment \lvec(-2 0) \lvec(-2 2)\lvec(0 2)\lvec(0
0)\ifill f:0.7 \esegment }

\nc{\dtrii}{ \bsegment \lvec(-2 0) \lvec(-2 1)\lvec(0 1)\lvec(0
0)\ifill f:0.7 \esegment }

\bsegment \move(16 26)\dtri \move(14 26)\dtri \move(12 26)\dtri
\move(12 24)\dtri \move(14 24)\dtri \move(12 22)\dtri \move(10
10)\dtri \move(10 12)\dtri \move(10 14)\dtri \move(8 10)\dtri
\move(8 12)\dtri \move(6 10)\dtri

\setgray 0 \move(-3 6) \lvec(18 6)\lvec(18 32) \move(10 30)
\rlvec(0 -8)\rlvec(-2 0)\rlvec(0 -2) \rlvec(-2 0)\rlvec(0
-2)\rlvec(-2 0)\rlvec(0 -2)\rlvec(-2 0)\rlvec(0 -2)\rlvec(-2
0)\rlvec(0 -8)

\move(18 28)\rlvec(-8 0)\move(18 26)\rlvec(-8 0) \move(10
12)\rlvec(-10 0)\move(10 10)\rlvec(-10 0)

\move(18 28)\rlvec(-2 -2)\move(16 28 )\rlvec(-2 -2)\move(14
28)\rlvec(-2 -2)\move(12 28)\rlvec(-2 -2)

\move(10 12)\rlvec(-2 -2)\move(8 12)\rlvec(-2 -2)\move(6
12)\rlvec(-2 -2)\move(4 12)\rlvec(-2 -2)\move(2 12)\rlvec(-2 -2)

\move(7 29)\avec(10 27)\move(18.5 23)\avec(15 26)

\move(3 7.5)\avec(5 10)\move(12 7.5)\avec(9 10)

\move(13 15)\avec(10 15)\move(12 19.5)\avec(11 22)

\htext(6.5 29.5){$b_L$}\htext(19 23){$b_R$}

\htext(2.5 7){$b'_L$}\htext(12.5 7){$b'_R$}

\htext(14 15){$b'$}\htext(12 19){$b$}

\htext(17.5 26.5){$n$} \htext(16.5 27.5){$n\!\!-\!\!1$}\htext(15.5
26.5){$n\!\!-\!\!1$} \htext(14.5 27.5){$n$}\htext(13.5 26.5){$n$}
\htext(12.5 27.5){$n\!\!-\!\!1$}\htext(11.5 26.5){$n\!\!-\!\!1$}
\htext(10.5 27.5){$n$}

\htext(9.5 10.5){$n$}\htext(8.5 11.5){$n\!\!-\!\!1$}\htext(7.5
10.5){$n\!\!-\!\!1$}\htext(6.5 11.5){$n$}\htext(5.5
10.5){$n$}\htext(4.5 11.5){$n\!\!-\!\!1$}\htext(3.5
10.5){$n\!\!-\!\!1$}\htext(2.5 11.5){$n$}\htext(1.5
10.5){$n$}\htext(0.5 11.5){$n\!\!-\!\!1$}

\esegment

\end{texdraw}
\end{center}

As we can see from above picture, ${\overline L}_{\omega_{i_k} +
\omega_{i_{k+1}}}$ contain two $(n-1,n)$-rows above ${\overline
H}_{\omega_{i_k}+\omega_{i_{k+1}}}$. We denote by $b_L$ the
left-most blocks in the upper $(n-1,n)$-row and $b_R$ the blocks
lying in the $(i-2)$-th column to the right of $b_L$. Then the
blocks $b_L$, $b_R$ and the block $b$ lying in the $(i-2)$-th row
below $b_L$ form a {\it right isosceles triangle}. We denote by
$L_{(\omega_{i_k},\omega_{i_{k+1}})}^-$ the part of ${\overline
L}_{\omega_{i_k} + \omega_{i_{k+1}}}$ consisting of this right
isosceles triangle. Note that the size of
$L_{(\omega_{i_k},\omega_{i_{k+1}})}^-$ in the case of $D_n$ is
smaller than that of $L_{(\omega_{i_k},\omega_{i_{k+1}})}^-$ in
the case of $C_n$ or $B_n$.

Similarly, let $b'_R$ the right-most blocks in the lower
$(n-1,n)$-row outside the highest weight vector ${\overline
H}_{\omega_{i_k} + \omega_{i_{k+1}}}$ and let $b'_L$ the blocks
lying in the $(i-2)$-th column to the left of $b'_R$. Then the
blocks $b'_R$, $b'_L$ and the block $b'$ lying in the $(i-2)$-th
row above $b'_R$ form another {\it right isosceles triangle}. We
denote by $L_{(\omega_{i_k},\omega_{i_{k+1}})}^+$ the part of
${\overline L}_{\omega_{i_k} + \omega_{i_{k+1}}}$ consisting of
this right isosceles triangle. Note that
$L_{(\omega_{i_k},\omega_{i_{k+1}})}^-$ and
$L_{(\omega_{i_k},\omega_{i_{k+1}})}^+$ are of the same size with
each base of length $i-1$. Now, we can also define
$L_{(\omega_{i_t},\la_{n-1})}^{\pm}$ and
$L_{(\omega_{i_t},\la_n)}^{\pm}$ in a similar way, and set
\begin{equation}
\begin{aligned}
& Y_{(\omega_{i_k},\omega_{i_{k+1}})}^{-} = Y \cap
L_{(\omega_{i_k},\omega_{i_{k+1}})}^{-},
\quad Y_{(\omega_{i_k},\omega_{i_{k+1}})}^{+} =
Y \cap L_{(\omega_{i_k},\omega_{i_{k+1}})}^{+}, \\
& Y_{(\omega_{i_t},\la_{n-1})}^{\pm} = Y \cap
L_{(\omega_{i_t},\la_{n-1})}^{\pm}, \quad
Y_{(\omega_{i_t},\la_{n})}^{\pm} = Y \cap
L_{(\omega_{i_t},\la_{n})}^{\pm}.
\end{aligned}
\end{equation}
As usual, let $|Y_{(\omega_{i_k},\omega_{i_{k+1}})}^{-}|$ (resp.
$|Y_{(\omega_{i_t},\la_{n-1})}^{-}|$ and
$|Y_{(\omega_{i_t},\la_n)}^{-}|$) be the wall obtained by
reflecting $Y_{(\omega_{i_k},\omega_{i_{k+1}})}^{-}$ (resp.
$Y_{(\omega_{i_t},\la_{n-1})}^{-}$ and
$Y_{(\omega_{i_t},\la_n)}^{-}$) with respect to the $(n-1,n)$-row
and shifting the blocks to the right as much as possible.

\savebox{\tmpfiga}{\begin{texdraw} \fontsize{6}{6}\selectfont
\drawdim em \setunitscale 0.9 \nc{\dtri}{ \bsegment \lvec(-2 0)
\lvec(-2 2)\lvec(0 2)\lvec(0 0)\ifill f:0.7 \esegment }

\nc{\dtrii}{ \bsegment \lvec(-2 0) \lvec(-2 -2)\lvec(0 0)\ifill
f:0.7 \esegment }

\nc{\dtriii}{ \bsegment \lvec(-2 0)\lvec(0 2)\lvec(0 0)\ifill
f:0.7 \esegment }

\bsegment \move(16 16)\dtriii \move(14 16)\dtriii \move(18
16)\dtriii

\move(16 14)\dtri \move(14 14)\dtri \move(14 12)\dtri

 \move(12 2)\dtri \dtrii \move(10 2)\dtri
\dtrii \move(8 2)\dtrii  \move(12 4)\dtri

\setgray 0.6 \move(-1 0) \lvec(20 0)\lvec(20 24)\move(20
24)\rlvec(-2 0)\rlvec(0 -24)\move(20 22)\rlvec(-4 0)\rlvec(0
-22)\move(20 20)\rlvec(-4 0)\rlvec(0 -20) \move(20 18) \rlvec(-6
0)\rlvec(0 -18) \move(20 16)\rlvec(-6 0)\move(20 14) \rlvec(-8
0)\rlvec(0 -14) \move(20 12) \rlvec(-8 0)\move(20 10)\rlvec(-10
0)\rlvec(0 -10) \move(20 8)\rlvec(-12 0)\rlvec(0 -8) \move(20
6)\rlvec(-12 0)\move(20 4)\rlvec(-14 0)\rlvec(0 -4) \move(20
2)\rlvec(-16 0)\rlvec(0 -2)

\move(20 18)\rlvec(-2 -2)\move(18 18)\rlvec(-2 -2)\move(16
18)\rlvec(-2 -2) \move(20 10)\rlvec(-2 -2) \move(18 10)\rlvec(-2
-2) \move(16 10)\rlvec(-2 -2)\move(14 10)\rlvec(-2 -2)\move(12
10)\rlvec(-2 -2) \move(20 2)\rlvec(-2 -2) \move(18 2)\rlvec(-2 -2)
\move(16 2)\rlvec(-2 -2)\move(14 2)\rlvec(-2 -2) \move(12
2)\rlvec(-2 -2)\move(10 2)\rlvec(-2 -2)\move(8 2)\rlvec(-2
-2)\move(6 2)\rlvec(-2 -2)\move(4 2)\rlvec(-2 -2)

\move(20 14)\linewd 0.15 \setgray 0 \rlvec(-2 0)\rlvec(0
-2)\rlvec(-2 0)\rlvec(0 -2)\rlvec(-4 0)\rlvec(0 -8)\rlvec(-2
-2)\rlvec(0 2)\rlvec(-2 -2)\rlvec(0 2)\rlvec(-2 -2)\rlvec(0
2)\rlvec(-2 -2)\rlvec(0 2)\rlvec(-2 -2)\lvec(20 0)\lvec(20
14)\move(14 10)\rlvec(-2 -2)

\htext(19 23){$2$} \htext(19 21){$3$} \htext(17 21){$3$} \htext(19
19){$4$} \htext(17 19){$4$} \htext(19.5 16.5){$6$} \htext(18.5
17.5){$5$}\htext(17.5 16.5){$5$} \htext(16.5 17.5){$6$}\htext(15.5
16.5){$6$} \htext(14.5 17.5){$5$} \htext(19 15){$4$} \htext(17
15){$4$} \htext(15 15){$4$} \htext(15 13){$3$}\htext(17 13){$3$}
\htext(19 13){$3$} \htext(13 13){$3$} \htext(15 11){$2$}\htext(17
11){$2$} \htext(19 11){$2$} \htext(13 11){$2$}

\htext(19.5 8.5){$0$}\htext(18.5 9.5){$1$} \htext(17.5
8.5){$1$}\htext(16.5 9.5){$0$} \htext(15.5 8.5){$0$}\htext(14.5
9.5){$1$} \htext(13.5 8.5){$1$}\htext(12.5 9.5){$0$} \htext(10.5
9.5){$1$}

\htext(9 7){$2$}\htext(19 7){$2$} \htext(17 7){$2$} \htext(15
7){$2$}\htext(13 7){$2$} \htext(11 7){$2$} \htext(9
5){$5$}\htext(19 5){$3$} \htext(17 5){$3$} \htext(15
5){$3$}\htext(13 5){$3$} \htext(11 5){$3$} \htext(7
3){$4$}\htext(19 3){$4$} \htext(17 3){$4$}\htext(15 3){$4$}
\htext(13 3){$4$} \htext(11 3){$4$}\htext(9 3){$4$}

\htext(19.5 0.5){$6$}\htext(18.5 1.5){$5$}\htext(17.5
0.5){$5$}\htext(16.5 1.5){$6$}\htext(15.5 0.5){$6$}\htext(14.5
1.5){$5$}\htext(13.5 0.5){$5$}\htext(12.5 1.5){$6$}\htext(11.5
0.5){$6$}\htext(10.5 1.5){$5$}\htext(9.5 0.5){$5$}\htext(8.5
1.5){$6$}\htext(7.5 0.5){$6$}\htext(6.5 1.5){$5$}\htext(5.5
0.5){$5$}\htext(4.5 1.5){$6$}\htext(3.5 0.5){$6$} \esegment
\end{texdraw}}
\savebox{\tmpfigb}{\begin{texdraw} \fontsize{6}{6}\selectfont
\drawdim em \setunitscale 0.9 \nc{\dtri}{ \bsegment \lvec(-2 0)
\lvec(-2 2)\lvec(0 2)\lvec(0 0)\ifill f:0.7 \esegment }

\nc{\dtrii}{ \bsegment \lvec(-2 0) \lvec(-2 -2)\lvec(0 0)\ifill
f:0.7 \esegment }

\nc{\dtriii}{ \bsegment \lvec(-2 0)\lvec(0 2)\lvec(0 0)\ifill
f:0.7 \esegment }

\rlvec(0 2)\rlvec(4 0)\rlvec(0 4)\rlvec(-2 -2)\rlvec(0
-4)\rlvec(-2 0)\move(2 4)\rlvec(4 0)\rlvec(0 2) \rlvec(-2 -2)

\htext(1 1){$3$}\htext(3 3){$4$} \htext(3.5 4.5){$6$} \htext(5.5
4.5){$5$}

\end{texdraw}}
\savebox{\tmpfigc}{\begin{texdraw} \fontsize{6}{6}\selectfont
\drawdim em \setunitscale 0.9 \nc{\dtri}{ \bsegment \lvec(-2 0)
\lvec(-2 2)\lvec(0 2)\lvec(0 0)\ifill f:0.7 \esegment }

\nc{\dtrii}{ \bsegment \lvec(-2 0) \lvec(-2 -2)\lvec(0 0)\ifill
f:0.7 \esegment }

\nc{\dtriii}{ \bsegment \lvec(-2 0)\lvec(0 2)\lvec(0 0)\ifill
f:0.7 \esegment }

\rlvec(6 0)\rlvec(0 6)\rlvec(-2 0)\rlvec(0 -6)\move(0 0)\rlvec(0
2) \rlvec(6 0)\move(2 0)\rlvec(0 4) \rlvec(4 0)\move(6 2)\rlvec(-2
-2) \move(4 2)\rlvec(-2 -2) \move(2 2)\rlvec(-2 -2)

\htext(0.5 1.5){$5$} \htext(2.5 1.5){$6$} \htext(3 3){$4$}
\htext(4.5 1.5){$5$} \htext(5 3){$4$} \htext(5 5){$3$}
\end{texdraw}}
\savebox{\tmpfigd}{\begin{texdraw} \fontsize{6}{6}\selectfont
\drawdim em \setunitscale 0.9 \nc{\dtri}{ \bsegment \lvec(-2 0)
\lvec(-2 2)\lvec(0 2)\lvec(0 0)\ifill f:0.7 \esegment }

\nc{\dtrii}{ \bsegment \lvec(-2 0) \lvec(-2 -2)\lvec(0 0)\ifill
f:0.7 \esegment }

\nc{\dtriii}{ \bsegment \lvec(-2 0)\lvec(0 2)\lvec(0 0)\ifill
f:0.7 \esegment }

\rlvec(4 0)\rlvec(0 6)\rlvec(-2 0)\rlvec(0 -6)\move(0 0)\rlvec(0
2) \rlvec(4 0)\move(2 4)\rlvec(2 0) \move(4 2)\rlvec(-2 -2)
\move(2 2)\rlvec(-2 -2)

\htext(0.5 1.5){$6$} \htext(2.5 1.5){$5$} \htext(3 3){$4$}\htext(3
5){$3$}

\end{texdraw}}

\begin{example}

If $\frak g = D_7$, $\la =\omega_5+\la_6$, and
$$Y= \raisebox{-0.4\height}{\usebox{\tmpfiga}}\in F(\la),$$
then we have
$$Y_{(\omega_5,\la_6)}^-=\raisebox{-0.4\height}{\usebox{\tmpfigb}}\,, \quad
Y_{(\omega_5,\la_6)}^+=\raisebox{-0.4\height}{\usebox{\tmpfigc}}
\quad \text{and} \quad
|Y_{(\omega_5,\la_6)}^-|=\raisebox{-0.4\height}{\usebox{\tmpfigd}}\,.$$
Here, the shaded parts represent $L_{(\omega_5,\la_6)}^-$ and
$L_{(\omega_5,\la_6)}^+$, respectively.

\end{example}

\vskip 5mm Assume that ${\overline Y}_{\omega_{i_k} +
\omega_{i_{k+1}}}$, ${\overline Y}_{\omega_{i_t}+\la_n}$ or
${\overline Y}_{\omega_{i_t}+\la_{n-1}}$ of $Y\in F(\la)$
satisfies {\bf (C1)}. Then we can define
$L_{\omega_{i_k}}^{\pm}(a;p,q)$,
$L_{\omega_{i_{k+1}}}^{\pm}(a;p,q)$,
$Y_{\omega_{i_k}}^{\pm}(a;p,q)$,
$Y_{\omega_{i_{k+1}}}^{\pm}(a;p,q)$ and
$|Y_{\omega_{i_{k+1}}}^{\pm}(a;p,q)|$ as in the case of $\frak
g=C_n$ or $B_n$.

\vskip 2mm Now, suppose that ${\overline Y}_{\omega_{i_k} +
\omega_{i_{k+1}}}$, ${\overline Y}_{\omega_{i_t}+\la_n}$ or
${\overline Y}_{\omega_{i_t}+\la_{n-1}}$ of $Y\in F(\la)$ having
the following configuration:

\vskip 3mm

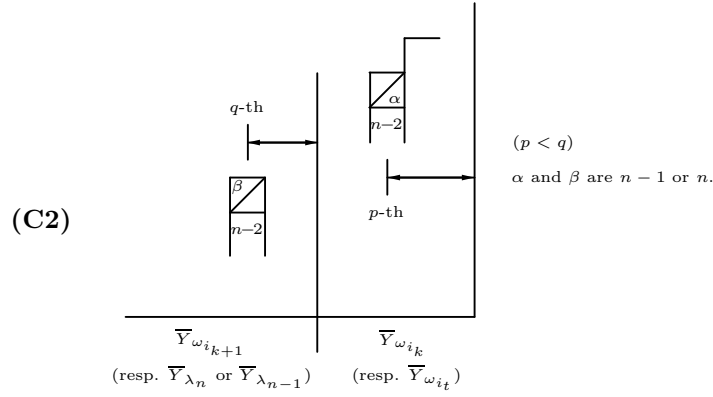
\savebox{\tmpfiga}{\begin{texdraw} \fontsize{6}{6}\selectfont
\drawdim em \setunitscale 0.9 \nc{\dtri}{ \bsegment \lvec(-2 0)
\lvec(-2 2)\lvec(0 2)\lvec(0 0)\ifill f:0.7 \esegment }

\bsegment \setgray 0 \move(0 0) \lvec(0 18) \move(-2 16 )\rlvec(-2
0)\move(-4 14)\rlvec(-2 0)\move(-4 12)\rlvec(-2 0) \move(-12
8)\rlvec(-2 0) \move(-12 6)\rlvec(-2 0) \move(-4 10)\rlvec(0
6)\move(-6 10)\rlvec(0 4) \move(-12 3.5)\rlvec(0 4.5)\move(-14
3.5)\rlvec(0 4.5) \move(0 0)\rlvec(-20 0) \move(-9 -2)\rlvec(0
16)\move(-4 14)\rlvec(-2 -2)\move(-12 8)\rlvec(-2 -2)

\move(-5 8)\ravec(5 0) \move(0 8)\ravec(-5 0) \move(-5 9)\rlvec(0
-2)\htext(-5 6){$p$-th}

\move(-13 10)\ravec(4 0) \move(-9 10)\ravec(-4 0) \move(-13
9)\rlvec(0 2)\htext(-13 12){$q$-th}

\htext(-4.5 12.5){$\alpha$}\htext(-5 11){$n\!\!-\!\!2$}
\htext(-13.5 7.5){$\beta$} \htext(-13 5){$n\!\!-\!\!2$}

\htext(-15 -1.5){$\overline{Y}_{\omega_{i_{k+1}}}$} \htext(-15
-3.5){(resp. $\overline{Y}_{\la_n}$ or
$\overline{Y}_{\la_{n-1}}$)}

\htext(-4 -1.5){$\overline{Y}_{\omega_{i_k}}$}\htext(-4 -3.5){
(resp. $\overline{Y}_{\omega_{i_t}}$)} \htext(4
10){($p<q$)}\htext(8 8){$\alpha$ and $\beta$ are $n-1$ or $n$.}

\esegment
\end{texdraw}}

\begin{center}
{\bf (C2)}\qquad $\raisebox{-0.4\height}{\usebox{\tmpfiga}}$
\end{center}

\savebox{\tmpfiga}{\begin{texdraw} \fontsize{6}{6}\selectfont
\drawdim em \setunitscale 0.9

\rlvec(0 3.5)\rlvec(2 0)\rlvec(0 -3.5)\move(0 1.5)\rlvec(2
0)\move(2 3.5)\rlvec(-2 -2)

\htext(0.5 2.75){$n$}\htext(1 0.75){$n\!\!-\!\!2$}
\end{texdraw}}
\savebox{\tmpfigb}{\begin{texdraw} \fontsize{6}{6}\selectfont
\drawdim em \setunitscale 0.9

\rlvec(0 3.5)\rlvec(2 0)\rlvec(0 -3.5)\move(0 1.5)\rlvec(2
0)\move(0 1.5)\rlvec(1.3 1.3)

\htext(1 3){$n\!\!-\!\!1$}\htext(1 0.75){$n\!\!-\!\!2$}
\end{texdraw}}
\savebox{\tmpfigc}{\begin{texdraw} \fontsize{6}{6}\selectfont
\drawdim em \setunitscale 0.9

\rlvec(0 3.5)\rlvec(2 0)\rlvec(0 -3.5)\move(0 1.5)\rlvec(2
0)\move(2 3.5)\rlvec(-2 -2)

\htext(1.5 2){$n$}\htext(1 0.75){$n\!\!-\!\!2$}
\end{texdraw}}
\savebox{\tmpfigd}{\begin{texdraw} \fontsize{6}{6}\selectfont
\drawdim em \setunitscale 0.9

\rlvec(0 3.5)\rlvec(2 0)\rlvec(0 -3.5)\move(0 1.5)\rlvec(2
0)\move(2 3.5)\rlvec(-1.3 -1.3)

\htext(1 2){$n\!\!-\!\!1$}\htext(1 0.75){$n\!\!-\!\!2$}
\end{texdraw}}

\vskip 3mm \noindent That is, the top of the $p$-th column of
$\overline{Y}_{\omega_{i_k}}$ from the right is
$\raisebox{-0.4\height}{\usebox{\tmpfiga}}$\,,
$\raisebox{-0.4\height}{\usebox{\tmpfigc}}$\,,
$\raisebox{-0.4\height}{\usebox{\tmpfigb}}$ or
$\raisebox{-0.4\height}{\usebox{\tmpfigd}}$\,. In the case of
$\raisebox{-0.4\height}{\usebox{\tmpfiga}}$ or
$\raisebox{-0.4\height}{\usebox{\tmpfigc}}$ (resp.
$\raisebox{-0.4\height}{\usebox{\tmpfigb}}$ or
$\raisebox{-0.4\height}{\usebox{\tmpfigd}}$\,), the top of the
$q$-th column of $\overline{Y}_{\omega_{i_{k+1}}}$ from the right
is $\raisebox{-0.4\height}{\usebox{\tmpfiga}}$ or
$\raisebox{-0.4\height}{\usebox{\tmpfigc}}$ (resp.
$\raisebox{-0.4\height}{\usebox{\tmpfigb}}$ or
$\raisebox{-0.4\height}{\usebox{\tmpfigd}}$\,) when $q-p$ is odd
and the top of the $q$-th column of
$\overline{Y}_{\omega_{i_{k+1}}}$ from the right is
$\raisebox{-0.4\height}{\usebox{\tmpfigb}}$ or
$\raisebox{-0.4\height}{\usebox{\tmpfigd}}$ (resp.
$\raisebox{-0.4\height}{\usebox{\tmpfiga}}$ or
$\raisebox{-0.4\height}{\usebox{\tmpfigc}}$\,) when $q-p$ is even.

\savebox{\tmpfiga}{\begin{texdraw} \fontsize{6}{6}\selectfont
\drawdim em \setunitscale 0.9

\rlvec(0 2)\rlvec(2 0)\rlvec(0 -2)\rlvec(-2 0)

\htext(1 1){$n\!\!-\!\!q$}
\end{texdraw}}
\savebox{\tmpfigb}{\begin{texdraw} \fontsize{6}{6}\selectfont
\drawdim em \setunitscale 0.9

\rlvec(0 2)\rlvec(2 0)\rlvec(0 -2)\rlvec(-2 0)

\htext(1 1.5){$n\!\!-\!\!q$}\htext(1 0.5){$+p\!\!-\!\!1$}
\end{texdraw}}
\savebox{\tmpfigc}{\begin{texdraw}\fontsize{6}{6}\selectfont
\drawdim em \setunitscale 0.9

\rlvec(0 2)\rlvec(2 0)\rlvec(0 -2)\rlvec(-2 0)

\htext(1 1){$n\!\!-\!\!i$}
\end{texdraw}}
\savebox{\tmpfigd}{\begin{texdraw} \fontsize{6}{6}\selectfont
\drawdim em \setunitscale 0.9

\rlvec(0 2)\rlvec(2 0)\rlvec(0 -2)\rlvec(-2 0)

\htext(1 1.5){$n\!\!-\!\!i$}\htext(1 0.5){$+p\!\!-\!\!1$}
\end{texdraw}}

\vskip 3mm We define $L_{\omega_{i_k}}(n-1,n;p,q)$ to be the
parallelogram formed by the $(n-q)$-block and $(n-q+p-1)$-block in
the $q$-th column, and $(n-i)$-block and $(n-i+p-1)$-block in the
$i$-th column lying below the ($n-1,n$)-row. We will denote by
$L_{\omega_{i_{k+1}}}(n-1,n;p,q)$ the parallelogram formed by the
$(n-q)$-block and $(n-i)$-block in the first column, and
$(n-q+p-1)$-block and $(n-i+p-1)$-block in the $p$-th column lying
above the ($n-1,n$)-row. Similarly, we can define the
parallelograms $L_{\omega_{i_t}}(n-1,n;p,q)$,
$L_{\la_n}(n-1,n;p,q)$ and $L_{\la_{n-1}}(n-1,n;p,q)$. The shaded
parts in the following picture represent
$L_{\omega_{i_k}}(n-1,n;p,q)$ and
$L_{\omega_{i_{k+1}}}(n-1,n;p,q)$.

\vskip 5mm
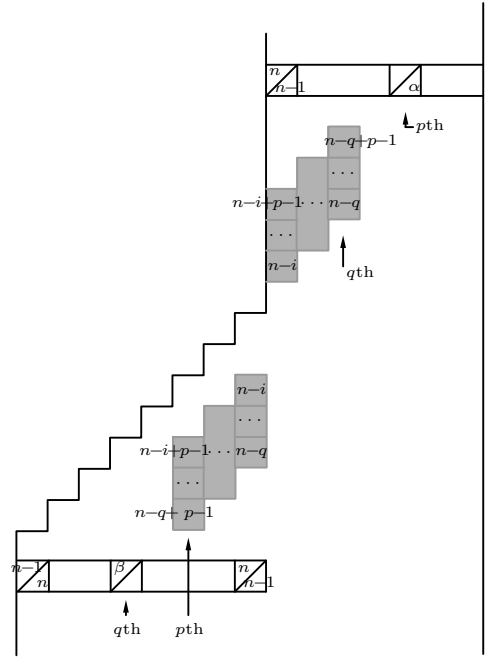
\begin{center}
\begin{texdraw}
\fontsize{6}{6}\selectfont \drawdim em \setunitscale 0.8
\nc{\dtri}{ \bsegment \lvec(-2 0) \lvec(-2 2)\lvec(0 2)\lvec(0
0)\ifill f:0.7 \esegment }

\nc{\dtrii}{ \bsegment \lvec(-2 0) \lvec(-2 1)\lvec(0 1)\lvec(0
0)\ifill f:0.7 \esegment }

\bsegment  \move(22 32)\dtri \move(22 34)\dtri \move(22 36)\dtri
\move(20 30)\dtri \move(20 32)\dtri \move(20 34)\dtri \move(18
28)\dtri \move(18 30)\dtri \move(18 32)\dtri

\move(16 20)\dtri \move(16 18)\dtri \move(16 16)\dtri \move(14
18)\dtri \move(14 16)\dtri \move(14 14)\dtri \move(12 16)\dtri
\move(12 14)\dtri \move(12 12)\dtri

\setgray 0 \move(30 4)\rlvec(0 42)\move(16 44)\rlvec(0 -18)

\move(16 26)\rlvec(-2 0)\rlvec(0 -2)\rlvec(-2 0)\rlvec(0
-2)\rlvec(-2 0)\rlvec(0 -2)\rlvec(-2 0)\rlvec(0 -2)\rlvec(-2
0)\rlvec(0 -2)\rlvec(-2 0)\rlvec(0 -2)\rlvec(-2 0)\rlvec(0
-2)\rlvec(-2 0)\rlvec(0 -8)

\move(16 40)\rlvec(14 0)\move(16 42)\rlvec(14 0)\move(16
40)\rlvec(2 2)\rlvec(0 -2)\move(24 40)\rlvec(2 2)\rlvec(0
-2)\move(24 40)\rlvec(0 2)

\move(0 8)\rlvec(16 0)\move(0 10)\rlvec(16 0)\move(0 8)\rlvec(2
2)\rlvec(0 -2)\move(6 8)\rlvec(2 2)\rlvec(0 -2)\move(6 8)\rlvec(0
2)\move(14 8)\rlvec(2 2)\rlvec(0 -2)\move(14 8)\rlvec(0 2)

\move(25 38)\ravec(0 1)\move(25 38)\rlvec(0.5 0)\htext(26.5
38){$p$th} \move(21 29)\ravec(0 2)\htext(22 28.5){$q$th}

\move(11 6.5)\ravec(0 5)\htext(11 5.5){$p$th} \move(7 6.5)\ravec(0
1)\htext(7 5.5){$q$th}

\htext(25.5 40.5){$\alpha$}\htext(22
37){$n\!\!-\!\!q\!\!+\!\!p\!\!-\!\!1$} \htext(21
35){$\cdots$}\htext(21 33){$n\!\!-\!\!q$}\htext(16
33){$n\!\!-\!\!i\!\!+\!\!p\!\!-\!\!1$}\htext(17
31){$\cdots$}\htext(17 29){$n\!\!-\!\!i$}\htext(16.5 41.5){$n$}
\htext(17.5 40.5){$n\!\!-\!\!1$}\htext(19 33){$\cdots$}

\htext(6.5 9.5){$\beta$}\htext(15 21){$n\!\!-\!\!i$} \htext(15
19){$\cdots$}\htext(15 17){$n\!\!-\!\!q$}\htext(10
17){$n\!\!-\!\!i\!\!+\!\!p\!\!-\!\!1$}\htext(11
15){$\cdots$}\htext(10 13){$n\!\!-\!\!q\!\!+p\!\!-\!\!1$}
\htext(0.5 9.5){$n\!\!-\!\!1$} \htext(1.5 8.5){$n$}\htext(14.5
9.5){$n$} \htext(15.5 8.5){$n\!\!-\!\!1$}\htext(13 17){$\cdots$}

\setgray 0.6 \move(16 34)\rlvec(2 0)\rlvec(0 2)\rlvec(2 0)\rlvec(0
2)\rlvec(2 0)\rlvec(0 -6)\rlvec(-2 0)\rlvec(0 -2)\rlvec(-2
0)\rlvec(0 -2)\rlvec(-2 0)\move(16 32)\rlvec(2 0)\move(16
30)\rlvec(2 0) \move(20 34)\rlvec(2 0)\move(18 34)\rlvec(0
-4)\move(20 36)\rlvec(0 -4)\move(20 36)\rlvec(2 0)

\move(10 18)\rlvec(2 0)\rlvec(0 2)\rlvec(2 0)\rlvec(0 2)\rlvec(2
0)\rlvec(0 -6)\rlvec(-2 0)\rlvec(0 -2)\rlvec(-2 0)\rlvec(0
-2)\rlvec(-2 0)\rlvec(0 6)\move(10 16)\rlvec(2 0) \move(14
20)\rlvec(2 0)\move(14 18)\rlvec(2 0) \move(12 18)\rlvec(0
-4)\move(14 20)\rlvec(0 -4)\move(10 14)\rlvec(2 0)

\esegment
\end{texdraw}
\end{center}
\noindent Here, $\alpha$ and $\beta$ are $n$ or $n-1$ in the
hypothesis.

\vskip 3mm  We define
$$\aligned
Y_{\omega_{i_k}}(n-1,n;p,q)&=L_{\omega_{i_k}}(n-1,n;p,q)\cap Y, \\
Y_{\omega_{i_{k+1}}}(n-1,n;p,q)&=L_{\omega_{i_{k+1}}}(n-1,n;p,q)\cap
Y,\\
Y_{\omega_{i_t}}(n-1,n;p,q)&=L_{\omega_{i_t}}(n-1,n;p,q)\cap Y, \\
Y_{\la_{n-1}}(n-1,n;p,q)&=L_{\la_{n-1}}(n-1,n;p,q)\cap Y, \\
Y_{\la_{n}}(n-1,n;p,q)&=L_{\la_{n}}(n-1,n;p,q)\cap Y,
\endaligned$$
and let $|Y_{\omega_{i_k}}(n-1,n;p,q)|$ be the wall obtained by
reflecting $Y_{\omega_{i_k}}(n-1,n;p,q)$ with respect to the
$(n-1,n)$-arrow and shifting the blocks to the right as much as
possible and let $Y_{\omega_{i_{k+1}}}^t(n-1,n;p,q)$ be the wall
obtained by shifting the blocks of
$Y_{\omega_{i_{k+1}}}(n-1,n;p,q)$ to the right as much as
possible.

\savebox{\tmpfiga} {\begin{texdraw} \fontsize{6}{6}\selectfont
\drawdim em \setunitscale 0.9 \nc{\dtri}{ \bsegment \lvec(-2 0)
\lvec(-2 2)\lvec(0 2)\lvec(0 0)\ifill f:0.85 \esegment }

\nc{\dtrii}{ \bsegment \lvec(-2 0) \lvec(-2 -2)\lvec(0 0)\ifill
f:0.7 \esegment }

\nc{\dtriii}{ \bsegment \lvec(-2 0) \lvec(0 2)\lvec(0 0)\ifill
f:0.6 \esegment }

\bsegment \move(18 36)\dtriii \move(6 12)\dtriii  \move(16
34)\dtri   \move(16 32)\dtri \move(14 32)\dtri \move(16 30)\dtri
\move(14 30)\dtri \move(14 28)\dtri

\move(12 20)\dtri \move(12 18)\dtri   \move(10 18)\dtri  \move(10
16)\dtri  \move(8 16)\dtri \move(8 14)\dtri

\setgray 0.6 \move(-3 0) \lvec(22 0)\rlvec(0 40)\move(20
24)\rlvec(0 16) \move(18 24)\rlvec(0 14)\move(16 24)\rlvec(0 12)
\move(14 0)\rlvec(0 34) \move(12 0)\rlvec(0 32) \move(10
0)\rlvec(0 22)\move(8 0)\rlvec(0 20)\move(6 0)\rlvec(0 16)\move(4
0)\rlvec(0 12)\move(2 0)\rlvec(0 10)\move(0 0)\rlvec(0 8)\move(20
0)\rlvec(0 2)\move(18 0)\rlvec(0 2)\move(16 0)\rlvec(0 2)

\move(22 40)\rlvec(-2 0)\move(22 38)\rlvec(-4 0)\move(22
36)\rlvec(-6 0)\move(22 34)\rlvec(-8 0)\move(22 32)\rlvec(-10
0)\move(22 30)\rlvec(-10 0)\move(22 28)\rlvec(-10 0)\move(22
26)\rlvec(-10 0)\move(22 24)\rlvec(-10 0) \move(14 22)\rlvec(-4 0)
\move(14 20)\rlvec(-6 0)\move(14 18)\rlvec(-6 0)\move(14
16)\rlvec(-8 0)\move(14 14)\rlvec(-8 0)\move(14 12)\rlvec(-10
0)\move(14 10)\rlvec(-12 0)\move(14 8)\rlvec(-14 0)\move(14
6)\rlvec(-14 0)\move(14 4)\rlvec(-14 0)\move(22 2)\rlvec(-22 0)

\move(22 38)\rlvec(-2 -2)\move(20 38)\rlvec(-2 -2)\move(18
38)\rlvec(-2 -2)\move(22 26)\rlvec(-2 -2)\move(20 26)\rlvec(-2
-2)\move(18 26)\rlvec(-2 -2)\move(16 26)\rlvec(-2 -2)\move(14
26)\rlvec(-2 -2) \move(14 14)\rlvec(-2 -2) \move(12 14)\rlvec(-2
-2)\move(10 14)\rlvec(-2 -2)\move(8 14)\rlvec(-2 -2)\move(6
14)\rlvec(-2 -2) \move(22 2)\rlvec(-2 -2)\move(20 2)\rlvec(-2
-2)\move(18 2)\rlvec(-2 -2)\move(16 2)\rlvec(-2 -2)\move(14
2)\rlvec(-2 -2) \move(12 2)\rlvec(-2 -2)\move(10 2)\rlvec(-2
-2)\move(8 2)\rlvec(-2 -2)\move(6 2)\rlvec(-2 -2)\move(4
2)\rlvec(-2 -2) \move(2 2)\rlvec(-2 -2)

\move(22 22)\rlvec(-8 0)\move(22 20)\rlvec(-8 0)\move(22
18)\rlvec(-8 0)\move(22 16)\rlvec(-8 0) \move(22 14)\rlvec(-8
0)\move(22 12)\rlvec(-8 0)\move(22 10)\rlvec(-8 0)\move(22
8)\rlvec(-8 0)\move(22 6)\rlvec(-8 0)\move(22 4)\rlvec(-8 0)
\move(20 2)\rlvec(0 22)\move(18 2)\rlvec(0 22)\move(16 2)\rlvec(0
22) \move(22 14)\rlvec(-2 -2)\move(20 14)\rlvec(-2 -2)\move(18
14)\rlvec(-2 -2)\move(16 14)\rlvec(-2 -2)

\move(22 32)\linewd 0.15 \setgray 0 \rlvec(-2 0)\rlvec(0
-2)\rlvec(-2 0)\rlvec(0 -2)\rlvec(-2 0)\rlvec(0 -2)\rlvec(-2
0)\rlvec(-2 -2)\rlvec(0 -14)\rlvec(-2 0)\rlvec(0 -2)\rlvec(-2
0)\rlvec(0 -2)\rlvec(-2 0)\rlvec(0 -2)\rlvec(-2 0)\rlvec(0
-2)\rlvec(-2 0)\rlvec(-2 -2)\rlvec(22 0)\rlvec(0 32)

\htext(21 39){$6$}

\htext(21.5 36.5){$8$}\htext(20.5 37.5){$7$} \htext(19.5
36.5){$7$}\htext(18.5 37.5){$8$}\htext(17.5 36.5){$8$}

\htext(21 35){$6$}\htext(19 35){$6$}\htext(17 35){$6$} \htext(21
33){$5$} \htext(19 33){$5$}\htext(17 33){$5$}\htext(15 33){$5$}
\htext(21 31){$4$} \htext(19 31){$4$}\htext(17 31){$4$}\htext(15
31){$4$}\htext(13 31){$4$} \htext(21 29){$3$}\htext(19
29){$3$}\htext(17 29){$3$}\htext(15 29){$3$}\htext(13 29){$3$}
\htext(21 27){$2$} \htext(19 27){$2$}\htext(17 27){$2$}\htext(15
27){$2$}\htext(13 27){$2$}

\htext(21.5 24.5){$0$}\htext(20.5 25.5){$1$} \htext(19.5
24.5){$1$}\htext(18.5 25.5){$0$}\htext(17.5 24.5){$0$}\htext(16.5
25.5){$1$}\htext(15.5 24.5){$1$}\htext(14.5 25.5){$0$}\htext(13.5
24.5){$0$} \htext(12.5 25.5){$1$}

\htext(13 23){$2$}\htext(15 23){$2$}\htext(17 23){$2$}\htext(19
23){$2$}\htext(21 23){$2$} \htext(13 21){$3$} \htext(11
21){$3$}\htext(15 21){$3$} \htext(17 21){$3$}\htext(19 21){$3$}
\htext(21 21){$3$} \htext(13 19){$4$}\htext(11 19){$4$} \htext(9
19){$4$} \htext(15 19){$4$}\htext(17 19){$4$} \htext(19
19){$4$}\htext(21 19){$4$}\htext(13 17){$5$}\htext(11 17){$5$}
\htext(9 17){$5$}\htext(15 17){$5$}\htext(17 17){$5$} \htext(19
17){$5$}\htext(21 17){$5$} \htext(13 15){$6$}\htext(11 15){$6$}
\htext(9 15){$6$}\htext(7 15){$6$}\htext(15 15){$6$}\htext(17
15){$6$} \htext(19 15){$6$}\htext(21 15){$6$}

\htext(21.5 12.5){$8$} \htext(20.5 13.5){$7$}\htext(19.5
12.5){$7$} \htext(18.5 13.5){$8$}\htext(17.5 12.5){$8$}
\htext(16.5 13.5){$7$}\htext(15.5 12.5){$7$} \htext(14.5
13.5){$8$}\htext(13.5 12.5){$8$} \htext(12.5 13.5){$7$}\htext(11.5
12.5){$7$} \htext(10.5 13.5){$8$}\htext(9.5 12.5){$8$} \htext(8.5
13.5){$7$} \htext(7.5 12.5){$7$} \htext(6.5 13.5){$8$}\htext(5.5
12.5){$8$}

\htext(13 11){$6$}\htext(11 11){$6$} \htext(9 11){$6$}\htext(7
11){$6$} \htext(5 11){$6$}\htext(15 11){$6$}\htext(17 11){$6$}
\htext(19 11){$6$}\htext(21 11){$6$} \htext(13 9){$5$}\htext(11
9){$5$} \htext(9 9){$5$}\htext(7 9){$5$} \htext(5 9){$5$} \htext(3
9){$5$} \htext(15 9){$5$}\htext(17 9){$5$} \htext(19
9){$5$}\htext(21 9){$5$} \htext(13 7){$4$}\htext(11 7){$4$}
\htext(9 7){$4$}\htext(7 7){$4$} \htext(5 7){$4$} \htext(3
7){$4$}\htext(1 7){$4$} \htext(15 7){$4$}\htext(17 7){$4$}
\htext(19 7){$4$}\htext(21 7){$4$} \htext(13 5){$3$}\htext(11
5){$3$} \htext(9 5){$3$}\htext(7 5){$3$} \htext(5 5){$3$} \htext(3
5){$3$}\htext(1 5){$3$}\htext(15 5){$3$}\htext(17 5){$3$}
\htext(19 5){$3$}\htext(21 5){$3$} \htext(13 3){$2$}\htext(11
3){$2$} \htext(9 3){$2$}\htext(7 3){$2$} \htext(5 3){$2$} \htext(3
3){$2$}\htext(1 3){$2$}\htext(15 3){$2$}\htext(17 3){$2$}
\htext(19 3){$2$}\htext(21 3){$2$}

\htext(21.5 0.5){$0$}\htext(20.5 1.5){$1$} \htext(19.5
0.5){$1$}\htext(18.5 1.5){$0$}\htext(17.5 0.5){$0$}\htext(16.5
1.5){$1$}\htext(15.5 0.5){$1$}\htext(14.5 1.5){$0$}\htext(13.5
0.5){$0$} \htext(12.5 1.5){$1$}\htext(11.5 0.5){$1$}\htext(10.5
1.5){$0$}\htext(9.5 0.5){$0$} \htext(8.5 1.5){$1$}\htext(7.5
0.5){$1$}\htext(6.5 1.5){$0$}\htext(5.5 0.5){$0$} \htext(4.5
1.5){$1$}\htext(3.5 0.5){$1$}\htext(2.5 1.5){$0$}\htext(1.5
0.5){$0$} \htext(0.5 1.5){$1$}

\esegment
\end{texdraw}}
\savebox{\tmpfigd}{\begin{texdraw}\fontsize{6}{6}\selectfont
\drawdim em \setunitscale 0.9

\rlvec(2 0)\rlvec(0 2)\rlvec(2 0)\rlvec(0 4)\rlvec(-2 0)\rlvec(0
-2)\rlvec(-2 0)\rlvec(0 -4) \move(2 2)\rlvec(0 2)\move(0
2)\rlvec(2 0)\move(2 4)\rlvec(2 0)\move(2 6)\rlvec(2 0)

\htext(1 1){$3$}\htext(1 3){$4$} \htext(3 3){$4$}\htext(3 5){$5$}
\end{texdraw}}
\savebox{\tmpfige}{\begin{texdraw} \fontsize{6}{6}\selectfont
\drawdim em \setunitscale 0.9

\rlvec(2 0)\rlvec(0 2)\rlvec(2 0)\rlvec(0 2)\rlvec(2 0)\rlvec(0 4)
\rlvec(-2 0)\rlvec(0 -2)\rlvec(-2 0)\rlvec(0 -4)\rlvec(-2
0)\rlvec(0 -2)\move(4 4)\rlvec(0 2)\move(2 4)\rlvec(2 0)\move(4
6)\rlvec(2 0)

\htext(1 1){$6$}\htext(3 3){$5$} \htext(3 5){$4$}\htext(5
5){$4$}\htext(5 7){$3$}
\end{texdraw}}
\savebox{\tmpfigf}{\begin{texdraw} \fontsize{6}{6}\selectfont
\drawdim em \setunitscale 0.9

\rlvec(2 0)\rlvec(0 6)\rlvec(-2 0)\rlvec(0 -2)\rlvec(-2 0)\rlvec(0
-2)\rlvec(2 0)\rlvec(0 -2) \move(0 2)\rlvec(0 2)\move(0 2)\rlvec(2
0)\move(0 4)\rlvec(2 0)

\htext(1 1){$5$}\htext(1 3){$4$}\htext(1 5){$3$} \htext(-1 3){$4$}
\end{texdraw}}
\savebox{\tmpfigg}{\begin{texdraw} \fontsize{6}{6}\selectfont
\drawdim em \setunitscale 0.9

\rlvec(2 0)\rlvec(0 8)\rlvec(-2 0)\rlvec(0 -2)\rlvec(-2 0)\rlvec(0
-2)\rlvec(2 0)\rlvec(0 -4) \move(0 2)\rlvec(2 0)\move(0 4)\rlvec(2
0)\move(0 6)\rlvec(2 0)\move(0 4)\rlvec(0 2)

\htext(1 1){$6$}\htext(1 3){$5$}\htext(1 5){$4$}\htext(1 7){$3$}
\htext(-1 5){$4$}
\end{texdraw}}

\vskip 3mm

\begin{example}
If $\frak g=D_8$, $\lambda=\omega_5+\omega_6$ and
$$Y=\raisebox{-0.4\height}{\usebox{\tmpfiga}}\in F(\la),$$
then we have
$$\aligned &Y_{\omega_5}(7,8;3,4)=\raisebox{-0.4\height}{\usebox{\tmpfigd}}\,,\,\,\,\,\,
Y_{\omega_6}(7,8;3,4)=\raisebox{-0.4\height}{\usebox{\tmpfige}}\,,\,\,\,\,\,\\
&|Y_{\omega_5}(7,8;3,4)|=\raisebox{-0.4\height}{\usebox{\tmpfigf}}\,\,\,\,\,\
\text{and}\,\,\,\,\,
Y_{\omega_6}^t(7,8;3,4)=\raisebox{-0.4\height}{\usebox{\tmpfigg}}\,.
\endaligned$$
Here, the shaded parts represent $L_{\omega_5}(7,8;3,4)$ and
$L_{\omega_6}(7,8;3,4)$.
\end{example}

\vskip 3mm

\begin{thm}\label{thm:Dn}

Let $\la \in P^{+}$ be a dominant integral weight for $\frak g = D_n$
and write
$$\lambda=\omega_{i_1}+\cdots+\omega_{i_t}+b_1\la_{n-1}+b_2\la_n,$$
where $1\le i_1\le \cdots\le i_t\le n\!+\!1$ and $(b_1,b_2)=(1,0)$
or $(0,1)$.

Define $Y(\la)$ to be the set of all reduced proper Young walls
in $F(\la)$ satisfying the following conditions {\rm :}

{\bf (Y1)} For each $k = 1, \cdots, t$, we have
$Y_{\omega_{i_k}}^{+} \subset |Y_{\omega_{i_k}}^{-}|$.

{\bf (Y2)} For each $k=1, \cdots, t-1$, we have
$$\text{$Y^{\omega_{i_k}}
\subset Y^{\omega_{i_{k+1}}}$ in ${\overline
Y}_{\omega_{i_k}+\omega_{i_{k+1}}}$, $Y^{\omega_{i_t}} \subset
Y^{\la_{n-1}}$ in ${\overline Y}_{\omega_{i_t}+\la_{n-1}}$,
$Y^{\omega_{i_t}} \subset Y^{\la_{n}}$ in ${\overline
Y}_{\omega_{i_t}+\la_{n}}$.}$$

{\bf (Y3)} For each $k=1, \cdots, t-1$, we have
$$|Y_{(\omega_{i_k},\omega_{i_{k+1}})}^{-}| \subset
Y_{(\omega_{i_k},\omega_{i_{k+1}})}^{+}, \quad
|Y_{(\omega_{i_t},\la_{n-1})}^{-}| \subset
Y_{(\omega_{i_t},\la_{n-1})}^{+}, \quad
|Y_{(\omega_{i_t},\la_n)}^{-}| \subset
Y_{(\omega_{i_t},\la_n)}^{+}.$$

{\bf (Y4)} For each $k=1, \cdots, t-1$, if ${\overline
Y_{\omega_{i_k}+\omega_{i_{k+1}}}}$, ${\overline
Y_{\omega_{i_t}+\la_{n-1}}}$ or ${\overline
Y_{\omega_{i_t}+\la_n}}$ satisfies {\bf (C1)},

\qquad\hskip 2mm  then we have
$$\aligned
&Y_{\omega_{i_k}}^{+}(a;p,q) \subset
|Y_{\omega_{i_k}}^{-}(a;p,q)|, \quad
Y_{\omega_{i_{k+1}}}^{+}(a;p,q) \subset
|Y_{\omega_{i_{k+1}}}^{-}(a;p,q)|,\\
&Y_{\omega_{i_t}}^{+}(a;p,q) \subset
|Y_{\omega_{i_t}}^{-}(a;p,q)|.
\endaligned$$

{\bf (Y5)} For each $k=1, \cdots, t-1$, if ${\overline
Y_{\omega_{i_k}+\omega_{i_{k+1}}}}$, ${\overline
Y_{\omega_{i_t}+\la_{n-1}}}$ or ${\overline
Y_{\omega_{i_t}+\la_n}}$ satisfies {\bf (C2)},

\qquad\hskip 2mm then we have
$$\aligned
&|Y_{\omega_{i_k}}(n-1,n;p,q)| \subset
Y_{\omega_{i_{k+1}}}^t(n-1,n;p,q),\\
&|Y_{\omega_{i_t}}(n-1,n;p,q)| \subset
Y_{\la_{n-1}}^t(n-1,n;p,q),\\
&|Y_{\omega_{i_t}}(n-1,n;p,q)| \subset Y_{\la_n}^t(n-1,n;p,q).
\endaligned$$

Then there exists an isomorphism of $U_q(D_n)$-crystals
\begin{equation}
Y(\lambda) \stackrel{\sim} \longrightarrow B(\la) \quad
\text{given by} \ \ H_{\la} \longmapsto u_{\la},
\end{equation}
where $u_\la$ is the highest weight vector in $B(\la)$.
\end{thm}

\vskip 3mm
\begin{example}
Let $\frak g=D_4$ and $\la = \omega_3 + \la_4$. Then, in the
following picture, the first Young wall belongs to $Y(\la)$, but
the other ones don't. They do not satisfy the conditions {\bf
(Y1)}, {\bf (Y3)}, {\bf (Y4)} and {\bf (Y5)}, respectively.

\vskip 3mm

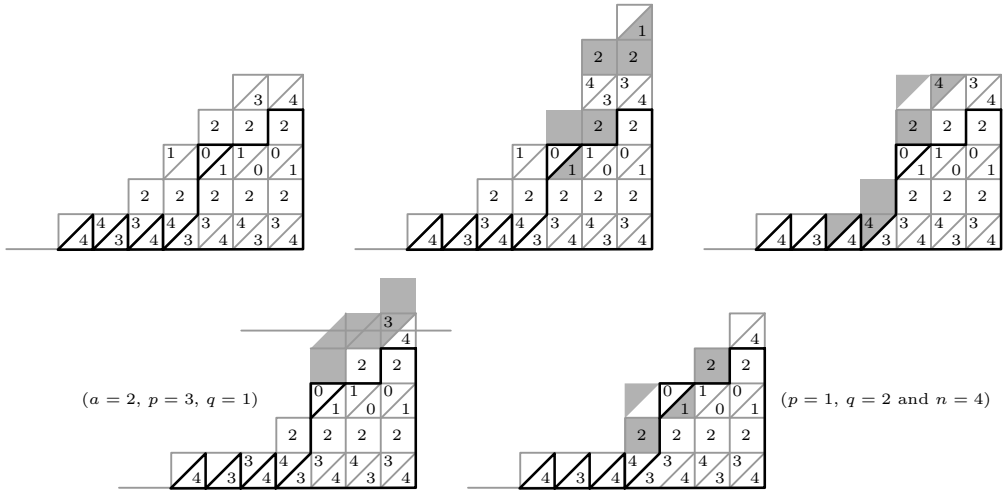
\begin{center}
\begin{texdraw}
\fontsize{6}{6}\selectfont \drawdim em \setunitscale 0.9
\nc{\dtri}{ \bsegment \lvec(-2 0) \lvec(-2 2)\lvec(0 2)\lvec(0
0)\ifill f:0.7 \esegment }

\nc{\dtrii}{ \bsegment \lvec(-2 0)\lvec(0 2)\lvec(0 0)\ifill f:0.7
\esegment }

\nc{\dtriii}{ \bsegment \lvec(-2 0) \lvec(-2 -2)\lvec(0 0)\ifill
f:0.7 \esegment }

\bsegment

\setgray 0.6 \move(-3 0) \lvec(14 0)\lvec(14 10) \move(14
10)\rlvec(-2 0)\rlvec(0 -10)\move(12 10)\rlvec(-2 0)\rlvec(0 -10)
\move(14 8)\rlvec(-6 0)\rlvec(0 -8)

\move(14 6)\rlvec(-8 0)\rlvec(0 -6)\move(14 4)\rlvec(-10
0)\rlvec(0 -4)\move(14 2)\rlvec(-12 0)\rlvec(0 -2)

\move(14 10)\rlvec(-2 -2)\move(12 10)\rlvec(-2 -2) \move(14
2)\rlvec(-2 -2)\move(12 2)\rlvec(-2 -2)\move(10 2)\rlvec(-2
-2)\move(8 2)\rlvec(-2 -2)\move(6 2)\rlvec(-2 -2)\move(4
2)\rlvec(-2 -2) \move(2 2)\rlvec(-2 -2)\move(14 6)\rlvec(-2
-2)\move(12 6)\rlvec(-2 -2)\move(10 6)\rlvec(-2 -2) \move(8
6)\rlvec(-2 -2)

\move(0 0)\rlvec(0 2)\rlvec(2 0)

\move(14 8)\linewd 0.15 \setgray 0 \rlvec(-2 0)\rlvec(0 -2)
\rlvec(-4 0)\rlvec(0 -4)\rlvec(-2 -2)\rlvec(0 2)\rlvec(-2
-2)\rlvec(0 2)\rlvec(-2 -2)\rlvec(0 2)\rlvec(-2 -2)\lvec(14
0)\lvec(14 8)\move(10 6)\rlvec(-2 -2)

\htext(13.5 8.5){$4$}\htext(11.5 8.5){$3$} \htext(13 7){$2$}
\htext(11 7){$2$}\htext(9 7){$2$} \htext(13.5 4.5){$1$}
\htext(12.5 5.5){$0$} \htext(11.5 4.5){$0$} \htext(10.5 5.5){$1$}
\htext(9.5 4.5){$1$} \htext(8.5 5.5){$0$} \htext(6.5 5.5){$1$}
\htext(13 3){$2$} \htext(11 3){$2$} \htext(9 3){$2$} \htext(7
3){$2$} \htext(5 3){$2$}

\htext(13.5 0.5){$4$}\htext(12.5 1.5){$3$} \htext(11.5
0.5){$3$}\htext(10.5 1.5){$4$} \htext(9.5 0.5){$4$}\htext(8.5
1.5){$3$}\htext(7.5 0.5){$3$} \htext(6.5 1.5){$4$} \htext(5.5
0.5){$4$}\htext(4.5 1.5){$3$}\htext(3.5 0.5){$3$} \htext(2.5
1.5){$4$} \htext(1.5 0.5){$4$}

\esegment

\move(20 0) \bsegment \move(12 10)\dtri \move(14 10)\dtri \move(12
6)\dtri\move(10 6)\dtri \move(14 12)\dtrii\move(10 4)\dtrii

\setgray 0.6 \move(-3 0) \lvec(14 0)\lvec(14 14) \move(14
14)\rlvec(-2 0)\rlvec(0 -14) \move(14 12)\rlvec(-4 0)\rlvec(0
-12)\move(14 10)\rlvec(-4 0) \move(14 8)\rlvec(-4 0)

\move(14 6)\rlvec(-8 0)\rlvec(0 -6)\move(14 4)\rlvec(-10
0)\rlvec(0 -4)\move(14 2)\rlvec(-10 0)\rlvec(-2 -2)\rlvec(0
2)\rlvec(-2 -2)\move(8 0)\rlvec(0 6)

\move(14 14)\rlvec(-2 -2)\move(14 10)\rlvec(-2 -2)\move(12
10)\rlvec(-2 -2) \move(14 2)\rlvec(-2 -2)\move(12 2)\rlvec(-2
-2)\move(10 2)\rlvec(-2 -2)\move(8 2)\rlvec(-2 -2)\move(6
2)\rlvec(-2 -2)\move(4 2)\rlvec(-2 -2) \move(2 2)\rlvec(-2
-2)\move(14 6)\rlvec(-2 -2)\move(12 6)\rlvec(-2 -2)\move(10
6)\rlvec(-2 -2) \move(8 6)\rlvec(-2 -2)

\move(0 0)\rlvec(0 2)\rlvec(4 0)

\move(14 8)\linewd 0.15 \setgray 0 \rlvec(-2 0)\rlvec(0 -2)
\rlvec(-4 0)\rlvec(0 -4)\rlvec(-2 -2)\rlvec(0 2)\rlvec(-2
-2)\rlvec(0 2)\rlvec(-2 -2)\rlvec(0 2)\rlvec(-2 -2)\lvec(14
0)\lvec(14 8)\move(10 6)\rlvec(-2 -2)

\htext(13.5 12.5){$1$} \htext(13 11){$2$}\htext(11 11){$2$}
\htext(12.5 9.5){$3$}\htext(10.5 9.5){$4$} \htext(13.5
8.5){$4$}\htext(11.5 8.5){$3$} \htext(13 7){$2$} \htext(11 7){$2$}
\htext(13.5 4.5){$1$} \htext(12.5 5.5){$0$} \htext(11.5 4.5){$0$}
\htext(10.5 5.5){$1$} \htext(9.5 4.5){$1$} \htext(8.5 5.5){$0$}
\htext(6.5 5.5){$1$} \htext(13 3){$2$} \htext(11 3){$2$} \htext(9
3){$2$} \htext(7 3){$2$} \htext(5 3){$2$}

\htext(13.5 0.5){$4$}\htext(12.5 1.5){$3$} \htext(11.5
0.5){$3$}\htext(10.5 1.5){$4$} \htext(9.5 0.5){$4$}\htext(8.5
1.5){$3$}\htext(7.5 0.5){$3$} \htext(6.5 1.5){$4$} \htext(5.5
0.5){$4$}\htext(4.5 1.5){$3$}\htext(3.5 0.5){$3$} \htext(1.5
0.5){$4$}

\esegment

\move(40 0) \bsegment \move(10 6)\dtri \move(10 10)\dtriii
\move(12 10)\dtriii \move(8 2)\dtriii \move(6 2)\dtriii \move(8
2)\dtri

\setgray 0.6 \move(-3 0) \lvec(14 0)\lvec(14 10) \move(14
10)\rlvec(-2 0)\rlvec(0 -10)\move(12 10)\rlvec(-2 0)\rlvec(0 -10)
\move(14 8)\rlvec(-6 0)\rlvec(0 -8)

\move(14 6)\rlvec(-6 0)\move(14 4)\rlvec(-6 0)\move(14 2)\rlvec(-8
0)\rlvec(0 -2)

\move(14 10)\rlvec(-2 -2)\move(12 10)\rlvec(-2 -2) \move(14
2)\rlvec(-2 -2)\move(12 2)\rlvec(-2 -2)\move(10 2)\rlvec(-2
-2)\move(8 2)\rlvec(-2 -2)\move(6 2)\rlvec(-2 -2)\move(4
2)\rlvec(-2 -2) \move(2 2)\rlvec(-2 -2)\move(14 6)\rlvec(-2
-2)\move(12 6)\rlvec(-2 -2)\move(10 6)\rlvec(-2 -2)

\move(0 0)\rlvec(0 2)\rlvec(6 0)

\move(14 8)\linewd 0.15 \setgray 0 \rlvec(-2 0)\rlvec(0 -2)
\rlvec(-4 0)\rlvec(0 -4)\rlvec(-2 -2)\rlvec(0 2)\rlvec(-2
-2)\rlvec(0 2)\rlvec(-2 -2)\rlvec(0 2)\rlvec(-2 -2)\lvec(14
0)\lvec(14 8)\move(10 6)\rlvec(-2 -2)

\htext(13.5 8.5){$4$}\htext(12.5 9.5){$3$} \htext(10.5
9.5){$4$}\htext(13 7){$2$} \htext(11 7){$2$}\htext(9 7){$2$}
\htext(13.5 4.5){$1$} \htext(12.5 5.5){$0$} \htext(11.5 4.5){$0$}
\htext(10.5 5.5){$1$} \htext(9.5 4.5){$1$} \htext(8.5 5.5){$0$}
\htext(13 3){$2$} \htext(11 3){$2$} \htext(9 3){$2$}

\htext(13.5 0.5){$4$}\htext(12.5 1.5){$3$} \htext(11.5
0.5){$3$}\htext(10.5 1.5){$4$} \htext(9.5 0.5){$4$}\htext(8.5
1.5){$3$}\htext(7.5 0.5){$3$} \htext(6.5 1.5){$4$} \htext(5.5
0.5){$4$}\htext(3.5 0.5){$3$} \htext(1.5 0.5){$4$}

\esegment

\end{texdraw}
\end{center}

\begin{center}
\begin{texdraw}
\fontsize{6}{6}\selectfont \drawdim em \setunitscale 0.9
\nc{\dtri}{ \bsegment \lvec(-2 0) \lvec(-2 2)\lvec(0 2)\lvec(0
0)\ifill f:0.7 \esegment }

\nc{\dtrii}{ \bsegment \lvec(-2 0) \lvec(-2 -2)\lvec(0 0)\ifill
f:0.7 \esegment }

\nc{\dtriii}{ \bsegment \lvec(-2 -2) \lvec(0 -2)\lvec(0 0)\ifill
f:0.7 \esegment }

\bsegment \move(14 10)\dtri \move(12 8)\dtri \move(10 6)\dtri
\move(14 10)\dtrii \move(10 10)\dtriii

\setgray 0.6 \move(-3 0) \lvec(14 0)\lvec(14 10) \move(14
10)\rlvec(-2 0)\rlvec(0 -10)\move(14 8)\rlvec(-4 0)\rlvec(0 -8)

\move(14 6)\rlvec(-6 0)\rlvec(0 -6)\move(14 4)\rlvec(-8 0)\rlvec(0
-4)\move(14 2)\rlvec(-10 0)\rlvec(0 -2)

\move(14 10)\rlvec(-2 -2)\move(14 2)\rlvec(-2 -2)\move(12
2)\rlvec(-2 -2)\move(10 2)\rlvec(-2 -2)\move(8 2)\rlvec(-2
-2)\move(6 2)\rlvec(-2 -2)\move(4 2)\rlvec(-2 -2) \move(2
2)\rlvec(-2 -2)\move(14 6)\rlvec(-2 -2)\move(12 6)\rlvec(-2
-2)\move(10 6)\rlvec(-2 -2)\move(12 10)\rlvec(-2 -2)\move(10
10)\rlvec(0 -2)\rlvec(-2 0)

\move(4 9)\rlvec(12 0)

\move(0 0)\rlvec(0 2)\rlvec(4 0)

\move(14 8)\linewd 0.15 \setgray 0 \rlvec(-2 0)\rlvec(0 -2)
\rlvec(-4 0)\rlvec(0 -4)\rlvec(-2 -2)\rlvec(0 2)\rlvec(-2
-2)\rlvec(0 2)\rlvec(-2 -2)\rlvec(0 2)\rlvec(-2 -2)\lvec(14
0)\lvec(14 8)\move(10 6)\rlvec(-2 -2)

\htext(13.5 8.5){$4$}\htext(12.5 9.5){$3$}

\htext(13 7){$2$} \htext(11 7){$2$} \htext(13.5 4.5){$1$}
\htext(12.5 5.5){$0$} \htext(11.5 4.5){$0$} \htext(10.5 5.5){$1$}
\htext(9.5 4.5){$1$} \htext(8.5 5.5){$0$} \htext(13 3){$2$}
\htext(11 3){$2$} \htext(9 3){$2$} \htext(7 3){$2$}

\htext(13.5 0.5){$4$}\htext(12.5 1.5){$3$} \htext(11.5
0.5){$3$}\htext(10.5 1.5){$4$} \htext(9.5 0.5){$4$}\htext(8.5
1.5){$3$}\htext(7.5 0.5){$3$} \htext(6.5 1.5){$4$} \htext(5.5
0.5){$4$}\htext(4.5 1.5){$3$}\htext(3.5 0.5){$3$} \htext(1.5
0.5){$4$}

\htext(0 5){($a=2$, $p=3$, $q=1$)}

\esegment

\move(20 0) \bsegment \move(10 6)\dtriii \move(12 6)\dtri \move(8
2)\dtri \move(8 6)\dtrii

\setgray 0.6 \move(-3 0) \lvec(14 0)\lvec(14 10) \move(14
10)\rlvec(-2 0)\rlvec(0 -10) \move(14 8)\rlvec(-4 0)\rlvec(0 -8)

\move(14 6)\rlvec(-6 0)\move(14 4)\rlvec(-8 0)\rlvec(0 -4)\move(14
2)\rlvec(-8 0)\rlvec(0 -2)

\move(14 10)\rlvec(-2 -2)\move(14 2)\rlvec(-2 -2)\move(12
2)\rlvec(-2 -2)\move(10 2)\rlvec(-2 -2)\move(8 2)\rlvec(-2
-2)\move(6 2)\rlvec(-2 -2)\move(4 2)\rlvec(-2 -2) \move(2
2)\rlvec(-2 -2)\move(14 6)\rlvec(-2 -2)\move(12 6)\rlvec(-2
-2)\move(10 6)\rlvec(-2 -2)

\move(0 0)\rlvec(0 2)\rlvec(6 0)

\move(14 8)\linewd 0.15 \setgray 0 \rlvec(-2 0)\rlvec(0 -2)
\rlvec(-4 0)\rlvec(0 -4)\rlvec(-2 -2)\rlvec(0 2)\rlvec(-2
-2)\rlvec(0 2)\rlvec(-2 -2)\rlvec(0 2)\rlvec(-2 -2)\lvec(14
0)\lvec(14 8)\move(10 6)\rlvec(-2 -2)

\htext(13.5 8.5){$4$} \htext(13 7){$2$} \htext(11 7){$2$}
\htext(13.5 4.5){$1$} \htext(12.5 5.5){$0$} \htext(11.5 4.5){$0$}
\htext(10.5 5.5){$1$} \htext(9.5 4.5){$1$} \htext(8.5 5.5){$0$}
\htext(13 3){$2$} \htext(11 3){$2$} \htext(9 3){$2$}\htext(7
3){$2$}

\htext(13.5 0.5){$4$}\htext(12.5 1.5){$3$} \htext(11.5
0.5){$3$}\htext(10.5 1.5){$4$} \htext(9.5 0.5){$4$}\htext(8.5
1.5){$3$}\htext(7.5 0.5){$3$} \htext(6.5 1.5){$4$} \htext(5.5
0.5){$4$}\htext(3.5 0.5){$3$} \htext(1.5 0.5){$4$}

\htext(21 5){($p=1$, $q=2$ and $n=4$)} \esegment
\end{texdraw}
\end{center}

\noindent Here, the shaded parts represent $L_{\omega_3}^{\pm}$,
$L_{(\omega_3,\la_4)}^{\pm}$, $L_{\omega_3}^{\pm}(2;3,1)$, and
$L_{\omega_3}(3,4;1,2)$ and $L_{\la_4}(3,4;1,2)$ in the second,
third, fourth and fifth Young walls, respectively.
\end{example}

\vskip 1cm

\section{The Proof of Main Theorem}

In this section, we will give a proof of our main theorems. In
fact, we will only prove the case $\frak g = B_n$ because the
remaining cases can be proved in a similar manner. Observe that it
suffices to prove the following statements:

(1) For all $i=1, \cdots, n$, we have
$$\tilde e_i Y(\la) \subset Y(\la) \cup \{0\}, \quad
\tilde f_i Y(\la) \subset Y(\la) \cup \{0\}.$$

(2) If $Y \in Y(\la)$ satisfies $\tilde e_i Y =0$ for all
$i=1,\cdots, n$, then $Y=H_{\la}$.

\vskip 5mm

\noindent{\bf The Proof of Theorem 3.9 :} \hskip 3mm  We first
prove the statement (1). Let $Y\in Y(\lambda)$ and suppose that
$\fit Y \neq 0$ but $\fit Y \notin Y(\la)$ for some  $i\in I$.
Then $\fit Y$ would violate at least one of the conditions {\bf
(Y1)}-{\bf(Y4)}.

\vskip 3mm {\bf (Case 1)} Suppose $\fit Y$ does not satisfy {\bf
(Y1)}. Then there is an $i$-admissible slot in some
$Y_{\omega_{i_k}}^+$, where an $i$-block can be added to get
$\fit Y$ such that $(\fit Y)_{\omega_{i_k}}^+ \nsubseteq (\fit
Y)_{\omega_{i_k}}^-$. Note that $i\neq n$ because $(\tilde{f_n}
Y)_{\omega_{i_k}}^{\pm}=Y_{\omega_{i_k}}^{\pm}$ for all
$k=1,\cdots,t$. For simplicity, we denote by $N_j^{\pm}$ the
number of $j$-blocks in $Y_{\omega_{i_k}}^{\pm}$. Then $N_j^+ \le
N_j^-$ for all $j=1,\cdots,n-1$ because
$Y_{\omega_{i_k}}^{+}\subset Y_{\omega_{i_k}}^{-}$. Since $Y$ is
proper, we have
$$N_{i-1}^+=N_i^+\,\,\,\,\text{and}\,\,\,\, N_{i+1}^+=N_i^{+}+1.$$
Moreover, since $Y_{\omega_{i_k}}^{+}\subset Y_{\omega_{i_k}}^{-}$
but $(\fit Y)_{\omega_{i_k}}^+ \nsubseteq (\fit
Y)_{\omega_{i_k}}^-$, we can deduce $N_i^-=N_i^+=N_{i-1}^+\le
N_{i-1}^-$. Observe that ${\overline H}_{\omega_{i_k}}$ is of
staircase shape below the $i_k$-row. Then we have $N_{i+1}^-\le
N_i^{-}+1$. But we know that $N_i^{-}+1=N_i^{+}+1=N_{i+1}^{+}\le
N_{i+1}^{-}$. Therefore, $N_{i+1}^-=N_i^{-}+1$ and $Y$ must have
the following form:

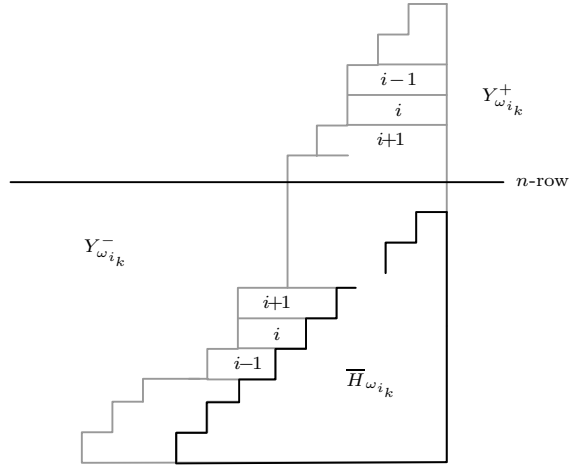
\savebox{\tmpfiga}{\begin{texdraw} \fontsize{7}{7}\selectfont
\drawdim em \setunitscale 0.9 \nc{\dtri}{ \bsegment \lvec(-2 0)
\lvec(-2 2)\lvec(0 2)\lvec(0 0)\ifill f:0.7 \esegment }

\bsegment

\setgray 0.6 \move(0 0) \lvec(0 17) \move(0 17)\rlvec(-2
0)\rlvec(0 -1.6)\rlvec(-1.6 0)\rlvec(0 -1.6)\rlvec(-1.6 0)\rlvec(0
-3.2)\rlvec(-1.6 0)\rlvec(0 -1.6)

\move(0 13.8)\rlvec(-3.6 0)\move(0 12.2)\rlvec(-5.2 0)\move(0
10.6)\rlvec(-5.2 0)\move(-5.2 9)\rlvec(-3.2 0)\rlvec(0 -7)

\move(0 6)\rlvec(-1.6 0)\rlvec(0 -1.6)\rlvec(-1.6 0)\rlvec(0
-1.6)\move(-4.8 2)\rlvec(-1 0)\rlvec(0 -1.6)\rlvec(-1.6 0)\rlvec(0
-1.6)

\move(-5.8 2)\rlvec(-5.2 0)\rlvec(0 -3.2)\rlvec(3.6 0)\move(-11
0.4)\rlvec(3.6 0)\move(-11 -1.2)\rlvec(-1.6 0)\rlvec(0
-1.6)\rlvec(3.6 0)\rlvec(0 1.6)\move(-12.6 -2.8)\rlvec(-1 0)

\move(-16 -4)\rlvec(-1.6 0)\rlvec(0 -1.6)\rlvec(-1.6 0)\rlvec(0
-1.6)\lvec(0 -7.2)\lvec(0 0)

\move(-16 -4)\rlvec(0 1.2)\rlvec(3 0)

\move(-11 -4)\rlvec(-1.6 0)\rlvec(0 -1.6)\rlvec(-1.6 0)\rlvec(0
-1.6)

\htext(-2.6 13){$i\!-\!1$}\htext(-2.6 11.4){$i$}\htext(-3
9.8){$i\!\!+\!\!1$}

\htext(-9 1.2){$i\!\!+\!\!1$} \htext(-9 -0.4){$i$} \htext(-10.5
-2){$i\!\!-\!\!1$} \htext(-4 -3){${\overline H}_{\omega_{i_k}}$}

\setgray 0 \move(0 6)\rlvec(-1.6 0)\rlvec(0 -1.6)\rlvec(-1.6
0)\rlvec(0 -1.6)\move(-4.8 2)\rlvec(-1 0)\rlvec(0 -1.6)\rlvec(-1.6
0)\rlvec(0 -1.6)\rlvec(-1.6 0)\rlvec(0 -1.6)\rlvec(-1.9 0)\rlvec(0
-1.2)\rlvec(-1.7 0)\rlvec(0 -1.6)\rlvec(-1.6 0)\rlvec(0
-1.6)\lvec(0 -7.2)\lvec(0 6)

\move(3 7.6)\rlvec(-26 0)\htext(5 7.6){$n$-row}\htext(3
12){$Y_{\omega_{i_k}}^+$}\htext(-18 4){$Y_{\omega_{i_k}}^-$}
\esegment
\end{texdraw}}

\vskip 3mm
\begin{center}
$\raisebox{-0.4\height}{\usebox{\tmpfiga}}$

\end{center}

\vskip 2mm That is, there is another $i$-admissible slot in
$Y_{\omega_{i_k}}^-$. Then, by the tensor product rule for the
Kashiwara operators, $\fit$ would have acted on the $i$-admissible
slot in $Y_{\omega_{i_k}}^-$, not on the one in
$Y_{\omega_{i_k}}^+$, which is a contradiction. Hence, $\fit Y$
must satisfy the condition {\bf (Y1)}.


%
\savebox{\tmpfiga}{\begin{texdraw} \fontsize{6}{6}\selectfont
\drawdim em \setunitscale 0.9 \nc{\dtri}{ \bsegment \lvec(-2
0)\lvec(-2 1.25)\lvec(0 1.25)\lvec(0 0)\ifill f:0.7 \esegment }

\rlvec(4 0)\rlvec(0 6)\rlvec(-2 0)\rlvec(0 -4)\rlvec(-2 0)\rlvec(0
-2)\move(2 0)\rlvec(0 2)\move(2 2)\rlvec(2 0)\move(2 4)\rlvec(2 0)

\htext(1 1){$i\!\!+\!\!1$}\htext(3 1){$i\!\!+\!\!1$}\htext(3
3){$i$}\htext(3 5){$i\!\!-\!\!1$}
\end{texdraw}}
\savebox{\tmpfigb}{\begin{texdraw} \fontsize{6}{6}\selectfont
\drawdim em \setunitscale 0.9 \nc{\dtri}{ \bsegment \lvec(-2
0)\lvec(-2 1.25)\lvec(0 1.25)\lvec(0 0)\ifill f:0.7 \esegment }

\rlvec(2 0)\rlvec(0 4)\rlvec(2 0)\rlvec(0 -2)\rlvec(-2 0)

\htext(3 3){$i\!\!+\!\!1$}
\end{texdraw}}
\savebox{\tmpfigc}{\begin{texdraw} \fontsize{6}{6}\selectfont
\drawdim em \setunitscale 0.9 \nc{\dtri}{ \bsegment \lvec(-2
0)\lvec(-2 1.25)\lvec(0 1.25)\lvec(0 0)\ifill f:0.7 \esegment }

\rlvec(2 0)\rlvec(0 6)\rlvec(2 0)\rlvec(0 -4)\rlvec(-4 0)\rlvec(0
-2)\move(2 4)\rlvec(2 0)

\htext(1 1){$i\!\!-\!\!1$}\htext(3 3){$i$}\htext(3
5){$i\!\!+\!\!1$}
\end{texdraw}}
\savebox{\tmpfigd}{\begin{texdraw} \fontsize{6}{6}\selectfont
\drawdim em \setunitscale 0.9 \nc{\dtri}{ \bsegment \lvec(-2
0)\lvec(-2 1.25)\lvec(0 1.25)\lvec(0 0)\ifill f:0.7 \esegment }

\rlvec(4 0)\rlvec(0 6)\rlvec(-2 0)\rlvec(0 -4)\rlvec(-2 0)\rlvec(0
-2)\move(2 0)\rlvec(0 2)\move(2 2)\rlvec(2 0)\move(2 4)\rlvec(2 0)

\htext(1 1){$i\!\!-\!\!1$}\htext(3 1){$i\!\!-\!\!1$}\htext(3
3){$i$}\htext(3 5){$i\!\!+\!\!1$}
\end{texdraw}}

\vskip 3mm {\bf (Case 2)} Suppose $\fit Y$ does not satisfy {\bf
(Y2)}. Then there exists an $i$-admissible slot in some
$Y^{\omega_{i_k}}$ (or $Y^{\omega_{i_t}}$), where an $i$-block can
be added to get $\fit Y$ such that $(\fit
Y)^{\omega_{i_k}}\nsubseteq (\fit Y)^{\omega_{i_{k+1}}}$ (or
$(\fit Y)^{\omega_{i_t}}\nsubseteq (\fit Y)^{\lambda_n}$). If
$1\leq i \leq n-1$, since $Y$ is proper, $Y$ has a subwall of the
form $\raisebox{-0.4\height}{\usebox{\tmpfiga}}$,
$\raisebox{-0.4\height}{\usebox{\tmpfigb}}$,
$\raisebox{-0.4\height}{\usebox{\tmpfigc}}$ or
$\raisebox{-0.4\height}{\usebox{\tmpfigd}}$ in $Y^{\omega_{i_k}}$
(or $Y^{\omega_{i_t}}$). Since $Y$ satisfies {\bf (Y2)} and just
adding an $i$-block to $ Y^{\omega_{i_k}}$ would violate {\bf
(Y2)}, the same part appears in $ Y^{\omega_{i_{k+1}}}$ (or $
Y^{\lambda_n}$), as is shown in the following picture.

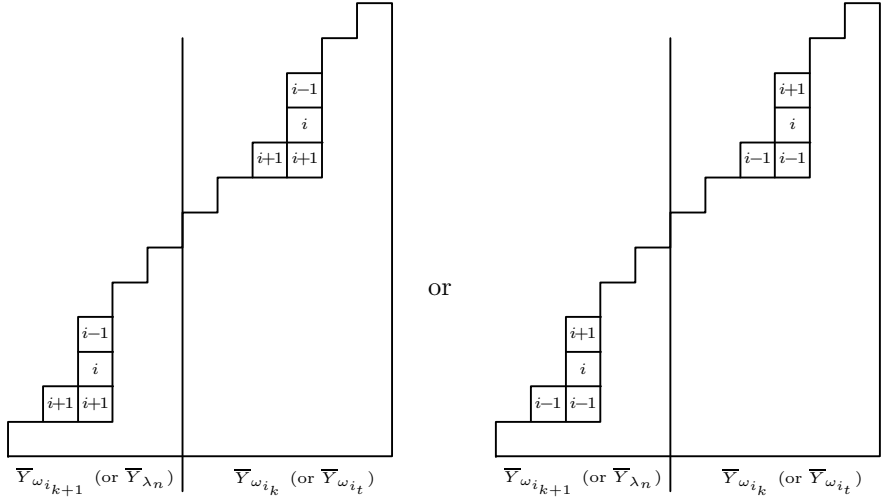
\savebox{\tmpfiga}{
\begin{texdraw}
\fontsize{6}{6}\selectfont \drawdim em \setunitscale 0.9
\nc{\dtri}{ \bsegment \lvec(-2 0) \lvec(-2 2)\lvec(0 2)\lvec(0
0)\ifill f:0.7 \esegment }

\rlvec(22 0)\rlvec(0 26)\rlvec(-2 0)\rlvec(0 -2) \rlvec(-2
0)\rlvec(0 -8)\rlvec(-6 0)\rlvec(0 -2)\rlvec(-2 0)\rlvec(0
-2)\rlvec(-2 0)\rlvec(0 -2)\rlvec(-2 0)\rlvec(0 -8)\rlvec(-6
0)\rlvec(0 -2)

\move(2 2)\rlvec(0 2)\rlvec(4 0)\move(4 2)\rlvec(0 6)\rlvec(2 0)
\move(4 6)\rlvec(2 0)

\move(14 16)\rlvec(0 2)\rlvec(4 0)\move(16 16)\rlvec(0 6)\rlvec(2
0)\move(16 20)\rlvec(2 0)

\move(10 -2)\rlvec(0 26)

\htext(3 3){$i\!\!+\!\!1$}\htext(5 3){$i\!\!+\!\!1$}\htext(5 5
){$i$}\htext(5 7){$i\!\!-\!\!1$}

\htext(15 17){$i\!\!+\!\!1$}\htext(17 17){$i\!\!+\!\!1$}\htext(17
19){$i$}\htext(17 21){$i\!\!-\!\!1$}

\htext(5
-1.5){$\overline{Y}_{\omega_{i_{k+1}}}$\,\,$(\text{or}\,\,\overline{Y}_{\la_n})$}
\htext(17
-1.5){$\overline{Y}_{\omega_{i_k}}$\,\,$(\text{or}\,\,\overline{Y}_{\omega_{i_t}})$}

\end{texdraw}}
\savebox{\tmpfigb}{
\begin{texdraw}
\fontsize{6}{6}\selectfont \drawdim em \setunitscale 0.9
\nc{\dtri}{ \bsegment \lvec(-2 0) \lvec(-2 2)\lvec(0 2)\lvec(0
0)\ifill f:0.7 \esegment }

\rlvec(22 0)\rlvec(0 26)\rlvec(-2 0)\rlvec(0 -2) \rlvec(-2
0)\rlvec(0 -8)\rlvec(-6 0)\rlvec(0 -2)\rlvec(-2 0)\rlvec(0
-2)\rlvec(-2 0)\rlvec(0 -2)\rlvec(-2 0)\rlvec(0 -8)\rlvec(-6
0)\rlvec(0 -2)

\move(2 2)\rlvec(0 2)\rlvec(4 0)\move(4 2)\rlvec(0 6)\rlvec(2 0)
\move(4 6)\rlvec(2 0)

\move(14 16)\rlvec(0 2)\rlvec(4 0)\move(16 16)\rlvec(0 6)\rlvec(2
0)\move(16 20)\rlvec(2 0)

\move(10 -2)\rlvec(0 26)

\htext(3 3){$i\!\!-\!\!1$}\htext(5 3){$i\!\!-\!\!1$}\htext(5 5
){$i$}\htext(5 7){$i\!\!+\!\!1$}

\htext(15 17){$i\!\!-\!\!1$}\htext(17 17){$i\!\!-\!\!1$}\htext(17
19){$i$}\htext(17 21){$i\!\!+\!\!1$}

\htext(5 -1.5){$\overline{Y}_{\omega_{i_{k+1}}}$\,\,(or $
\overline{Y}_{\la_n}$)} \htext(17
-1.5){$\overline{Y}_{\omega_{i_k}}$\,\,(or
$\overline{Y}_{\omega_{i_t}}$)}

\end{texdraw}}

\vskip 3mm
\begin{center}
$\raisebox{-0.4\height}{\usebox{\tmpfiga}}$ \quad or \quad
$\raisebox{-0.4\height}{\usebox{\tmpfigb}}$
\end{center}

\vskip 2mm \noindent We claim that there is no removable $i$-block
between these two parts. Then by the tensor product rule, $\fit$
would have acted on the $i$-admissible slot in
$Y^{\omega_{i_{k+1}}}$ (or $Y^{\lambda_n}$), not on the one in $
Y^{\omega_{i_k}}$ (or $ Y^{\omega_{i_t}}$), which is a
contradiction. Hence $\fit Y$ must satisfy {\bf (Y2)}. To prove
our claim, assume first that there exists a removable $i$-block in
$Y^{\omega_{i_k}}$. Then $Y$ must have the following shape:

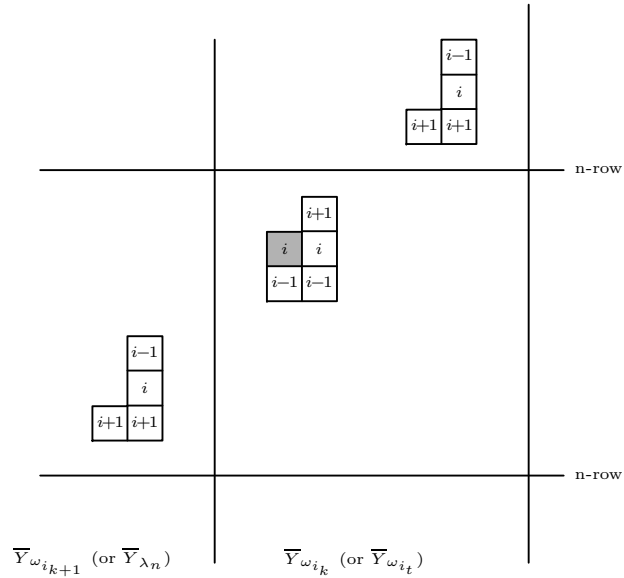
\begin{center}
\begin{texdraw}
\fontsize{6}{6}\selectfont \drawdim em \setunitscale 0.9
\nc{\dtri}{ \bsegment \lvec(-2 0) \lvec(-2 2)\lvec(0 2)\lvec(0
0)\ifill f:0.7 \esegment }

\bsegment \move(15 15)\dtri

\move(28 -2)\rlvec(0 32) \move(3 5)\rlvec(4 0)\rlvec(0 6)\rlvec(-2
0)\rlvec(0 -4)\rlvec(-2 0)\rlvec(0 -2)

\move(5 5)\rlvec(0 2)\move(5 7)\rlvec(2 0)\move(5 9)\rlvec(2 0)

\move(0 3)\rlvec(30 0)

\htext(4 6){$i\!\!+\!\!1$}\htext(6 6){$i\!\!+\!\!1$}\htext(6 8
){$i$}\htext(6 10){$i\!\!-\!\!1$}


\move(13 13)\rlvec(4 0)\rlvec(0 6)\rlvec(-2 0)\rlvec(0
-2)\rlvec(-2 0)\rlvec(0 -4)

\move(15 13)\rlvec(0 4)\move(13 15)\rlvec(4 0)\move(15 17)\rlvec(2
0)

\htext(14 14){$i\!\!-\!\!1$}\htext(16 14 ){$i\!\!-\!\!1$}\htext(14
16){$i$}\htext(16 16){$i$}\htext(16 18){$i\!\!+\!\!1$}

\move(0 20.5)\rlvec(30 0)

\move(21 22)\rlvec(4 0)\rlvec(0 6)\rlvec(-2 0)\rlvec(0
-4)\rlvec(-2 0)\rlvec(0 -2)

\move(23 22)\rlvec(0 2)\move(23 24)\rlvec(2 0)\move(23 26)\rlvec(2
0)

\htext(22 23){$i\!\!+\!\!1$}\htext(24 23){$i\!\!+\!\!1$}\htext(24
25){$i$}\htext(24 27){$i\!\!-\!\!1$}

\htext(32 3){n-row}\htext(32 20.5){n-row}

\move(10 -2)\rlvec(0 30)

\htext(3 -2){$\overline{Y}_{\omega_{i_{k+1}}}$\,\,(or
$\overline{Y}_{\la_n}$)} \htext(18
-2){$\overline{Y}_{\omega_{i_k}}$\,\,(or
$\overline{Y}_{\omega_{i_t}}$)}

\esegment
\end{texdraw}
\end{center}

\savebox{\tmpfiga}{\begin{texdraw} \fontsize{6}{6}\selectfont
\drawdim em \setunitscale 0.9 \nc{\dtri}{ \bsegment \lvec(-2
0)\lvec(-2 1.25)\lvec(0 1.25)\lvec(0 0)\ifill f:0.7 \esegment }

\rlvec(4 0)\rlvec(0 6)\rlvec(-2 0)\rlvec(0 -4)\rlvec(-2 0)\rlvec(0
-2)\move(2 0)\rlvec(0 2)\move(2 2)\rlvec(2 0)\move(2 4)\rlvec(2 0)

\htext(1 1){$i\!\!+\!\!1$}\htext(3 1){$i\!\!+\!\!1$}\htext(3
3){$i$}\htext(3 5){$i\!\!-\!\!1$}
\end{texdraw}}
\savebox{\tmpfigb}{\begin{texdraw} \fontsize{6}{6}\selectfont
\drawdim em \setunitscale 0.9 \nc{\dtri}{ \bsegment \lvec(-2
0)\lvec(-2 2)\lvec(0 2)\lvec(0 0)\ifill f:0.7 \esegment }

\bsegment\move(2 2)\dtri

\move(0 0)\rlvec(4 0)\rlvec(0 6)\rlvec(-2 0)\rlvec(0 -4)\rlvec(-2
0)\rlvec(0 -2)\move(2 0)\rlvec(0 2)\move(2 2)\rlvec(2 0)\move(2
4)\rlvec(2 0) \move(2 4)\rlvec(-2 0)\rlvec(0 -2)

\htext(1 1){$i\!\!-\!\!1$}\htext(1 3){$i$}\htext(3
1){$i\!\!-\!\!1$}\htext(3 3){$i$}\htext(3 5){$i\!\!+\!\!1$}
\esegment
\end{texdraw}}
\savebox{\tmpfigc}{\begin{texdraw} \fontsize{6}{6}\selectfont
\drawdim em \setunitscale 0.9 \nc{\dtri}{ \bsegment \lvec(-2
0)\lvec(-2 2)\lvec(0 2)\lvec(0 0)\ifill f:0.7 \esegment }

\bsegment \move(2 0)\dtri

\move(0 0)\rlvec(2 0)\rlvec(0 4)\rlvec(2 0)\rlvec(0 -2)\rlvec(-4
0)\rlvec(0 -2)

\htext(1 1){$i$}\htext(3 3){$i\!\!+\!\!1$}
\esegment
\end{texdraw}}
\savebox{\tmpfigd}{\begin{texdraw} \fontsize{6}{6}\selectfont
\drawdim em \setunitscale 0.9 \nc{\dtri}{ \bsegment \lvec(-2
0)\lvec(-2 2)\lvec(0 2)\lvec(0 0)\ifill f:0.7 \esegment }

\bsegment  \move(2 2)\dtri

\move(0 0)\rlvec(2 0)\rlvec(0 6)\rlvec(2 0)\rlvec(0 -4)\rlvec(-4
0)\rlvec(0 -2)\move(0 2)\rlvec(0 2)\rlvec(4 0)

\htext(1 1){$i\!\!-\!\!1$}\htext(1 3){$i$}\htext(3 3){$i$}\htext(3
5){$i\!\!+\!\!1$}
\esegment
\end{texdraw}}
\vskip 2mm \noindent Thus $Y^{+}_{(\omega_{i_k},
\omega_{i_{k+1}})}$ would contain the part
$\raisebox{-0.4\height}{\usebox{\tmpfiga}}$ and
$Y^{-}_{(\omega_{i_k}, \omega_{i_{k+1}})}$ would contain one of
the parts $\raisebox{-0.4\height}{\usebox{\tmpfigc}}$,
$\raisebox{-0.4\height}{\usebox{\tmpfigd}}$ or
$\raisebox{-0.4\height}{\usebox{\tmpfigb}}$\,. This implies that
$Y$ does not satisfy the condition {\bf (Y3)} or {\bf (Y4)},
which is a contradiction. Hence there is no removable $i$-block in
$Y^{\omega_{i_{k}}}.$ Next, assume that there exists a removable
$i$-block in $Y^{\omega_{i_{k+1}}}.$ Then $Y$ must have the
following shape:

\vskip 1mm
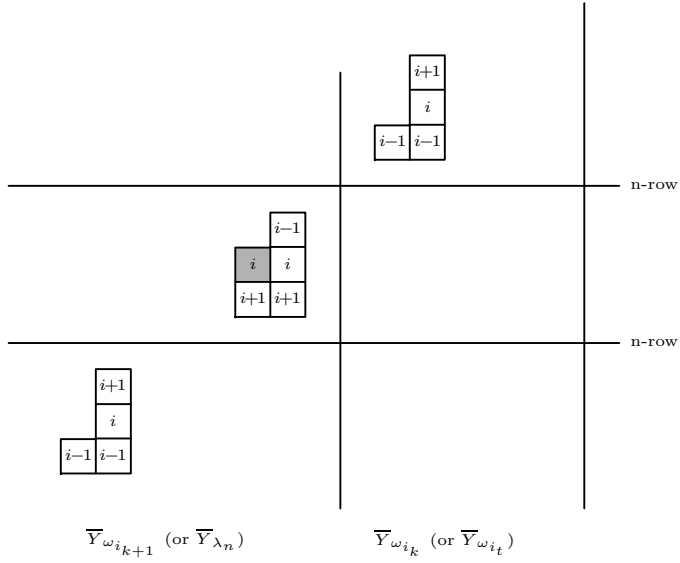
\begin{center}
\begin{texdraw}
\fontsize{6}{6}\selectfont \drawdim em \setunitscale 0.9
\nc{\dtri}{ \bsegment \lvec(-2 0) \lvec(-2 2)\lvec(0 2)\lvec(0
0)\ifill f:0.7 \esegment }

\bsegment \move(15 13)\dtri

\move(33 0)\rlvec(0 29)

\move(3 2)\rlvec(4 0)\rlvec(0 6)\rlvec(-2 0)\rlvec(0 -4)\rlvec(-2
0)\rlvec(0 -2)

\move(5 2)\rlvec(0 2)\move(5 4)\rlvec(2 0)\move(5 6)\rlvec(2 0)

\move(0 9.5)\rlvec(35 0)

\htext(4 3){$i\!\!-\!\!1$}\htext(6 3){$i\!\!-\!\!1$}\htext(6 5
){$i$}\htext(6 7){$i\!\!+\!\!1$}


\move(13 11)\rlvec(4 0)\rlvec(0 6)\rlvec(-2 0)\rlvec(0
-2)\rlvec(-2 0)\rlvec(0 -4)

\move(15 11)\rlvec(0 4)\move(13 13)\rlvec(4 0)\move(15 15)\rlvec(2
0)

\htext(14 12){$i\!\!+\!\!1$}\htext(16 12){$i\!\!+\!\!1$}\htext(14
14){$i$}\htext(16 14){$i$}\htext(16 16){$i\!\!-\!\!1$}

\move(0 18.5)\rlvec(35 0)

\move(21 20)\rlvec(4 0)\rlvec(0 6)\rlvec(-2 0)\rlvec(0
-4)\rlvec(-2 0)\rlvec(0 -2)

\move(23 20)\rlvec(0 2)\move(23 22)\rlvec(2 0)\move(23 24)\rlvec(2
0)

\htext(22 21){$i\!\!-\!\!1$}\htext(24 21){$i\!\!-\!\!1$}\htext(24
23){$i$}\htext(24 25){$i\!\!+\!\!1$}

\htext(37 9.5){n-row}\htext(37 18.5){n-row}

\move(19 0)\rlvec(0 25)

\htext(9 -2){$\overline{Y}_{\omega_{i_{k+1}}}$\,\,(or
$\overline{Y}_{\la_n}$)} \htext(25
-2){$\overline{Y}_{\omega_{i_k}}$\,\,(or
$\overline{Y}_{\omega_{i_t}}$)}

\esegment
\end{texdraw}
\end{center}

\savebox{\tmpfiga}{\begin{texdraw} \fontsize{6}{6}\selectfont
\drawdim em \setunitscale 0.9 \nc{\dtri}{ \bsegment \lvec(-2
0)\lvec(-2 1.25)\lvec(0 1.25)\lvec(0 0)\ifill f:0.7 \esegment }

\rlvec(2 0)\rlvec(0 3)\rlvec(2 0)\rlvec(0 -1)\rlvec(-4 0)\rlvec(0
-2)

\htext(1 1){$n\!\!-\!\!1$}\htext(3 2.5){$n$}
\end{texdraw}}
\savebox{\tmpfigb}{\begin{texdraw} \fontsize{6}{6}\selectfont
\drawdim em \setunitscale 0.9 \nc{\dtri}{ \bsegment \lvec(-2
0)\lvec(-2 1.25)\lvec(0 1.25)\lvec(0 0)\ifill f:0.7 \esegment }

\rlvec(2 0)\rlvec(0 6)\rlvec(2 0)\rlvec(0 -4)\rlvec(-4 0)\rlvec(0
-2)\move(2 3)\rlvec(2 0)\move(2 4)\rlvec(2 0)

\htext(1 1){$n\!\!-\!\!1$}\htext(3 2.5){$n$}\htext(3
3.5){$n$}\htext(3 5){$n\!\!-\!\!1$}
\end{texdraw}}
\savebox{\tmpfigc}{\begin{texdraw} \fontsize{6}{6}\selectfont
\drawdim em \setunitscale 0.9 \nc{\dtri}{ \bsegment \lvec(-2
0)\lvec(-2 1.25)\lvec(0 1.25)\lvec(0 0)\ifill f:0.7 \esegment }

\rlvec(4 0)\rlvec(0 4)\rlvec(-2 0)\rlvec(0 -3)\rlvec(-2 0)\rlvec(0
-1)\move(2 0)\rlvec(0 1)\move(2 2)\rlvec(2 0)\move(2 1)\rlvec(2 0)

\htext(1 0.5){$n$} \htext(3 0.5){$n$}\htext(3 1.5){$n$}\htext(3
3){$n\!\!-\!\!1$}
\end{texdraw}}
\savebox{\tmpfigd}{\begin{texdraw} \fontsize{6}{6}\selectfont
\drawdim em \setunitscale 0.9 \nc{\dtri}{ \bsegment \lvec(-2
0)\lvec(-2 1.25)\lvec(0 1.25)\lvec(0 0)\ifill f:0.7 \esegment }

\rlvec(2 0)\rlvec(0 6)\rlvec(2 0)\rlvec(0 -4)\rlvec(-4 0)\rlvec(0
-2)\move(0 2)\rlvec(4 0)\move(2 3)\rlvec(2 0)\move(2 4)\rlvec(2
0)\move(0 2)\rlvec(0 1)\rlvec(2 0)

\htext(1 1){$n\!\!-\!\!1$}\htext(1 2.5){$n$} \htext(3
2.5){$n$}\htext(3 3.5){$n$}\htext(3 5){$n\!\!-\!\!1$}
\end{texdraw}}

\vskip 2mm \noindent Now, by a similar argument as above, one can
see that $Y$ does not satisfy {\bf (Y3)} or {\bf (Y4)}, which is a
contradiction. Thus there is no removable $i$-block between these
two parts as we claimed.

\vskip 2mm If $i=n$, $Y$ has a subwall of the form

$$\aligned
&(\text{a})\,\,\quad
\raisebox{-0.4\height}{\usebox{\tmpfiga}}\,\,\,\,\text{in}\,\,\,\,
Y^{\omega_{i_k}} \,\,\,\,\,\text{and}\,\,\,\,
\raisebox{-0.4\height}{\usebox{\tmpfigb}}\,\,\,\,\text{in}\,\,\,\,
Y^{\omega_{i_{k+1}}}\,\, \,\,\text{or}\\
&(\text{b})\,\,\quad
\raisebox{-0.4\height}{\usebox{\tmpfigc}}\,\,\,\,\text{in}\,\,\,\,
Y^{\omega_{i_k}} \,\,\,\,\,\text{and}\,\,\,\,
\raisebox{-0.4\height}{\usebox{\tmpfigd}}\,\,\,\,\text{in}\,\,\,\,
Y^{\omega_{i_{k+1}}}.
\endaligned$$

\vskip 2mm \noindent The case (b) does not occur because $Y$ would
not satisfy the condition {\bf (Y3)}. In the case of (a),
$\tilde{f_n}$ would have acted on the $n$-admissible slot in
$Y^{\omega_{i_{k+1}}}$ not on the one in $Y^{\omega_{i_k}}$, which
is a contradiction. Therefore, $\fit Y$ must satisfy the condition
{\bf (Y2)}.


\vskip 3mm {\bf (Case 3)} If  $\fit Y$ does not satisfy (Y3),
then by a similar argument to (Case 1) and (Case 2), we can
derive a contradiction. Hence, $\fit Y$ must satisfy the
condition {\bf (Y3)}.

\vskip 3mm  {\bf (Case 4)}  Suppose that $\fit Y$ ($1\le i \le
n-1$) has the configuration {\bf (C1)} but does not satisfy {\bf
(Y4)}. (If $i=n$, $\tilde f_n Y$ does not have the configuration
{\bf (C1)} by the condition {\bf (Y3)}.)  Then we have the
following two possibilities:

(a) $Y$ has the configuration {\bf (C1)}, $Y$ satisfies {\bf
(Y4)}, but $\fit Y$ does not satisfy {\bf (Y4)}.

(b) $Y$ does not have the configuration {\bf (C1)}, $\fit Y$ has
the configuration {\bf (C1)}, but $\fit Y$ does not satisfy {\bf
(Y4)}.

\savebox{\tmpfiga}{\begin{texdraw} \fontsize{6}{6}\selectfont
\drawdim em \setunitscale 0.9 \nc{\dtri}{ \bsegment \lvec(-2 0)
\lvec(-2 2)\lvec(0 2)\lvec(0 0)\ifill f:0.7 \esegment }

\move(10.5 12)\dtri

\setgray 0 \move(0 0)\rlvec(24 0)\rlvec(0 29.5) \move(17.5
-2)\rlvec(0 30)

\move(3 3)\rlvec(0 3)\rlvec(2 0)\rlvec(0 4)\rlvec(2 0)\rlvec(0 -7)
\move(3 4)\rlvec(2 0)\move(5 6)\rlvec(0 -3) \move(5 6)\rlvec(2 0)
\move(5 8)\rlvec(2 0)\move(3 6)\rlvec(0 2)\rlvec(2 0)

\move(8.5 12)\rlvec(0 2)\rlvec(2 0)\rlvec(0 2)\rlvec(2 0)\rlvec(0
-4)\move(10.5 14)\rlvec(0 -2)\move(10.5 14)\rlvec(2 0)

\move(14 18)\rlvec(0 4)\rlvec(2 0)\rlvec(0 -4)\move(14 20)\rlvec(2
0)

\move(19 25)\rlvec(0 2)\rlvec(2 0)\rlvec(0 -2)

\htext(4 5){$i\!\!-\!\!1$}\htext(4 7){$\circledast$}\htext(6
7){$i$}\htext(6 9){$i\!\!+\!\!1$}

\htext(9.5 13){$i$}\htext(11.5 15){$i\!\!-\!\!1$}

\htext(15 19){$a\!\!+\!\!1$}\htext(15 21){$a$}

\htext(20 26){$a\!\!-\!\!1$}

\move(2 1.2)\rlvec(0 5)\htext(2 0.7){$p$-th}

\move(15 14)\rlvec(0 3)\htext(15 13){$q$-th}

\move(20 21)\rlvec(0 3)\htext(20 20){$p$-th}

\htext(9 -1.5){$\overline{Y}_{\omega_{i_{k+1}}}$ (or
$\overline{Y}_{\la_n}$)}

\htext(22 -1.5){$\overline{Y}_{\omega_{i_{k}}}$ (or
$\overline{Y}_{\omega_{i_t}}$)}
\end{texdraw}}
\savebox{\tmpfigb}{\begin{texdraw} \fontsize{6}{6}\selectfont
\drawdim em \setunitscale 0.9 \nc{\dtri}{ \bsegment \lvec(-2 0)
\lvec(-2 2)\lvec(0 2)\lvec(0 0)\ifill f:0.7 \esegment }

\move(16 20)\dtri

\setgray 0 \move(-1 0)\rlvec(23 0)\rlvec(0 29.5) \move(3.5
-2)\rlvec(0 30)

\move(8.5 11)\rlvec(0 3)\rlvec(2 0)\rlvec(0 4)\rlvec(2 0)\rlvec(0
-7) \move(8.5 12)\rlvec(2 0)\move(10.5 14)\rlvec(0 -3) \move(10.5
14)\rlvec(2 0) \move(10.5 16)\rlvec(2 0)\move(8.5 14)\rlvec(0
2)\rlvec(2 0)

\move(14 20)\rlvec(0 2)\rlvec(2 0)\rlvec(0 2)\rlvec(2 0)\rlvec(0
-4)\move(16 22)\rlvec(0 -2)\move(16 22)\rlvec(2 0)

\move(0.5 3)\rlvec(0 4)\rlvec(2 0)\rlvec(0 -4)\move(0.5 5)\rlvec(2
0)

\move(5 7)\rlvec(0 2)\rlvec(2 0)\rlvec(0 -2)

\htext(9.5 13){$i\!\!-\!\!1$}\htext(9.5
15){$\circledast$}\htext(11.5 15){$i$}\htext(11.5
17){$i\!\!+\!\!1$}

\htext(15 21){$i$}\htext(17 23){$i\!\!-\!\!1$}

\htext(1.5 4){$a\!\!+\!\!1$}\htext(1.5 6){$a$}

\htext(6 8){$a\!\!-\!\!1$}

\move(6 3)\rlvec(0 4)\htext(6 2.5){$p$-th}

\move(1.5 1)\rlvec(0 2)\htext(1.5 0.5){$q$-th}

\move(20 5)\rlvec(0 20)\htext(20 4.5){$q$-th}

\htext(-1 -1.5){$\overline{Y}_{\omega_{i_{k+1}}}$ (or
$\overline{Y}_{\la_n}$)}

\htext(12 -1.5){$\overline{Y}_{\omega_{i_{k}}}$ (or
$\overline{Y}_{\omega_{i_t}}$)}
\end{texdraw}}
\savebox{\tmpfigc}{\begin{texdraw} \fontsize{6}{6}\selectfont
\drawdim em \setunitscale 0.9 \nc{\dtri}{ \bsegment \lvec(-2 0)
\lvec(-2 2)\lvec(0 2)\lvec(0 0)\ifill f:0.7 \esegment }

\bsegment \move(2 0)\dtri

\move(0 0) \rlvec(0 2)\rlvec(2 0)\rlvec(0 -2)\rlvec(-2 0)

\htext(1 1){$i$}

\esegment
\end{texdraw}}

\vskip 2mm On the one hand, in the case of (a), $\fit Y$ must have
the form

\vskip 3mm
\begin{center}
$\raisebox{-0.4\height}{\usebox{\tmpfiga}}$ \qquad or \qquad
$\raisebox{-0.4\height}{\usebox{\tmpfigb}}$
\end{center}

\vskip 3mm\noindent because $Y$ satisfies {\bf (Y4)}. Then $\fit$
would have acted on $\circledast$, not on
$\raisebox{-0.2\height}{\usebox{\tmpfigc}}$\,, which is a
contradiction.

\vskip 2mm On the other hand, in the case of (b), observe that, by
adding an $i$-block to $Y$, $\fit Y$ can have the configuration
{\bf (C1)} with $a=i$ or $a=i+1$, which is shown in the following
figure:

\savebox{\tmpfiga}{\begin{texdraw} \fontsize{6}{6}\selectfont
\drawdim em \setunitscale 0.9 \nc{\dtri}{ \bsegment \lvec(-2 0)
\lvec(-2 2)\lvec(0 2)\lvec(0 0)\ifill f:0.7 \esegment }

\bsegment \move(-14 11)\dtri

\setgray 0 \move(0 0) \lvec(0 24) \move(-2 22)\rlvec(-2 0)\move(-2
20)\rlvec(-4 0)\move(-4 18)\rlvec(-2 0)\move(-4 16)\rlvec(-2
0)\move(-6 15)\rlvec(0 5)\move(-4 15)\rlvec(0 7) \move(-2
19)\rlvec(0 3)

\move(-14 13)\rlvec(-2 0)\move(-14 11)\rlvec(-2 0) \move(-14
9)\rlvec(-2 0)\move(-14 7)\rlvec(-2 0) \move(-14 6.5)\rlvec(0
6.5)\move(-16 6.5)\rlvec(0 6.5)\move(0 0)\rlvec(-18 0)\move(-14
13)\rlvec(0 1.5)\rlvec(0.5 0)

\move(-11 -2)\rlvec(0 20)

\move(-5 15)\rlvec(0 -2)\move(-5 14)\ravec(5 0)\move(0
14)\ravec(-5 0)\htext(-5 12){$p$-th}

\move(-15 4.5)\rlvec(0 2)\move(-15 5.5)\ravec(4 0)\move(-11
5.5)\ravec(-4 0)\htext(-14 4){$q$-th}

\htext(-3 21){$i$}\htext(-5 19){$i\!\!-\!\!1$}\htext(-5
17){$i\!\!-\!\!2$}

\htext(-15 12){$i$}\htext(-15 10){$i\!\!+\!\!1$} \htext(-15
8){$i\!\!+\!\!2$}

\htext(-16 -1.5){$\overline{Y}_{\omega_{i_{k+1}}}$\,\,(or
$\overline{Y}_{\la_n}$)} \htext(-5
-1.5){$\overline{Y}_{\omega_{i_k}}$\,\,(or
$\overline{Y}_{\omega_{i_t}}$)}

\htext(-10 -4){$a=i$} \esegment
\end{texdraw}}
\savebox{\tmpfigb}{\begin{texdraw} \fontsize{6}{6}\selectfont
\drawdim em \setunitscale 0.9 \nc{\dtri}{ \bsegment \lvec(-2 0)
\lvec(-2 2)\lvec(0 2)\lvec(0 0)\ifill f:0.7 \esegment }

\bsegment \move(-4 18)\dtri

\setgray 0 \move(0 0) \lvec(0 24) \move(-2 22)\rlvec(-2 0)\move(-2
20)\rlvec(-4 0)\move(-4 18)\rlvec(-2 0)\move(-4 16)\rlvec(-2
0)\move(-6 15)\rlvec(0 5)\move(-4 15)\rlvec(0 7) \move(-2
19)\rlvec(0 3)

\move(-14 13)\rlvec(-2 0)\move(-14 11)\rlvec(-2 0) \move(-14
9)\rlvec(-2 0)\move(-14 7)\rlvec(-2 0) \move(-14 6.5)\rlvec(0
6.5)\move(-16 6.5)\rlvec(0 6.5)\move(0 0)\rlvec(-18 0)\move(-14
13)\rlvec(0 1.5)\rlvec(1.5 0)

\move(-11 -2)\rlvec(0 20)

\move(-5 15)\rlvec(0 -2)\move(-5 14)\ravec(5 0)\move(0
14)\ravec(-5 0)\htext(-5 12){$p$-th}

\move(-15 4.5)\rlvec(0 2)\move(-15 5.5)\ravec(4 0)\move(-11
5.5)\ravec(-4 0)\htext(-14 4){$q$-th}

\htext(-3 21){$i\!\!+\!\!1$}\htext(-5 19){$i$}\htext(-5
17){$i\!\!-\!\!1$}

\htext(-13.3 13.75){$i$}

\htext(-15 12){$i\!\!+\!\!1$}\htext(-15 10){$i\!\!+\!\!2$}
\htext(-15 8){$i\!\!+\!\!3$}

\htext(-16 -1.5){$\overline{Y}_{\omega_{i_{k+1}}}$\,\,(or
$\overline{Y}_{\la_n}$)} \htext(-5
-1.5){$\overline{Y}_{\omega_{i_k}}$\,\,(or
$\overline{Y}_{\omega_{i_t}}$)}

\htext(-10 -4){$a=i+1$} \esegment
\end{texdraw}}

\vskip 5mm \begin{center} (i)
$\raisebox{-0.5\height}{\usebox{\tmpfiga}}$\qquad  (ii)
$\raisebox{-0.5\height}{\usebox{\tmpfigb}}$
\end{center}

%
\savebox{\tmpfiga}{\begin{texdraw} \fontsize{6}{6}\selectfont
\drawdim em \setunitscale 0.9 \nc{\dtri}{ \bsegment \lvec(-2 0)
\lvec(-2 2)\lvec(0 2)\lvec(0 0)\ifill f:0.7 \esegment }

\bsegment \move(-11 14)\dtri

\setgray 0 \move(0 0) \lvec(0 26) \move(-2 20)\rlvec(0 5)\rlvec(-2
0)\rlvec(0 -5)\move(-2 23)\rlvec(-4 0)\rlvec(0 -3) \move(-2
21)\rlvec(-4 0)

\move(-9 11)\rlvec(0 7)\rlvec(-2 0)\rlvec(0 -7)\move(-9
16)\rlvec(-4 0)\rlvec(0 -5) \move(-9 14)\rlvec(-4 0)\move(-9
12)\rlvec(-4 0)

\move(-16 0)\rlvec(0 6)\rlvec(-2 0)\rlvec(0 -6)\move(-16
4)\rlvec(-4 0)\rlvec(0 -4) \move(-16 2)\rlvec(-4 0)\move(-16
0)\rlvec(-4 0)

\move(-7.5 -4.3)\rlvec(0 27)

\move(-22 -3)\rlvec(22 0)\lvec(0 0)

\htext(-3 24){$i$}\htext(-5 22){$i\!\!-\!\!1$}\htext(-3
22){$i\!\!-\!\!1$}

\htext(-10 17){$i\!\!-\!\!1$}\htext(-10 15){$i$}\htext(-10
13){$i\!\!+\!\!1$}\htext(-12 15){$i$}\htext(-12 13){$i\!\!+\!\!1$}

\htext(-17 5){$i\!\!+\!\!1$}\htext(-17 3){$i$} \htext(-17
1){$i\!\!-\!\!1$}\htext(-19 3){$\circledast$} \htext(-19
1){$i\!\!-\!\!1$}

\move(-19 -1)\ravec(11.5 0)\move(-7.5 -1)\ravec(-11.5 0) \move(-19
-1.5)\rlvec(0 1)\htext(-19 -2){$p$-th}

\htext(-14 -4.4){$\overline{Y}_{\omega_{i_{k+1}}}$\,\,(or
$\overline{Y}_{\la_n}$)} \htext(-3
-4.4){$\overline{Y}_{\omega_{i_k}}$\,\,(or
$\overline{Y}_{\omega_{i_t}}$)}

\esegment
\end{texdraw}}
\savebox{\tmpfigb}{\begin{texdraw} \fontsize{6}{6}\selectfont
\drawdim em \setunitscale 0.9 \nc{\dtri}{ \bsegment \lvec(-2 0)
\lvec(-2 2)\lvec(0 2)\lvec(0 0)\ifill f:0.7 \esegment }

\bsegment \move(2 0)\dtri

\move(0 0) \rlvec(0 2)\rlvec(2 0)\rlvec(0 -2)\rlvec(-2 0)

\htext(1 1){$i$}

\esegment
\end{texdraw}}
\savebox{\tmpfigc}{\begin{texdraw} \fontsize{6}{6}\selectfont
\drawdim em \setunitscale 0.9 \nc{\dtri}{ \bsegment \lvec(-2 0)
\lvec(-2 2)\lvec(0 2)\lvec(0 0)\ifill f:0.7 \esegment }

\bsegment \move(9 15)\dtri

\setgray 0 \move(0 -2) \rlvec(0 5) \rlvec(-2 0)\rlvec(0 -5)\move(0
1)\rlvec(-4 0)\rlvec(0 -3) \move(0 -1)\rlvec(-4 0)\move(0
3)\rlvec(0 2)\rlvec(-2 0)\rlvec(0 -2)

\move(6 7) \rlvec(0 5) \rlvec(-2 0)\rlvec(0 -5)\move(6
10)\rlvec(-4 0)\rlvec(0 -3) \move(6 8)\rlvec(-4 0)

\move(11 14) \rlvec(0 5) \rlvec(-2 0)\rlvec(0 -5)\move(11
17)\rlvec(-4 0)\rlvec(0 -3) \move(11 15)\rlvec(-4 0)

\move(17 21) \rlvec(0 5) \rlvec(-2 0)\rlvec(0 -5)\move(17
24)\rlvec(-4 0)\rlvec(0 -3) \move(17 22)\rlvec(-4 0)

\move(-5 -4.7)\rlvec(26 0)\rlvec(0 30)

\move(12 -6)\rlvec(0 26)

\move(-3 -3.5)\rlvec(0 1)\move(-3 -3)\ravec(15 0)\move(12
-3)\ravec(-15 0)\htext(-3 -4){$s$-th}

\move(3 7)\rlvec(0 -1)\move(3 6.5)\ravec(9 0)\move(12
6.5)\ravec(-9 0)\htext(3 5.5){$r$-th}

\htext(-3 0){$j\!\!-\!\!1$}\htext(-1 0){$j\!\!-\!\!1$}\htext(-1
2){$j$}\htext(-1 4){$j\!\!+\!\!1$}

\htext(3 9){$j$}\htext(5 9){$j$}\htext(5 11){$j\!\!-\!\!1$}

\htext(8 16){$i$}\htext(10 16){$i$}\htext(10 18){$i\!\!-\!\!1$}

\htext(14 23){$i\!\!-\!\!1$}\htext(16 23){$i\!\!-\!\!1$}\htext(16
25){$i$}

\htext(4.5 -6.2){$\overline{Y}_{\omega_{i_{k+1}}}$\,\,(or
$\overline{Y}_{\la_n}$)} \htext(17
-6.2){$\overline{Y}_{\omega_{i_k}}$\,\,(or
$\overline{Y}_{\omega_{i_t}}$)}

\esegment
\end{texdraw}}
\savebox{\tmpfigd}{\begin{texdraw} \fontsize{6}{6}\selectfont
\drawdim em \setunitscale 0.9 \nc{\dtri}{ \bsegment \lvec(-2 0)
\lvec(-2 2)\lvec(0 2)\lvec(0 0)\ifill f:0.7 \esegment }

\bsegment

\setgray 0 \move(0 -2) \rlvec(0 5) \rlvec(-2 0)\rlvec(0 -5)\move(0
1)\rlvec(-4 0)\rlvec(0 -3) \move(0 -1)\rlvec(-4 0)\move(0
3)\rlvec(0 2)\rlvec(-2 0)\rlvec(0 -2)

\move(6 7) \rlvec(0 5) \rlvec(-2 0)\rlvec(0 -5)\move(6
10)\rlvec(-4 0)\rlvec(0 -3) \move(6 8)\rlvec(-4 0)

\move(11 14) \rlvec(0 5) \rlvec(-2 0)\rlvec(0 -5)\move(11
17)\rlvec(-4 0)\rlvec(0 -3) \move(11 15)\rlvec(-4 0)

\move(17 21) \rlvec(0 5) \rlvec(-2 0)\rlvec(0 -5)\move(17
24)\rlvec(-4 0)\rlvec(0 -3) \move(17 22)\rlvec(-4 0)

\move(-5 -4.7)\rlvec(26 0)\rlvec(0 30)

\move(12 -6)\rlvec(0 26)

\move(-3 -3.5)\rlvec(0 1)\move(-3 -3)\ravec(15 0)\move(12
-3)\ravec(-15 0)\htext(-3 -4){$s$-th}

\move(3 7)\rlvec(0 -1)\move(3 6.5)\ravec(9 0)\move(12
6.5)\ravec(-9 0)\htext(3 5.5){$r$-th}

\move(8 14)\rlvec(0 -1)\move(8 13.5)\ravec(4 0)\move(12
13.5)\ravec(-4 0)\htext(8 12.5){$u$-th}

\move(14 21)\rlvec(0 -1)\move(14 20.5)\ravec(7 0)\move(21
20.5)\ravec(-7 0)\htext(14 19.5){$t$-th}

\htext(-3 0){$j\!\!-\!\!1$}\htext(-1 0){$j\!\!-\!\!1$}\htext(-1
2){$j$}\htext(-1 4){$j\!\!+\!\!1$}

\htext(3 9){$j$}\htext(5 9){$j$}\htext(5 11){$j\!\!-\!\!1$}

\htext(8 16){$k$}\htext(10 16){$k$}\htext(10 18){$k\!\!-\!\!1$}

\htext(14 23){$k\!\!-\!\!1$}\htext(16 23){$k\!\!-\!\!1$}\htext(16
25){$k$}

\htext(4.5 -6.2){$\overline{Y}_{\omega_{i_{k+1}}}$\,\,(or
$\overline{Y}_{\la_n}$)} \htext(17
-6.2){$\overline{Y}_{\omega_{i_k}}$\,\,(or
$\overline{Y}_{\omega_{i_t}}$)}

\esegment
\end{texdraw}}
\vskip 3mm In the case of (i), if $\fit Y$ violates the condition
{\bf (Y4)}, then $\fit Y$ must have the form
$$\raisebox{-0.3\height}{\usebox{\tmpfiga}}\qquad\text{or}\qquad
\raisebox{-0.3\height}{\usebox{\tmpfigc}}$$  In the first case,
$\fit$ would have acted on $\circledast$, not on
$\raisebox{-0.3\height}{\usebox{\tmpfigb}}$\,, which is a
contradiction. In the second case, since $Y$ satiesfies {\bf (Y1)}
and {\bf (Y2)}, we can observe that $Y$ must have the form
$$\raisebox{-0.4\height}{\usebox{\tmpfigd}}$$ where $s\leq t<p$,\,\, $q<u<r$ and $i<k<j$.
That is, $Y$ also has the configuration {\bf (C1)} and does not
satisfy {\bf (Y4)}, which is a contradiction.

\vskip 3mm In the case of (ii), we note that the top of $(q-1)$-th
column in $\overline{Y}_{\omega_{i_{k+1}}}$ (or
$\overline{Y}_{\la_n})$ must be an $i$-block by the tensor product
rule for the Kashiwara operators. Then we know that $Y$ also has a
configuration {\bf (C1)}.

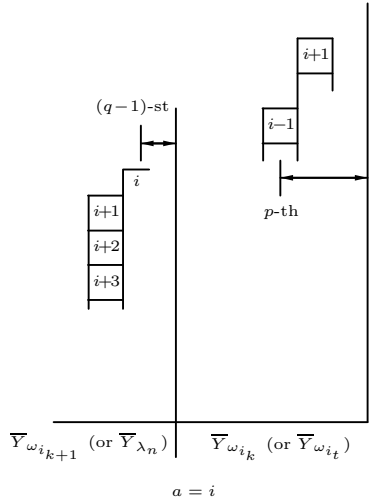
\savebox{\tmpfiga}{\begin{texdraw} \fontsize{6}{6}\selectfont
\drawdim em \setunitscale 0.9 \nc{\dtri}{ \bsegment \lvec(-2 0)
\lvec(-2 2)\lvec(0 2)\lvec(0 0)\ifill f:0.7 \esegment }

\bsegment

\setgray 0 \move(0 0) \lvec(0 24) \move(-2 22)\rlvec(-2 0)\move(-2
20)\rlvec(-2 0)\move(-4 18)\rlvec(-2 0)\move(-4 16)\rlvec(-2
0)\move(-6 15)\rlvec(0 3)\move(-4 15)\rlvec(0 7) \move(-2
19)\rlvec(0 3)

\move(-14 13)\rlvec(-2 0)\move(-14 11)\rlvec(-2 0) \move(-14
9)\rlvec(-2 0)\move(-14 7)\rlvec(-2 0) \move(-14 6.5)\rlvec(0
6.5)\move(-16 6.5)\rlvec(0 6.5)\move(0 0)\rlvec(-18 0)\move(-14
13)\rlvec(0 1.5)\rlvec(1.5 0)

\move(-11 -2)\rlvec(0 20)

\move(-5 15)\rlvec(0 -2)\move(-5 14)\ravec(5 0)\move(0
14)\ravec(-5 0)\htext(-5 12){$p$-th}

\move(-13 17)\rlvec(0 -2)\move(-13 16)\ravec(2 0)\move(-11
16)\ravec(-2 0)\htext(-13.5 18){$(q\!-\!1)$-st}

\htext(-3 21){$i\!\!+\!\!1$}\htext(-5 17){$i\!\!-\!\!1$}

\htext(-13.3 13.75){$i$}

\htext(-15 12){$i\!\!+\!\!1$}\htext(-15 10){$i\!\!+\!\!2$}
\htext(-15 8){$i\!\!+\!\!3$}

\htext(-16 -1.5){$\overline{Y}_{\omega_{i_{k+1}}}$\,\,(or
$\overline{Y}_{\la_n}$)} \htext(-5
-1.5){$\overline{Y}_{\omega_{i_k}}$\,\,(or
$\overline{Y}_{\omega_{i_t}}$)}

\htext(-10 -4){$a=i$} \esegment
\end{texdraw}}
$$\raisebox{-0.3\height}{\usebox{\tmpfiga}}$$
But, since $Y$ satisfies the condition {\bf (Y4)}, adding an
$i$-block to $\overline {Y}_{\omega_{i_k}}$ does not creat a Young
wall which violates the condition {\bf (Y4)}; i.e., the second
case does not occur.

\vskip 2mm  Similarly, if $Y\in Y(\la)$, then we can show that
$\eit Y\in Y(\la)\cup \{0\}$.

%
\savebox{\tmpfiga}{\begin{texdraw} \fontsize{6}{6}\selectfont
\drawdim em \setunitscale 0.9 \nc{\dtri}{ \bsegment \lvec(-2 0)
\lvec(-2 2)\lvec(0 2)\lvec(0 0)\ifill f:0.7 \esegment }

\bsegment

\move(0 0) \rlvec(0 2)\rlvec(2 0)\rlvec(0 -2)\rlvec(-2 0)

\htext(1 1){$i$}

\esegment
\end{texdraw}}
\savebox{\tmpfigb}{\begin{texdraw} \fontsize{6}{6}\selectfont
\drawdim em \setunitscale 0.9 \nc{\dtri}{ \bsegment \lvec(-2 0)
\lvec(-2 2)\lvec(0 2)\lvec(0 0)\ifill f:0.7 \esegment }

\bsegment \move(-8 18)\dtri

\move(0 0)\lvec(0 27) \move(0 26)\rlvec(-2 0)\rlvec(0 -2)
\rlvec(-2 0)\rlvec(0 -2)\rlvec(-2 0)\rlvec(0 -2)\rlvec(-4
0)\rlvec(0 -14)

\move(0 14)\rlvec(-2 0)\rlvec(0 -2)\rlvec(-2 0)\rlvec(0 -2)
\rlvec(-2 0)\rlvec(0 -2)\rlvec(-2 0)\rlvec(0 -2)\rlvec(-2
0)\rlvec(0 -2)\rlvec(-2 0)\rlvec(0 -2)\rlvec(-2 0)\rlvec(0 -2)

\move(-18 0)\lvec(0 0)

\move(-8 0)\rlvec(0 20)\move(-10 0)\rlvec(0 4)

\move(-13 11)\avec(-9.2 12)

\move(-18 16)\lvec(3 16)

\htext(-9 19){$i$}\htext(-11 5){$i$}

\htext(-14 11){$Y_s$}

\htext(4.5 16){n-row}

\htext(-8 -1.5){$\overline{Y}_{\omega_{i_k}}$}

\htext(-5 5){${\overline H}_{\omega_{i_k}}$} \esegment
\end{texdraw}}

\vskip 3mm Now, it remains to prove the statement (2). Suppose
$Y\in Y(\la)$ and $\eit Y=0$ for all $i=1,\cdots,n$. If $Y\neq
H_{\la}$, then there is a column in $Y$ which is higher than
$H_{\la}$. Consider the left-most column $Y_s$ among them, which
would belong to $\overline{Y}_{\omega_{i_k}}$ or
$\overline{Y}_{\la_n}$. Let
$\raisebox{-0.2\height}{\usebox{\tmpfiga}}$ be the block lying in
the top of the column $Y_s$. If there is an $i$-admissible slot to
the left of $\raisebox{-0.2\height}{\usebox{\tmpfiga}}$\,, then
$Y$ has the form
$$\raisebox{-0.3\height}{\usebox{\tmpfigb}}.$$
However, in this case, $Y_{\omega_{i_k}}^+\nsubseteq
|Y_{\omega_{i_k}}^-|$, which violates the condition {\bf (Y1)}.
Hence, there is no admissible $i$-slot to the left of $Y_s$, which
implies $\eit
Y=Y\nearrow\raisebox{-0.3\height}{\usebox{\tmpfiga}}\neq 0$, a
contradiction. Therefore, $Y$ must be equal to $H_{\la}$. $\quad
\square$

\vskip 15mm


\providecommand{\bysame}{\leavevmode\hbox
to5em{\hrulefill}\thinspace}

\end{document}